\documentclass[11pt]{amsart}

\usepackage{amsmath}
\usepackage{amssymb}
\usepackage{amsthm}
\usepackage{enumitem}
\usepackage{microtype}
\usepackage{array}
\usepackage{dcolumn}
\usepackage[margin=1.0in]{geometry}
\usepackage{float}
\usepackage{graphicx}
\usepackage{tikz-cd}
\usepackage[style=alphabetic]{biblatex}
\usepackage[hidelinks]{hyperref}
\usepackage{mathtools}
\usepackage{todonotes}
\usepackage{mathrsfs}
\usepackage{lipsum}

\renewcommand \v{\vert}

\renewcommand \Im{\mathrm{Im}}

\newcommand \R{\mathbb{R}}

\newcommand \Z{\mathbb{Z}}
\newcommand \C{\mathbb{C}}

\newcommand \CP{\mathbb{CP}}

\newcommand \ep{\varepsilon}
\newcommand \la{\langle}
\newcommand \ra{\rangle}
\newcommand \del{\partial}

\DeclareMathOperator \id{id}
\DeclareMathOperator \rel{rel}

\DeclareMathOperator \Ker{Ker}
\DeclareMathOperator \Int{Int}

\DeclareMathOperator \Aut{Aut}
\DeclareMathOperator \Out{Out}
\DeclareMathOperator \Inn{Inn}

\DeclareMathOperator \Homeo{Homeo}
\DeclareMathOperator \Diff{Diff}

\DeclareMathOperator \LDiff{LDiff}
\DeclareMathOperator \Mod{Mod}

\DeclareMathOperator \GL{GL}

\DeclareMathOperator \SO{SO}

\renewcommand \L{\mathcal{L}}

\newcommand \F{\mathcal{F}}
\renewcommand \P{\mathcal{P}}
\DeclareMathOperator \LMod{LMod}

\DeclareMathOperator \Conf{Conf}

\DeclareMathOperator \UConf{UConf}
\DeclareMathOperator \Emb{Emb}
\DeclareMathOperator \UEmb{UEmb}
\DeclareMathOperator \Fr{Fr}
\DeclareMathOperator \UFr{UFr}

\DeclareMathOperator \Fix{Fix}
\DeclareMathOperator \Stab{Stab}

\DeclareMathOperator \Sub{Sub}

\DeclareMathOperator \SymAut{SymAut}
\DeclareMathOperator \SymOut{SymOut}
\DeclareMathOperator \PSymOut{PSymOut}
\newcommand \calQ{\mathcal{Q}}

\DeclareMathOperator \DS{DS}
\newcommand \calDS{\mathcal{DS}}

\DeclareMathOperator \orb{orb}
\DeclareMathOperator \lift{lift}

\newcommand \calC{\mathcal{C}}

\DeclareMathOperator \ev{ev}
\DeclareMathOperator \Iso{Iso}
\newcommand \calO{\mathcal{O}}
\newcommand \calG{\mathcal{G}}
\newcommand \scrG{\mathscr{G}}

\theoremstyle{definition}
\newtheorem{definition}{Definition}[section]

\theoremstyle{plain}

\newtheorem{lemma}[definition]{Lemma}
\newtheorem{proposition}[definition]{Proposition}

\newtheorem{question}[definition]{Question}

\newtheorem{problem}[definition]{Problem}
\newtheorem{mainTheorem}{Theorem}

\newtheorem{mainCorollary}[mainTheorem]{Corollary}
\theoremstyle{remark}
\newtheorem{remark}[definition]{Remark}

\numberwithin{equation}{section}

\title{Isotopy Versus Equivariant Isotopy in Dimensions Three and Higher}
\author{Trent Lucas}

\begin{document}

\begin{abstract}
    Given a finite group action on a smooth manifold, we study the following question: if two equivariant diffeomorphisms are isotopic, must they be equivariantly isotopic?  
    Birman--Hilden and Maclachlan--Harvey proved the answer is ``yes'' for most surfaces.
    By contrast, we give a general criterion in higher dimensions under which there are many equivariant diffeomorphisms which are isotopic but not equivariantly isotopic.
    Examples satisfying this criterion include branched covers of split links and ``stabilized'' branched covers.
    We prove the result by constructing an invariant valued in the homology of a certain infinite cover of the manifold.
    We give applications to outer automorphism groups of free products and to group actions on manifolds which fiber over the circle.
\end{abstract}

\maketitle


\section{Introduction}

\subsection{The main question and result}
Let $M$ be a closed oriented smooth manifold, and $G$ a finite group of orientation-preserving diffeomorphisms of $M$.
We call a diffeomorphism of $M$ \emph{equivariant} if it commutes with $G$.
Given two equivariant diffeomorphisms of $M$, we can ask whether they are isotopic, or we can ask whether they are \emph{equivariantly} isotopic, i.e.\ isotopic through equivariant diffeomorphisms.
Our goal in this paper is to measure the difference between these two equivalence relations.

More formally, let $\Diff(M)$ denote the group of orientation-preserving diffeomorphisms of $M$, and let $\Diff(M)^G$ denote the subgroup of equivariant diffeomorphisms.
Let $\Mod(M)$ denote $\pi_0(\Diff(M))$, i.e.\ the group of isotopy classes of orientation-preserving diffeomorphisms, also called the \emph{mapping class group}
of $M$.
Let $\Gamma_G(M)$ denote $\pi_0(\Diff(M)^G)$, i.e.\ the group of equivariant isotopy classes of equivariant diffeomorphisms.
Then the inclusion $\Diff(M)^G \hookrightarrow \Diff(M)$ induces a map $\P_G:\Gamma_G(M) \rightarrow \Mod(M)$.
We then ask:

\begin{question}\label{ques:main-question}
    Is the map $\P_G:\Gamma_G(M) \rightarrow \Mod(M)$ injective?  If not, how large is the kernel?
\end{question}


A remarkable theorem of Birman--Hilden \cite{birman-hilden} and Maclachlan--Harvey \cite{maclachlan-harvey} says that if $M$ is a closed oriented surface of genus $g \geq 2$, then $\P_G$ is injective.
By contrast, we showed in previous work \cite{bh-3manifolds} that for many finite group actions on closed oriented $3$-manifolds, the kernel $\Ker(\P_G)$ is infinite, answering a question of Margalit--Winarski \cite[Ques~11.4]{margalit-winarski}.
Our main theorem strengthens this result and extends it to all dimensions $n \geq 3$.

\begin{mainTheorem}\label{mainthm:ker-not-fg}
    Let $M$ be a closed oriented connected smooth manifold of dimension $n \geq 3$, and let $G$ be a finite group of orientation-preserving diffeomorphisms.
    Let $M^\circ \subseteq M$ denote the set of points whose $G$-stabilizer is trivial.
    Assume that:
    \begin{itemize}
        \item For some $k > 0$, the quotient manifold $M^\circ/G$ is diffeomorphic to a $k$-fold connected sum $Q_0 \# \cdots \# Q_{k-1}$ where $\pi_1(Q_i)$ is nontrivial for each $i$.
        \item For some $g_0 \in G$, the set of fixed points $\Fix(g_0) \subseteq M$ has a nonempty codimension $2$ orientable component $B$ whose homology class in $H_{n-2}(M;\Z/\ell\Z)$ is trivial for some $\ell > 1$.
    \end{itemize}
    Then if $k \geq 3$:
    \begin{enumerate}[label=(\roman*)]
        \item the kernel of the map $\P_G:\Gamma_G(M) \rightarrow \Mod(M)$ contains a virtual free product, and
        \item the kernel $\Ker(\P_G)$ is not finitely generated.
    \end{enumerate}
\end{mainTheorem}

In part (i), a \emph{virtual free product} is a group which is virtually isomorphic to a free product of two nontrivial groups; here we say that two groups $G_1$ and $G_2$ are \emph{virtually isomorphic} if there are finite index subgroups $H_i \leq G_i$ and finite normal subgroups $F_i \leq H_i$ such that $H_1/F_1 \cong H_2/F_2$.
Note that Theorem \ref{mainthm:ker-not-fg} reduces to the case $k=3$, since for higher $k$ the quotient can be written as $Q_0 \# Q_1 \# Q$ where $Q = Q_2 \# \cdots \# Q_{k-1}$.



The surprising aspect of Theorem \ref{mainthm:ker-not-fg} is not just that $\Ker(\P_G)$ is nontrivial, but also that $\Ker(\P_G)$ is so large at this level of generality.
Indeed, the bulk of this paper is devoted to proving part (ii).
The proof requires us to find an invariant that can distinguish many equivariant diffeomorphisms which are all isotopically trivial; we accomplish this with an ``obstruction map'' valued in the homology of a certain infinite cover of $M$.
This invariant is simple to describe, but finding even one example where this invariant is nontrivial requires a careful and explicit computation.
The proof of part (i) is simpler, and serves as a warmup to the proof of part (ii).

We give more details of the proof in Section \ref{subsec:pf-sketch} below.
Before that, we will describe examples of group actions satisfying the hypotheses of Theorem \ref{mainthm:ker-not-fg}, and give two applications of Theorem \ref{mainthm:ker-not-fg}.

\subsection{Example: Branched covers of links}

A simple way to build examples for Theorem \ref{mainthm:ker-not-fg} in all dimensions $n \geq 3$ is to take branched covers of links in $S^n$.
We say that a surjective map of smooth manifolds $p:M \rightarrow N$ is a \emph{finite regular branched cover} if there is a compact codimension $2$ submanifold $C \subseteq N$ such that
\begin{itemize}
    \item the restriction of $p$ to $M \setminus p^{-1}(C) \rightarrow N \setminus C$ is a finite-sheeted regular cover, and
    \item for each $x \in p^{-1}(C)$, we can choose coordinates around $x$ and $p(x)$ valued in $\C \times \R^{n-2}$ in which $p$ is of the form $(z,x) \mapsto (z^m, x)$ for some integer $m \geq 2$.
\end{itemize}
A \emph{link} is the image of a smooth embedding $S^{n-2} \sqcup \cdots \sqcup S^{n-2} \hookrightarrow S^n$.
Given a link $L$ in $S^n$ and an integer $d \geq 2$, there is a closed oriented smooth manifold $M_{L,d}$ and a finite regular branched cover $p_{L,d}:M_{L,d} \rightarrow S^n$ with branch set $L$ and deck group $G = \Z/d\Z$ (see e.g.\ \cite[\S 27]{ranicki}, which constructs $M_{L,d}$ in the topological category but can be easily adapted to the smooth category).


We say a link $L$ is \emph{$k$-split} if there are $k$ disjoint closed $n$-disks $D_0, \ldots, D_{k-1} \subseteq S^n$ such that $D_i \cap L$ is a nonempty union of components of $L$.
Given a $k$-split link $L$, we claim that the action of $G$ on $M = M_{L,d}$ satisfies the hypotheses of Theorem \ref{mainthm:ker-not-fg}.
Indeed, let $Q_i$ be the $n$-manifold obtained by capping the boundary of $D_i \setminus L$ with an $n$-disk.
Then $S^n \setminus L$ is diffeomorphic to the connected sum $Q_0 \# \cdots \# Q_{k-1}$.
Since $H_1(Q_i)$ is nontrivial for each $i$ by Alexander duality, the first hypothesis of Theorem \ref{mainthm:ker-not-fg} is satisfied.
For the second hypothesis, we can take $B$ to be the preimage of any component of $L$; then $B$ will be homologically trivial in $M$ since its meridian will be nontrivial in $H_1(M \setminus B)$ (see Lemma \ref{lem:meridian-vs-homology-class}).
Thus Theorem \ref{mainthm:ker-not-fg} implies:

\begin{mainCorollary}\label{maincor:covers-of-links}
    Let $L$ be a $k$-split link in $S^n$ for $n \geq 3$, and let $M = M_{L,d}$ be the associated branched cover with deck group $G = \Z/d\Z$.
    If $k \geq 3$, then $\Ker(\P_G)$ is not finitely generated.
\end{mainCorollary}

The simplest case of Corollary \ref{maincor:covers-of-links} is when $d = 2$ and $L$ is the $k$-component unlink in $S^n$ (i.e.\ each component is isotopic to the standard embedding $S^{n-2} \hookrightarrow S^n$).
In this case, $M \cong (S^1 \times S^{n-1})^{\# k-1}$ and the action of $G = \Z/2\Z$ on $M$ is a generalization of the hyperelliptic involution on a genus $k-1$ surface (which was the example first studied by Birman and Hilden).
In \cite{bh-3manifolds}, we showed that $\Ker(\P_G)$ is not finitely generated for the double cover of the $3$-component unlink in $S^3$ using an algebraic proof much different from that of Theorem \ref{mainthm:ker-not-fg}.

\subsection{Example: Stabilized actions and branched covers}

Another family of examples for Theorem \ref{mainthm:ker-not-fg} come from ``stabilizing'' group actions and branched covers, as follows.

Suppose $P$ is a closed smooth $n$-manifold and $G$ is a finite group of diffeomorphisms.
For $1 \leq i \leq k$, let $Q_i$ be any closed oriented connected smooth $n$-manifold with $\pi_1(Q_i)$ nontrivial, and let $P_i$ be the disjoint union of $d$ copies of $Q_i$, where $d = \v G \v$.
We then define a closed smooth $n$-manifold $M$ by taking an \emph{equivariant connected sum} of $P$ with each $P_i$, meaning that we remove an open disk from each component of $P_i$, remove a free $G$-orbit of open disks from $P$, and glue together the resulting boundary components (some care must be taken to preserve smoothness and orientability, as in Section \ref{subsec:summand-slides}).
The action of $G$ on $P$ naturally extends to an action on $M$ by permuting the components of each $P_i$.
We call $M$ a \emph{k-fold stabilization} of the action of $G$ on $P$.

If $q:P \rightarrow Q$ is a finite regular branched cover with deck group $G$, and $M$ is a $k$-fold stabilization of the action of $G$ on $P$, then we get an induced branched cover $p:M \rightarrow N$ where $N = Q \# Q_1 \cdots \# Q_k$.
We call the branched cover $p$ a \emph{$k$-fold stabilization} of the branched cover $q$.
Theorem \ref{mainthm:ker-not-fg} implies:

\begin{mainCorollary}\label{maincor:stabilized-branched-cover}
    Let $P$ and $Q$ be closed oriented connected smooth manifolds of dimension $n \geq 3$, and let $q:P \rightarrow Q$ be a finite regular branched cover with branch set $C \subseteq Q$ and deck group $G$.
    Assume that that $p^{-1}(C)$ represents a trivial element of $H_{n-2}(P;\Z/\ell\Z)$ for some $\ell > 1$.
    Let $p:M \rightarrow N$ be a $k$-fold stabilization of $q$.
    Then if $k \geq 2$, the kernel of $\P_G:\Gamma_G(M) \rightarrow \Mod(M)$ is not finitely generated.
\end{mainCorollary}

To deduce Corollary \ref{maincor:stabilized-branched-cover} from Theorem \ref{mainthm:ker-not-fg}, we use the fact that $Q \setminus C$ is not simply connected.
This is because the unbranched cover $P \setminus q^{-1}(C) \rightarrow Q \setminus C$ is classified by a surjection $\pi_1(Q \setminus C) \rightarrow G$.

For a concrete family of examples in even dimensions, we can take $P$ to be a degree $d$ cyclic cover of $Q = \CP^{n/2}$ branched over a smooth degree $m$ complex hypersurface $C$, where $d$ is a proper divisor of $m$.
To check that the assumption of Corollary \ref{maincor:stabilized-branched-cover} holds, recall first that $\pi_1(\CP^{n/2} \setminus $C$) \cong \Z/m\Z$ and is generated by a meridian of $C$ (see \cite{cogolludo}).
Then $\pi_1(P \setminus q^{-1}(C)) \cong \Z/\ell \Z$, where $\ell = m/d$.
Thus $H_1(P \setminus q^{-1}(C); \Z/\ell\Z)$ is isomorphic to $\Z/\ell\Z$ and is generated by a meridian of $q^{-1}(C)$.
This implies that $q^{-1}(C)$ must represent a trivial homology class in $H_{n-2}(P; \Z/\ell\Z)$ (this can be seen by adapting the proof of Lemma \ref{lem:meridian-vs-homology-class}).

\subsection{Application: Symmetric automorphism groups of free products}\label{subsec:appliction-to-symmetric-auts}
We can apply Theorem \ref{mainthm:ker-not-fg} to study certain natural maps arising in group theory.

Let $L = L_1 * \cdots * L_k$ be any free product of groups.
We define $\SymAut(L)$ to be the group of \emph{symmetric} automorphisms of $L$, i.e.\ the group of automorphisms of $L$ sending each $L_i$ to a conjugate of some $L_j$, and we define $\SymOut(L)$ to be the image of $\SymAut(L)$ in $\Out(L)$.
In the case that each $L_i$ is freely indecomposable and not infinite cyclic, $\SymAut(L)$ is in fact the full automorphism group $\Aut(L)$ by the Kurosh subgroup theorem.
In the case that $L$ is the free group $F_k$, the group $\SymAut(F_k)$ arises naturally in the study of the braid group $B_k$, as $B_k$ acts on the fundamental group of the punctured disk by symmetric automorphisms.

Given integers $k,d > 1$, we define the free product
\begin{equation*}
    H_{k,d}
    \coloneqq
    \underbrace{(\Z/d\Z) * \cdots * (\Z/d\Z)}_{\text{$k$ times}}.
\end{equation*}
Then the natural surjection $F_k \rightarrow H_{k,d}$ induces maps
\begin{align*}
    &\widehat{\calQ}_{k,d}: \SymAut(F_k) \rightarrow \SymAut(H_{k,d}) \\
    &\calQ_{k,d}: \SymOut(F_k) \rightarrow \SymOut(H_{k,d}).
\end{align*}
As a surprising application of their work on Question \ref{ques:main-question}, Birman--Hilden \cite{birman-hilden} showed that the map $\widehat{\calQ}_{k,d}$ restricts to an embedding $B_k \hookrightarrow \SymAut(H_{k,d})$, answering a question of Magnus
\footnote{
    Magnus later applied this result in \cite{magnus}.
}
(see also \cite[\S 7]{margalit-winarski}).
Apart from this, it appears that the kernels of the maps $\widehat{\calQ}_{k,d}$ and $\calQ_{k,d}$ are not well-studied.
These kernels contain some obvious elements (e.g.\ conjugating one generator by the $d$th power of another), but it is not clear which elements generate these kernels or what finiteness properties these kernels possess.
Using Theorem \ref{mainthm:ker-not-fg}, we show:

\begin{mainTheorem}\label{mainthm:symmetric-aut-kernel}
    Let $k \geq 3$ and $d \geq 2$.
    Then the groups $\Ker(\widehat{\calQ}_{k,d})$ and $\Ker({\calQ}_{k,d})$ are not finitely generated.
\end{mainTheorem}

Theorem \ref{mainthm:symmetric-aut-kernel} is sharp.
For $d =1$, the maps $\widehat{\calQ}_{k,d}$ and $\calQ_{k,d}$ are trivial, and the groups $\SymAut(F_k)$ and $\SymOut(F_k)$ are finitely generated (see e.g.\ \cite{mccool}).
For $k = 1$, we have that $\SymAut(\Z)$ and $\SymOut(\Z)$ are isomorphic to $\Z/2\Z$.  For $k=2$, the group $\Inn(F_2)$ has finite index in $\SymAut(F_2)$ (again see e.g.\ \cite{mccool}), which implies that $\Ker(\widehat{\calQ}_{k,d})$ is not finitely generated but $\Ker(\calQ_{k,d})$ is finite.

Following our previous work in \cite{bh-3manifolds}, we give a topological proof of Theorem \ref{mainthm:symmetric-aut-kernel}.
The key idea is to apply Theorem \ref{mainthm:ker-not-fg} to finite covers $p_{k,d}:M_{k,d} \rightarrow S^3$ branched over the $k$-component unlink $C_k$ with deck group $G = \Z/d\Z$.
Let $\Diff(S^3, C_k^{\pm})$ denote the group orientation-preserving diffeomorphisms of $S^3$ preserving $C_k$ setwise (not necessarily preserving the orientation of $C_k$), and let $\Mod(S^3,C_k^{\pm}) = \pi_0(\Diff(S^3, C_k^{\pm}))$.
Then the relationship between the cover $p = p_{k,d}$ and the group $\SymOut(F_k)$ comes from a theorem of Wattenberg \cite{wattenberg}, which says that the action of $\Mod(S^3, C_k^{\pm})$ on $\pi_1(S^3 \setminus C_k)$ induces an isomorphism
\begin{equation*}
    \Mod(S^3, C_k^{\pm}) \cong \SymOut(F_k)
\end{equation*}
(the analogous result holds in the topological category by work of Goldsmith \cite{goldsmith}, and both results build off the work of Dahm \cite{dahm}).
The cover $p$ induces an orbifold structure $\mathcal{O}_p$ on $S^3$, and $\Mod(S^3,C_k^{\pm})$ also acts on the orbifold fundamental group $\pi_1^{\orb}(\mathcal{O}_p) \cong H_{k,d}$.
This yields a relationship between $\Ker(\calQ_{k,d})$ and $\Ker(\P_G)$ which we can use to prove Theorem \ref{mainthm:symmetric-aut-kernel}.

In \cite{bh-3manifolds}, we computed a finite \emph{normal} generating set of $\Ker(\calQ_{k,2})$ using results of McCullough--Miller \cite{mccullough-miller} on the groups $\SymOut(F_k)$ and $\SymOut(H_{k,d})$. 
It would be interesting to determine whether the groups $\Ker(\calQ_{k,d})$ have a finite normal generating set for all $k$ and $d$ using our topological framework.

\subsection{Application: Trivial bundles with a nontrivial fiberwise action}

As explained in \cite{bh-3manifolds}, Theorem \ref{mainthm:ker-not-fg} provides many examples of fiber bundles $E \rightarrow S^1$ equipped with a fiberwise group action such that $E$ is smoothly trivial, but not equivariantly trivial.

Let $M$ be a closed oriented $n$-manifold and $G$ a finite group of diffeomorphisms satisfying the hypotheses of Theorem \ref{mainthm:ker-not-fg}.
Then the non-injectivity of the map $\P_G$ implies that the map
\begin{equation*}
    \widehat{\P}_G:\pi_1(B\Diff(M)^G) \rightarrow \pi_1(B\Diff(M))
\end{equation*}
is not injective, where $B$ denotes the classifying space.
Conjugacy classes in $\pi_1(B\Diff(M))$ correspond to isomorphism classes of $M$-bundles over $S^1$.
Conjugacy classes in $\pi_1(B\Diff(M)^G)$ correspond to isomorphism classes of $M$-bundles over $S^1$ with structure group $\Diff(M)^G$; these are bundles for which the action of $G$ on $M$ extends to a fiberwise action on the total space.
Thus, nontrivial elements in $\Ker(\widehat{P}_G)$ yield bundles which are smoothly trivial, but not equivariantly trivial, i.e.\ not equivariantly isomorphic to the product $S^1 \times M$ equipped with the standard action of $G$ in each fiber. 

The same reasoning tells us that elements in the kernel of the map $\pi_k(\Diff(M)^G) \rightarrow \pi_k(\Diff(M))$ correspond to $M$-bundles over $S^{k+1}$ which are smoothly but not equivariantly trivial.
Our construction of elements in the kernel for $k=0$ does not obviously extend to higher $k$, so we naturally pose:

\begin{problem}\label{prob:higher_BH}
    Let $G$ be a finite group of diffeomorphisms of a closed oriented smooth $n$-manifold $M$.
    For $k \geq 1$, construct elements in the kernel of the map $\pi_k(\Diff(M)^G) \rightarrow \pi_k(\Diff(M))$ (or show this map is injective).
\end{problem}

\subsection{The topological setting}\label{subsec:topological-setting}
Theorem \ref{mainthm:ker-not-fg} is stated and proved in the smooth category.
We can adapt the proof to the topological setting under the assumption that the action of $G$ is locally linear.

That is, suppose $M$ is a closed oriented topological manifold and $G$ is a finite group of orientation-preserving homeomorphisms.
We say that the action of $G$ is \emph{locally linear} if for every $x \in M$ with stabilizer $G_x \leq G$, there is 
\begin{itemize}
    \item a $G_x$-invariant open neighborhood $U \subseteq M$ of $x$,
    \item a representation $\rho_x:G_x \rightarrow \SO(n)$, and
    \item a $G_x$-equivariant homeomorphism $U \cong \R^n$ taking $x$ to the origin, where $G$ acts on $\R^n$ via the representation $\rho_x$.
\end{itemize} 

Observe that if the action of $G$ is locally linear, then for any $g \in G$, the fixed set $\Fix(g)$ is a locally flat submanifold of $M$.
Note also that if $M$ is smooth and $G$ acts by diffeomorphisms, then the action is locally linear; this follows by choosing a $G$-invariant Riemannian metric and using the exponential map (see Lemma \ref{lem:local-model-of-action}).

We can now state the topological version of Theorem \ref{mainthm:ker-not-fg} as follows.

\begin{mainTheorem}\label{mainthm:top-version}
    Let $M$ be a closed oriented connected topological manifold of dimension $n \geq 3$, and let $G$ be a finite group of orientation-preserving homeomorphisms whose action is locally linear. 
    Assume $M$ and $G$ satisfy the assumptions of Theorem \ref{mainthm:ker-not-fg} (replacing ``diffeomorphic'' with ``homeomorphic'').
    Then if $k \geq 3$, the kernel of the natural map $\pi_0(\Homeo(M)^G) \rightarrow \pi_0(\Homeo(M))$ is not finitely generated.
\end{mainTheorem}

We expect that conclusion (i) of Theorem \ref{mainthm:ker-not-fg} can also be proven in the topological setting; since part (i) mainly serves as a warmup to part (ii), we leave this adaptation to the reader.

\subsection{Proof idea}\label{subsec:pf-sketch}

As mentioned above, the surprising aspect of Theorem \ref{mainthm:ker-not-fg} is that $\Ker(\P_G)$ is so large.
In fact, if one only wishes to show that $\Ker(\P_G)$ is nontrivial, there is a simple strategy to do so.
Namely, for any $g_0 \in G$, any equivariant diffeomorphism must preserve the fixed set $\Fix(g_0) \subseteq M$ setwise, and an element of $\Ker(\P_G)$ may restrict to a nontrivial mapping class of $\Fix(g_0)$ (for example, this occurs for the hyperelliptic involution on the $2$-torus \cite[\S 2]{margalit-winarski}).
However, if the mapping class group of $\Fix(g_0)$ is finite (such as for branched covers of links in $S^n$ for $n \leq 5$), then this obstruction only gives a finite lower bound on the size of $\Ker(\P_G)$.
The elements of $\Ker(\P_G)$ that we construct in the proof of Theorem \ref{mainthm:ker-not-fg} will in fact act trivially on the fixed points of any $g_0 \in G$.

The elements of $\Ker(\P_G)$ that we use to prove Theorem \ref{mainthm:ker-not-fg} are lifts of ``summand slides,'' defined as follows.
Choose a splitting of $M^\circ/G$ as a connected sum $N_1 \# N_2$, and let $\Sigma \subseteq M^\circ/G$ be a separating $(n-1)$-sphere realizing this decomposition.
Let $\delta$ be a loop in $N_1$ based at $\Sigma$.
Informally, a summand slide is a diffeomorphism of $M^\circ/G$ obtained by dragging $\Sigma$ along $\delta$, similar to a point-pushing diffeomorphism of a surface.

Part (i) of Theorem \ref{mainthm:ker-not-fg} has a fairly direct proof.
We find loops $\delta_1$ and $\delta_2$ in $M^\circ/G$ which are nontrivial in $\pi_1(M^\circ/G)$, but lift to trivial elements of $\pi_1(M)$.
We then let $f_i$ be a summand slide of $M^\circ/G$ along $\delta_i$, and let $\widetilde{f}_i$ be a lift of $f_i$ to an equivariant diffeomorphism of $M$.
A key property of summand slides is that up to finite order, they are determined by their action on $\pi_1$.
We show that $f_1$ and $f_2$ generate a free product in $\Out(\pi_1(M^\circ/G))$, which implies that $\widetilde{f}_1$ and $\widetilde{f}_2$ generate a virtual free product in $\Gamma_G(M)$, but $\widetilde{f}_1$ and $\widetilde{f}_2$ act trivially on $\pi_1(M)$, which implies that they generate a finite subgroup of $\Mod(M)$.

Part (ii) is more difficult, and its proof comprises the majority of this paper.
Let $\Mod(M,B) = \pi_0(\Diff(M,B))$, where $\Diff(M,B)$ is the group of orientation-preserving diffeomorphisms that restrict to an orientation-preserving diffeomorphism of $B$.
Let $\Mod_G(M,B) \leq \Mod(M,B)$ and $\Mod_G(M) \leq \Mod(M)$ each denote the subgroup of mapping classes with an equivariant representative.
Then the natural forgetful map $\Mod(M,B) \rightarrow \Mod(M)$ restricts to a map
\begin{equation*}
    \F_G:\Mod_G(M,B) \rightarrow \Mod_G(M).
\end{equation*}
Then a finite index subgroup of $\Ker(\P_G)$ will surject onto $\Ker(\F_G)$, so it's enough to show that $\Ker(\F_G)$ is not finitely generated.

We construct an obstruction map
\begin{equation*}
    \varphi:\Ker(\F_G) \rightarrow H_{n-1}(\widetilde{M}, \widetilde{B}; \Z/\ell\Z)
\end{equation*}
where $\widetilde{M}$ is a certain infinite-sheeted cover of $M$ and $\widetilde{B}$ is a component of the preimage of $B$ in $\widetilde{M}$.
The homology group $H_{n-1}(\widetilde{M}, \widetilde{B}; \Z/\ell\Z)$ is non-finitely generated, so the strategy of the proof is to find an infinite linearly independent subset of the image of $\varphi$.

The map $\varphi$ can be computed as follows: suppose $\alpha \in \Ker(\F_G)$ is represented by an equivariant diffeomorphism $f:M \rightarrow M$ that preserves $B$ setwise.
Then there is an isotopy $f \simeq \id_M$.
By tracking $B$ under this isotopy, we get a map $B \times [0,1] \rightarrow M$ mapping $B \times \{0\}$ and $B \times \{1\}$ diffeomorphically onto $B$.
We show that this lifts to a map $B \times [0,1] \rightarrow \widetilde{M}$ mapping $B \times \{0\}$ and $B \times \{1\}$ diffeomorphically onto $\widetilde{B}$.
This lifted map yields a homology class in $H_{n-1}(\widetilde{M}, \widetilde{B}; \Z/\ell\Z)$.

The most difficult step is to find at least one element $\alpha_0 \in \Ker(\F_G)$ such that $\varphi(\alpha_0)$ is nonzero.
We can then obtain an infinite linearly independent subset of $H_{n-1}(\widetilde{M},\widetilde{B};\Z/\ell\Z)$ by conjugating $\alpha_0$.
To find $\alpha_0$, we construct an explicit equivariant diffeomorphism $f$ and an explicit isotopy from $f$ to $\id_M$.
As in the proof of part (i), the diffeomorphism $f$ is a lift of a summand slide.

The obstruction map $\varphi$ can be viewed as an invariant defined on the homology of the space of unparametrized embeddings $B \hookrightarrow M$ (see Remark \ref{rem:embedding-spaces}).
It would be interesting to adapt the map $\varphi$ to study the homology of embedding spaces in a broader context.


\subsection{Sharpness of Theorem \ref{mainthm:ker-not-fg} in dimension 3}

In dimension $n = 3$, the assumptions of Theorem \ref{mainthm:ker-not-fg} are fairly mild restrictions.
The Prime Decomposition Theorem implies that $M^\circ/G$ splits uniquely as a connected sum of prime $3$-manifolds, and the Poincar\'e conjecture implies that each summand is not simply connected.
Therefore, the first assumption reduces to the statement that $M^\circ/G$ has $k$ prime factors.
For the second assumption, we note that if $\Fix(g_0)$ is nonempty, then any component will be an embedded loop in $M$, and the assumption will be satisfied if any such loop is a torsion element (e.g.\ zero) in $H_1(M;\Z)$.
Under these assumptions, Theorem \ref{mainthm:ker-not-fg} is sharp in dimension $3$.
Namely, by \cite[Cor~D]{bh-3manifolds}, there is an action of $G = \Z/2\Z$ on $M = S^1 \times S^2$ satisfying the assumptions of Theorem \ref{mainthm:ker-not-fg} with $k=2$, and for this action, $\Ker(\P_G)$ is finite.

\subsection{Question: Normal generation}

In \cite[Thm~C]{bh-3manifolds}, we show that for a certain action of $G = \Z/2\Z$ on $M = (S^1 \times S^2)^{\#k-1}$, the group $\Ker(\P_G)$ is the $\Gamma_G(M)$-normal closure of a finite set (the theorem is stated in the topological category, but can easily be adapted to the smooth category).
The proof is special to this case; it exploits the isomorphism $\Mod(S^3, C_k^{\pm}) \cong \SymOut(F_k)$ described in Section \ref{subsec:appliction-to-symmetric-auts} and appeals to the action of $\SymOut(F_k)$ on a certain contractible simplicial complex constructed by McCullough--Miller \cite{mccullough-miller}.
We therefore ask:

\begin{question}\label{ques:normal-generation}
    Let $M$ be a closed oriented smooth manifold of dimension $n \geq 3$, and $G$ a finite group of orientation-preserving diffeomorphisms.
    Is $\Ker(\P_G)$ the $\Gamma_G(M)$-normal closure of a finite set?
\end{question}

\subsection{Free actions}

For non-free group actions, the fixed points of $G$ provide a natural obstruction to equivariant isotopy.
In the case of free actions, it is less obvious how to obstruct equivariant isotopy, and indeed, there are classes of free actions for which $\Ker(\P_G)$ is known to be trivial.

For instance, suppose $M$ is a closed surface of genus $g \geq 1$ and $G$ acts freely.
Then by studying the action of diffeomorphisms on $\pi_1$ and appealing to $K(\pi,1)$-theory, one can show that $\Ker(\P_G)$ is trivial (see e.g.\ \cite[\S 9]{margalit-winarski}, which gives an argument following Birman--Hilden \cite{birman-hilden} and Aramayona--Leininger--Souto \cite{aramayona-leininger-souto}).
This argument easily adapts to the case that $M$ is a closed hyperbolic $3$-manifold; the main technical ingredient is that homotopic diffeomorphisms are isotopic, which was proved by Gabai--Meyerhoff--Thurston \cite{gabai-meyerhoff-thurston}.

Beyond these examples, it appears little is known about $\Ker(\P_G)$ in dimension $n \geq 3$ when $G$ acts freely.
A recent theorem of Raman \cite{raman} uses Hodge theory to prove that $\Ker(\P_G)$ is trivial if $M$ is an irreducible hyperk\"ahler manifold and $G$ is a finite group of automorphisms acting freely.
This provides examples of free actions for which $\Ker(\P_G)$ is trivial in dimension $4n$ for $n \geq 1$; the simplest example comes from a free involution on a $K3$ surface (which has real dimension 4).
Vogt \cite{vogt} proved that $\Ker(\P_G)$ is trivial for certain free actions on Seifert fibered $3$-manifolds, and applied this result to the study of foliations of certain $5$-manifolds.

\subsection{Outline}

In Section \ref{sec:disk-slides}, we describe certain diffeomorphisms of manifolds that we call \emph{disk slides} and \emph{summand slides}; these are the diffeomorphisms that we use to prove Theorem \ref{mainthm:ker-not-fg}.
In Section \ref{sec:pf-part-i}, we prove part (i) of Theorem \ref{mainthm:ker-not-fg}.
Sections \ref{sec:model-of-action} through \ref{sec:pf-part-ii} are dedicated to the proof of part (ii) of Theorem \ref{mainthm:ker-not-fg}.
In Section \ref{sec:model-of-action}, we give local and global models for the action of $G$ on $M$, and construct a certain infinite cover $\widetilde{M}$ whose homology is not finitely generated.
In Section \ref{sec:obstruction-map}, we construct the obstruction map valued in the homology of $\widetilde{M}$.
In Section \ref{sec:nontrivial-obstruction}, we construct an equivariant diffeomorphism which maps to a nontrivial homology class under the obstruction map.
In Section \ref{sec:pf-part-ii}, we complete the proof of part (ii) of Theorem \ref{mainthm:ker-not-fg} (and Theorem \ref{mainthm:top-version}) by using the element in Section \ref{sec:nontrivial-obstruction} to build an infinite linearly independent subset of the image of the obstruction map.
In Section \ref{sec:pf-sym-aut-kernel}, we use Theorem \ref{mainthm:top-version} to prove Theorem \ref{mainthm:symmetric-aut-kernel}.

Throughout this paper, we primarily work in the smooth category.
At the end of Sections 2, 4, 5, 6, and 7, we explain how to adapt the results to the topological setting as described in Section \ref{subsec:topological-setting}.

\subsection{Acknowledgements}
We thank Bena Tshishiku for introducing us to this problem, for helpful comments on an earlier draft, and for many helpful conversations throughout this project.
We also thank Tom Goodwillie for a helpful discussion about the obstruction map constructed in Section \ref{sec:obstruction-map}.
Finally, we thank the anonymous referee for their careful reading of this paper and for many helpful corrections and suggestions.

\section{Disk Slides and Summand Slides}\label{sec:disk-slides}

In this section, we define certain ``slide diffeomorphisms'' which we will use to construct elements of $\Ker(\P_G)$.
We first define the notion of a ``disk slide'' in Section \ref{subsec:disk-slides}, and study their action on the fundamental group in Section \ref{subsec:disk-slides-inner-auts}. 
We then use disk slides to define ``summand slides'' in Section \ref{subsec:summand-slides}.

\subsection{Disk slides}\label{subsec:disk-slides}

Let $M$ be a compact oriented smooth manifold of dimension $n \geq 3$ (possibly with boundary).
An \emph{(ordered) multidisk} in $M$ is a tuple $\Delta = (D_1, \ldots, D_d)$, where $D_1, \ldots, D_d$ are disjoint closed smooth $n$-disks in the interior $\Int(M)$.
Given a multidisk $\Delta$, we define $\Diff_\del(M \rel \Delta)$ to be the group of diffeomorphisms of $M$ which fix $\del M$ and $\Delta$ pointwise, and we define $\Mod(M \rel \Delta) \coloneqq \pi_0(\Diff_\del(M \rel \Delta))$.
Observe that we have a natural forgetful map $\Mod(M \rel \Delta) \rightarrow \Mod(M)$.

\begin{definition}
  Given a smooth manifold $M$ and a multidisk $\Delta$ in $M$, we define the group
  \begin{equation*}
      \DS_\Delta(M) \coloneqq \Ker\left( \Mod(M \rel \Delta) \rightarrow \Mod(M) \right),
  \end{equation*}
  and we define a \emph{disk slide} to be any diffeomorphism of $M$ that represents an element of $\DS_\Delta$.
\end{definition}

Our goal in this subsection is to characterize disk slides in terms of the fundamental group of a certain bundle.
First, let $\Conf_d(M)$ denote the configuration space of $d$ ordered points in the interior of $M$, i.e.
\begin{equation*}
    \Conf_d(M) \coloneqq \{(x_1, \ldots, x_d) \in \Int(M)^d \mid x_i \neq x_j \text{ for all } 1 \leq i \neq j \leq d\}.
\end{equation*}
Then we define the \emph{multiframe bundle of degree $d$} to be the principal $\GL_n(\R)^d$-bundle $\Fr_d(M) \rightarrow \Conf_d(M)$ where the fiber over $(x_1, \ldots, x_n)$ is the set of tuples $(\omega_1, \ldots, \omega_n)$, where $\omega_i$ is a basis of the tangent space $T_{x_i}M$.
In the case $d=1$, $\Fr_1(M)$ is simply the frame bundle $\Fr(M) \rightarrow \Int(M)$.
To construct $\Fr_d(M)$ formally, one can start with the direct product bundle $(\Fr(M))^d$ on $M^d$ and restrict it to the subspace $\Conf_d(M) \subseteq M^d$.
Alternatively, one can construct $\Fr_d(M)$ as the subspace of $\Conf_d(\Fr(M))$ consisting of tuples of frames over distinct points.

Suppose now $\Delta$ is a \emph{framed} multidisk, meaning that $\Delta = (D_1, \ldots, D_d)$ is a multidisk equipped with a choice of frame $\omega_i$ at the center of each disk $D_i$.
Then we let $\Fr_\Delta(M)$ denote the frame bundle $\Fr_d(M)$, equipped with the base point $(\omega_1, \ldots, \omega_d)$.
The main result of this subsection is then the following.

\begin{proposition}\label{prop:frame-bundle-and-disk-slides}
    Let $\Delta$ be a framed multidisk in $M$.
    Then there is natural surjection
    \begin{equation*}
        \calDS_\Delta:\pi_1(\Fr_\Delta(M)) \twoheadrightarrow \DS_\Delta(M).
    \end{equation*}
\end{proposition}

We prove Proposition \ref{prop:frame-bundle-and-disk-slides} with two standard lemmas about embedding spaces of disks.
Given a multidisk $\Delta = (D_1, \ldots, D_d)$, we let $\Emb(\Delta, M)$ denote the space of embeddings $D_1 \sqcup \cdots \sqcup D_d \hookrightarrow \Int(M)$.
Then the following two lemmas follow from the work of Cerf \cite{cerf}.

\begin{lemma}\label{lem:disk-embedding-fibration}
    The natural map $\Diff_\del(M) \rightarrow \Emb(\Delta, M)$ is a locally trivial fibration, where $\Diff_\del(M)$ is the group of diffeomorphisms fixing $\del M$ pointwise.
\end{lemma}

Lemma \ref{lem:disk-embedding-fibration} is given in \cite[\S II.2.2.2]{cerf}.

\begin{lemma}\label{lem:embedding-space-frame-bundle-htpy-equivalent}
    The natural map $\Emb(\Delta, M) \rightarrow \Fr_\Delta(M)$ is a homotopy equivalence.
\end{lemma}

Lemma \ref{lem:embedding-space-frame-bundle-htpy-equivalent} is deduced in \cite[\S II.5.1.5]{cerf} in the case of single disk, but the proof can be adapted to the case of a multidisk with minor modifications.

Now, we can prove Proposition \ref{prop:frame-bundle-and-disk-slides}.

\begin{proof}[Proof of Proposition \ref{prop:frame-bundle-and-disk-slides}]
    By Lemma \ref{lem:disk-embedding-fibration}, we get an exact sequence
    \begin{equation*}
        \pi_1(\Emb(\Delta, M)) \rightarrow \Mod(M \rel \Delta) \rightarrow \Mod(M).
    \end{equation*}
    By exactness, the image of the first map is precisely $\DS_\Delta(M)$.
    By Lemma \ref{lem:embedding-space-frame-bundle-htpy-equivalent}, we have an isomorphism $\pi_1(\Fr_\Delta(M)) \rightarrow \pi_1(\Emb(\Delta, M))$.
    Thus we define $\calDS_\Delta$ as the composition
    \begin{equation*}
        \pi_1(\Fr_\Delta(M)) \xrightarrow{\cong} \pi_1(\Emb(\Delta, M)) \twoheadrightarrow \DS_\Delta(M).
    \end{equation*}
\end{proof}

\subsection{Disk slides, sphere twists, and inner automorphisms}\label{subsec:disk-slides-inner-auts}

Next, given a framed multidisk $\Delta = (D_1, \ldots, D_d)$ in $M$, we can study the behavior of disks slides in terms of $\pi_1(\Fr_\Delta(M))$.
Observe first that since $\Fr_\Delta(M)$ is a principal $\GL_n(\R)^d$-bundle, we have an exact sequence
\begin{equation*}
    \pi_1(\GL_n(\R)^d) \rightarrow \pi_1(\Fr_\Delta(M)) \rightarrow \pi_1(\Conf_d(M)) \rightarrow \pi_0(\GL_n(\R)^d).
\end{equation*}
Since $M$ is orientable, the last map is in fact trivial.
So, we can understand disk slides in terms of $\pi_1(\GL_n(\R)^d)$ and $\pi_1(\Conf_d(M))$.
We will summarize our findings as Proposition \ref{prop:disk-slides-twists-and-inner-auts} below.
For a similar discussion in the case of $3$-manifolds, see \cite[Rem~2.4]{hatcher-wahl}.

\subsubsection{Sphere twists}

First, recall that $\pi_1(\GL_n(\R)) \cong \Z/2\Z$, and it is generated by a $2\pi$-rotation in the $xy$-plane.
Then $\pi_1(\GL_n(\R)^d) \cong (\Z/2\Z)^d$; let $e_1, \ldots, e_d$ be the images of standard generators in $\pi_1(\Fr_\Delta(M))$. 
If we let $p_i$ denote the center of $D_i$, then $e_i$ is a $2\pi$-rotation of our chosen frame at $p_i$.
It follows that $\calDS_\Delta(e_i)$ is the isotopy class of \emph{sphere twist} about $\del D_i$.
This is a diffeomorphism of $M$ supported on a collar neighborhood $\del D_i \times [0,1] \subseteq M - \Int(D_i)$ which acts on the slice $\del D_i \times \{t\}$ by a $2\pi t$-rotation.
In particular, $\calDS_\Delta(e_i)$ acts trivially on $\pi_1(M, p_i)$.
We emphasize that a sphere twist may be isotopically trivial, i.e. the composition $(\Z/2\Z)^d \rightarrow \pi_1(\Fr_\Delta(M)) \rightarrow \DS_{\Delta}(M)$ need not be injective (in particular, if $e_i$ is in the image of the composition 
\begin{equation*}
    \pi_1(\Diff_\del(M)) \rightarrow \pi_1(\Emb(\Delta, M)) \xrightarrow{\cong} \pi_1(\Fr_\Delta(M)),
\end{equation*}
then $\calDS_\Delta(e_i)$ will be trivial).

\subsubsection{Inner automorphisms}

Next, we can study disk slides in terms of $\pi_1(\Conf_d(M))$.
We let $p_i$ denote the center of $D_i$, and fix $(p_1, \ldots, p_d)$ as a base point of $\Conf_d(M)$.
First, we can compute $\pi_1(\Conf_d(M))$.

\begin{lemma}\label{lem:pi1-of-configuration-space}
  The natural map 
  \begin{equation*}
      \pi_1(\Conf_d(M)) \rightarrow \prod_{i=1}^d \pi_1(M, p_i)
  \end{equation*}
  is an isomorphism.
\end{lemma}
\begin{proof}
    Fix local coordinates $x^1, \ldots, x^n$ on $\Int(M)$.
    Then $\Int(M)^d$ has local coordinates of the form
    \begin{equation*}
        x^1_1, \ldots, x^n_1, \ldots, x^1_d, \ldots, x^n_d,
    \end{equation*}
    and $\Conf_d(M)$ is given locally as the complement of the linear subspaces
    \begin{equation*}
        V_{i,j} = \{(x^1_1, \ldots, x^n_1, \ldots, x^1_d, \ldots, x^n_d) \in \R^{nd} \mid x_i^k = x_j^k \text{ for all $1 \leq k \leq n$}\}
    \end{equation*}
    for $1 \leq i \neq j \leq d$.
    Since the equation $x_i^k = x_j^k$ defines a codimension $1$ subspace, it follows that $V_{i,j}$ has codimension $n$.
    
    Thus, since $n \geq 3$, the space $\Conf_d(M)$ is obtained from $\Int(M)^d$ by removing a transverse intersection of closed submanifolds of codimension at least $3$, which means that the map $\pi_1(\Conf_d(M)) \rightarrow \pi_1(M^d)$ induces an isomorphism on $\pi_1$.
\end{proof}




By Lemma \ref{lem:pi1-of-configuration-space}, we have a natural map
\begin{equation*}
    \rho:\pi_1(\Fr_\Delta(M)) \rightarrow \prod_{i=1}^d \pi_1(M, p_i).
\end{equation*}
We saw above that $\calDS_\Delta$ sends elements in $\Ker(\rho)$ to sphere twists, which act trivially on $\pi_1(M)$.
Thus, we will use the map $\rho$ to study the action of disk slides on $\pi_1(M)$.

Observe that for each $j \in \{1, \ldots, d\}$, any element of $\Mod(M \rel \Delta)$ fixes the point $p_j \in D_j$, and hence $\Mod(M \rel \Delta)$ acts on $\pi_1(M, p_j)$.
If $\alpha \in \DS_\Delta(M)$ is represented by a disk slide $f$, then since $f$ is isotopic to $\id_M$, $\alpha$ must induce an inner automorphism on $\pi_1(M,p_j)$.
Indeed, if we choose an isotopy $h_t$ with $h_0 = \id$ and $h_1 = f$, then $f$ will conjugate $\pi_1(M, p_j)$ by the path of $p_j$ under $h_t$.

Now, consider the composition
\begin{equation*}
    \pi_1(\Fr_\Delta(M)) \xrightarrow{\calDS_\Delta} \Mod(M \rel \Delta) \rightarrow \Aut(\pi_1(M, p_j)).
\end{equation*}
It follows from the construction of $\calDS_\Delta$ and the discussion above that for $\gamma \in \pi_1(\Fr_\Delta(M))$, the image of $\gamma$ in $\Aut(\pi_1(M, p_j))$ is the inner automorphism $\Inn(\overline{\gamma})$, where $\overline{\gamma}$ is the image of $\gamma$ in $\pi_1(M,p_j)$.
Indeed, $\gamma$ determines a loop in $\Emb(\Delta, M)$, which we can extend to an ambient isotopy $h_t$ with $h_0 = \id_M$. 
Then $\calDS_\Delta(\gamma)$ is represented by $h_1$. 
By the discussion above, $\calDS_\Delta(\gamma)$ conjugates $\pi_1(M,p_j)$ by the path traveled by $p_j$ under $h_t$, which is precisely the loop $\overline{\gamma}$.

\subsubsection{Conclusion}

All together, we can summarize the discussion in this subsection as the following proposition.

\begin{proposition}\label{prop:disk-slides-twists-and-inner-auts}
    Let $\Delta = (D_1, \ldots, D_d)$ be a framed multidisk in $M$, and let $p_i$ be the center of the disk $D_i$.
    Then there is an exact sequence
    \begin{equation*}
        (\Z/2\Z)^d \rightarrow \pi_1(\Fr_\Delta(M)) \rightarrow \prod_{i=1}^d \pi_1(M, p_i) \rightarrow 1.
    \end{equation*}
    Moreover:
    \begin{enumerate}[label=(\roman*)]
        \item If $e_i$ is the image of the $i$th generator under the map $(\Z/2\Z)^d \rightarrow \pi_1(\Fr_\Delta(M))$, then $\calDS_\Delta(e_i)$ is the isotopy class of a sphere twist about $\del D_i$.
        \item For each $j \in \{1, \ldots, d\}$, the composite map
        \begin{equation*}
            \pi_1(\Fr_\Delta(M)) \xrightarrow{\calDS_\Delta} \Mod(M \rel \Delta) \rightarrow \Aut(\pi_1(M, p_j))
        \end{equation*}
        is equal to the composition
        \begin{equation*}
            \pi_1(\Fr_\Delta(M)) \rightarrow \prod_{i=1}^d \pi_1(M, p_i) \rightarrow \pi_1(M, p_j) \rightarrow \Inn(\pi_1(M, p_j)).
        \end{equation*}
    \end{enumerate}
\end{proposition}

\subsection{Summand slides}\label{subsec:summand-slides}

Next, we will use disk slides to define \emph{summand slides}, which are the diffeomorphisms we will use to prove Theorem \ref{mainthm:ker-not-fg}.
The main result of this subsection is Proposition \ref{prop:summand-slide-action-on-pi1} below, which computes the action of a summand slide on $\pi_1$.

First, we introduce the notion of gluing manifolds along multidisks (following \cite{kervaire-milnor}).
Suppose $M_1$ and $M_2$ are two oriented smooth $n$-manifolds, and $\Delta_i = (D_{i,1}, \ldots, D_{i,d})$ is a multidisk in the interior of $M_i$.
Let $p_{i,j}$ denote the center of the disk $D_{i,j}$.
We choose parametrizations $k_{i,j}:D^n \rightarrow D_{i,j}$, where $D^n$ is the standard closed $n$-disk in $\R^n$.
We require each $k_{1,j}$ to be orientation-preserving and $k_{2,j}$ to be orientation-reversing.
Define 
\begin{equation*}
    M_i^* \coloneqq M_i \setminus \{p_{i,1}, \ldots, p_{i,d}\}.
\end{equation*}

To \emph{glue} $M_1$ and $M_2$ along $\Delta_1$ and $\Delta_2$, we start with the disjoint union $M_1^* \sqcup M_2^*$.
Then, using our parametrizations $k_{i,j}$, we identify $\Int(D_{1,j}) \setminus \{p_{1,j}\}$ with $\Int(D_{2,j}) \setminus \{p_{2,j}\}$ via the diffeomorphism of $\Int(D^n) \setminus \{0\}$ defined by $tu \mapsto (1-t)u$ for each unit vector $u \in S^{n-1}$ and $t \in (0,1)$.
The resulting manifold $M$ is well-defined independent of the choice of $k_{i,j}$, but if $M_1$ and $M_2$ are both oriented, then $M$ may depend on the choice of orientation of $M_1$ and $M_2$ (see \cite{Palais-extending-diffeos}).
We will frequently identify $M_i^*$ with its image in $M$.

With this convention set, we can define summand slides.

\begin{definition}
    Suppose $M$ is obtained by gluing two manifolds $M_1$ and $M_2$ along multidisks $\Delta_1$ and $\Delta_2$.
    A \emph{summand slide} is a diffeomorphism of $M$ which acts on $M_1^*$ by the restriction of a disk slide, and which acts on $M_2^*$ by $\id_{M_2^*}$.
\end{definition}

Note that any summand slide is well-defined since a disk slide on $M_1$ fixes the multidisk $\Delta_1$ pointwise.

An important special case is that $d = 1$, so $M$ a connnected sum $M_1 \# M_2$.
In this case, we can easily compute the action of a summand slide on $\pi_1(M)$.
Namely, suppose $M_1$ and $M_2$ are glued along disks $D_1 \subseteq M_1$ and $D_2 \subseteq M_2$.
Let $p_i$ be the center of $D_i$.
Fix points $q_i \in D_i$ away from the center so that each $q_i$ corresponds to the same point $q \in M$.
By choosing arcs in each $D_i$ from $p_i$ to $q_i$, we get a natural free product decomposition
\begin{equation*}
    \pi_1(M, q) \cong \pi_1(M_1, p_1) * \pi_1(M_2, p_2).
\end{equation*}
Since any summand slide fixes $D_1$ and $D_2$ pointwise, it will induce well-defined automorphisms of $\pi_1(M,q)$ and $\pi_1(M_i, p_i)$ (in particular, it will act trivially on $\pi_1(M_2, p_2)$).
Then we get the following.

\begin{proposition}\label{prop:summand-slide-action-on-pi1}
  Suppose $M = M_1 \# M_2$ as above.
  Choose $\gamma \in \pi_1(\Fr(M_1))$ and let $\overline{\gamma} \in \pi_1(M_1)$ be its image.
  Suppose $f$ is a summand slide of $M$ where the action on $M_1^*$ represents $\calDS_{D_1}(\gamma)$.
  Then, with respect to the identification
  \begin{equation*}
      \pi_1(M,q) \cong \pi_1(M_1,p_1) * \pi_1(M_2,p_2),
  \end{equation*}
  the diffeomorphism $f$ acts on $\pi_1(M)$ by conjugating $\pi_1(M_1)$ by $\overline{\gamma}$ and fixing $\pi_1(M_2)$ pointwise.
\end{proposition}
\begin{proof}
  This follows directly from the definition of $f$ and Proposition \ref{prop:disk-slides-twists-and-inner-auts} above.
\end{proof}

\subsection{Adaptations to the topological setting}\label{subsec:disk-slides-adaptation}

For topological manifolds, our notion of gluing along disks carries over without issue so long as we restrict to locally flat disks.

The key result of this section, Proposition \ref{prop:disk-slides-twists-and-inner-auts}, carries over to the topological setting, so long as we replace $\Fr_\Delta(M)$ with the space of locally flat embeddings $\Emb^{\mathrm{lf}}(\Delta, M)$ for an (unframed) ordered multidisk $\Delta$.
Namely, following the proof of Lemma 2.12 in \cite{ayala-francis}, there is a fibration 
\begin{equation*}
    \Homeo(\R^n, 0)^d \rightarrow \Emb^{\mathrm{lf}}(\Delta, M) \rightarrow \Conf_d(M)
\end{equation*}
where $\Homeo(\R^n, 0)$ is the group of homeomorphisms of $\R^n$ fixing the origin.
Lemma \ref{lem:pi1-of-configuration-space} carries over verbatim in the topological category.
Thus to prove the topological analog of Proposition \ref{prop:disk-slides-twists-and-inner-auts}, it's enough to check that for $n \geq 3$, the group $\pi_1(\Homeo(\R^n, 0))$ is order 2 and generated by a rotation in the $xy$-plane. 
This follows from \cite[Essay~V,\S 5]{kirby-siebenmann}, as well as \cite[Cor~2.2.3]{quinn} for the case of dimension $n=4$.

\section{Proof of Part (i) of Theorem \ref{mainthm:ker-not-fg}}\label{sec:pf-part-i}

In this section, we prove part (i) of Theorem \ref{mainthm:ker-not-fg}.
The proof is a relatively straightforward application of Propositions \ref{prop:disk-slides-twists-and-inner-auts} and \ref{prop:summand-slide-action-on-pi1}, and is essentially a warmup to part (ii) of Theorem \ref{mainthm:ker-not-fg}.

In Section \ref{subsec:preliminaries}, we establish some notation and a model of the group action that we will use throughout this section.
In Sections \ref{subsec:lifting-disk-slides} and \ref{subsec:lifting-and-including}, we describe how disk slides lift along the quotient map $M \rightarrow M/G$.
In Section \ref{subsec:fixed-set-meridian}, we reinterpret the assumption on the homology class of $B$ in terms of $\pi_1(M \setminus B)$.
In Section \ref{subsec:constructing-subgroup} we construct the candidate subgroup of $\Ker(\P_G)$, and we complete the proof in Section \ref{subsec:completing-pf-part-i}. 

\subsection{Preliminaries}\label{subsec:preliminaries}

We begin by establishing some notation.
Let $M$ and $G$ be a manifold and group respectively which satisfy the hypotheses of Theorem \ref{mainthm:ker-not-fg} with $k = 3$.
Let $R \subseteq M$ denote the set of points whose $G$-stabilizer is nontrivial.
Then $R$ is a deformation retract of a $G$-invariant open submanifold $U \subseteq M$.
Indeed, by \cite{illman-finite-group}, $M$ admits a $G$-equivariant triangulation; the subset $R$ will be a $G$-invariant subcomplex,\footnote{
    To see this, it's enough to show that for each $g \in G$, the fixed set $\Fix(g)$ is a subcomplex of $M$.  This follows because, as noted in \cite[\S 1]{illman-finite-group}, we can pass to the second barycentric subdivision to ensure that $G$ acts ``without rotations,'' i.e.\ if $x \in \Fix(g)$ and $x$ lies in the interior of a simplex $s$, then $s$ is fixed pointwise by $g$.
} and so we can take $U$ to be the interior of a regular neighborhood of $R$.
Let $X = M \setminus U$, so $X$ is a manifold with boundary on which $G$ acts freely.
Let $Y = X/G$.
Note that if we let $M^\circ \subseteq M$ denote the set of points whose $G$-stabilizer is trivial, then $X$ is a deformation retract of $M^\circ$ and $Y$ is a deformation retract of $M^\circ/G$.

Next, we fix connected sum decompositions of $X$ and $Y$ as follows.
By assumption, we have a connected sum decomposition $M^\circ/G \cong Q_0 \# Q_1 \# Q_2$ where $\pi_1(Q_i)$ is nontrivial for each $i$.
Let $N_1 = Q_0 \# Q_1$ and $N_2 = Q_2$.
Following the convention of Section \ref{subsec:summand-slides}, we choose $n$-disks $D_1 \subseteq N_1$ and $D_2 \subseteq N_2$ and view $M^\circ/G$ as obtained by gluing $N_1$ and $N_2$ along these disks.
We let $N_i^*$ (resp.\ $D_i^*$) denote $N_i$ (resp.\ $D_i$) minus the center of $D_i$, and we identify $N_i^*$ with its image in $M^\circ/G$.
Without loss of generality, we may assume that each $D_i^*$ lies in the interior of $Y$.\footnote{
    To see this, let $E \subseteq M^\circ$ denote the preimage of $D_1^* \cup D_2^*$, so $E$ is diffeomorphic to $G \times (S^{n-1} \times [0,1])$.
    Since $R$ and $E$ are compact and disjoint, we can choose the subset $U \subseteq M$ sufficiently small (e.g.\ by passing to higher barycentric subdivisions of $M$) so that $E$ lies outside the closure of $U$.
}
We let $Y_i^* = Y \cap N_i^*$, and let $Y_i$ be the manifold obtained by filling in the puncture;
we view $Y$ as being obtained by gluing $Y_1$ and $Y_2$ along $D_1$ and $D_2$.
Let $X_i^* \subseteq X$ denote the preimage of $Y_i^*$, and let $\Delta_i^*$ denote the preimage of $D_i^*$.
Then let $X_i$ and $\Delta_i$ be obtained by filling in the punctures.
Thus we can view $X$ as obtained by gluing together $X_1$ and $X_2$ along $\Delta_1$ and $\Delta_2$ (note that each $X_i$ may be disconnected).
Note that we have regular (unbranched) $G$-coverings $X_i \rightarrow Y_i$ and $X_i^* \rightarrow Y_i^*$.

Recall that by assumption, $M$ has a distinguished codimension $2$ compact orientable submanifold $B$ which is fixed by an element $g_0 \in G$.
The submanifold $B$ is contained in $R$.
Let $U_B$ be the component of $U$ containing $B$.
Then $\del U_B$ is a boundary component of $X$; we assume up to relabeling that this boundary component lies on $X_1^*$.
In particular, $X_1$ and $Y_1$ both have nonempty boundary.

Note that we have a further splitting of $Y_1$.
Namely, since $N_1 = Q_0 \# Q_1$, a similar analysis as above\footnote{
    Namely, we can again shrink $U$ if necessary to assume that the punctured disks realizing the connected sum decomposition of $N_1$ lie in $Y_1$.
} shows that we have a splitting $Y_1 \cong Q_0' \# Q_1'$ where $Q_i'$ is a deformation retract of $Q_i$.
This implies in particular that 
\begin{equation*}
    \pi_1(Y_1) \cong \pi_1(Q_0) * \pi_1(Q_1).
\end{equation*}
If we let $(Q_i')^*$ be obtained by puncturing a point from $Q_i'$, then as above we can view $(Q_i')^*$ as a subset of $Y_1$.
Then the boundary component $\del U_B \subseteq \del X_1$ will project into either $(Q_0')^*$ or $(Q_1')^*$;
up to relabeling, we assume that it projects into $(Q_0')^*$.

\begin{figure}[h]
    \centering
    \includegraphics[scale=.6]{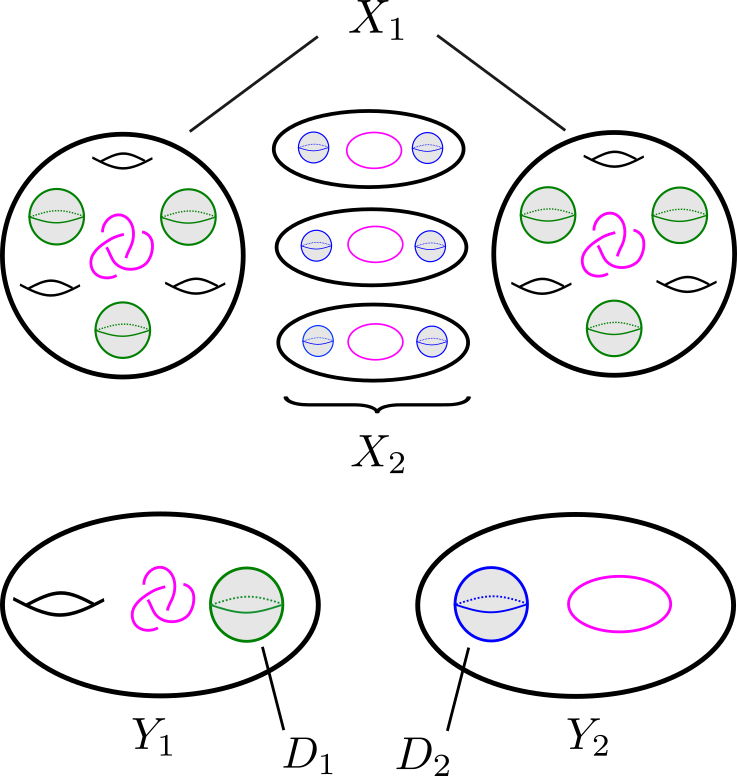}
    \caption{A schematic illustration of the manifolds $X$ and $Y$.  The trefoils and unknots (in magenta) represent the set $R$ and its image in $M/G$; the manifolds $X$ and $Y$ are the complement of a tubular neighborhood of these curves.}
    \label{fig:model-pf-pt-i}
\end{figure}

Observe that the inclusion $X \hookrightarrow M$ induces a way to fill in each boundary component of $X_1$ and $X_2$.
We call the resulting filled manifolds $M_1$ and $M_2$ respectively, so $M$ is obtained by gluing $M_1$ and $M_2$ along the multidisks $\Delta_1$ and $\Delta_2$.
See Figure \ref{fig:model-pf-pt-i} for an illustration.

Finally, we establish some notation for the various mapping class groups that arise throughout this section.
Given an oriented manifold $W$, we let $\Diff_\del(W)$ denote the group of orientation-preserving diffeomorphisms fixing $\del W$ pointwise, and let $\Mod(W) = \pi_0(\Diff_\del(W))$.
Similarly, given a closed submanifold $Z \subseteq W$, we let $\Diff_\del(W \rel Z)$ denote the subgroup of diffeomorphisms fixing $Z$ pointwise, and let $\Mod(W \rel Z) = \pi_0(\Diff_\del(W \rel Z))$.
Suppose $W' \rightarrow W$ is a finite cover.
Then we let $\LDiff_\del(W) \leq \Diff_\del(W)$ and $\LDiff_\del(W \rel Z) \leq \Diff_\del(W \rel Z)$ denote the finite index subgroups of diffeomorphisms that lift to a diffeomorphism of $W'$ that fixes $\del W'$ pointwise, and define $\LMod(W) = \pi_0(\LDiff_\del(W))$ and $\LMod(W \rel Z) = \pi_0(\LDiff_\del(W \rel Z))$.

\subsection{Lifting disk slides}\label{subsec:lifting-disk-slides}

Our next goal is to describe lifts of summand slides of $Y$ on the level of frame bundles.
It will be enough to describe lifts of disk slides on $Y_1$.


Fix an arbitrary framing on $D_1$.
Since the covering map $X_1 \rightarrow Y_1$ is a local diffeomorphism, this induces framings on $\Delta_1$.
Recall from Proposition \ref{prop:disk-slides-twists-and-inner-auts} that we have an exact sequence
\begin{equation*}
    \Z/2\Z \rightarrow \pi_1(\Fr_{D_1}(Y_1)) \rightarrow \pi_1(Y_1) \rightarrow 1.
\end{equation*}
Here $\Fr_{D_1}(Y_1)$ is simply the frame bundle $\Fr(Y_1)$ with a base point at our chosen frame at the center of $D_1$.
Let $\pi_1(Y_1)^{\mathrm{lift}} \leq \pi_1(Y_1)$ be the subgroup of loops that lift to loops in $X_1$ (equivalently, $\pi_1(Y_1)^{\mathrm{lift}}$ is the image of the fundamental group of any component of $X_1$).
Let $\pi_1(\Fr_{D_1}(Y_1))^{\lift} \leq \pi_1(\Fr_{D_1}(Y_1))$ be the preimage of $\pi_1(Y_1)^{\mathrm{lift}}$.

Now, observe that there's a natural map $\Fr_{\Delta_1}(X_1) \rightarrow \Fr_{D_1}(Y_1)$ which maps a tuple $(\omega_1, \ldots, \omega_d)$ to the image of $\omega_1$ in $Y_1$.
Then we get the following.

\begin{lemma}\label{lem:lifting-slides}
    The map $\pi_1(\Fr_{\Delta_1}(X_1)) \rightarrow \pi_1(\Fr_{D_1}(Y_1))$ above admits a section
    \begin{equation*}
        \pi_1(\Fr_{D_1}(Y_1))^{\lift} \rightarrow \pi_1(\Fr_{\Delta_1}(X_1))
    \end{equation*}
    which fits into the following commutative square:
    \begin{equation*}
        \begin{tikzcd}
            \pi_1(\Fr_{D_1}(Y_1))^{\lift} & {\pi_1(\Fr_{\Delta_1}(X_1))} \\
            {\LMod(Y_1 \rel D_1)} & {\Mod(X_1 \rel \Delta_1)}
            \arrow[from=1-1, to=1-2]
            \arrow["{\calDS_{D_1}}"', from=1-1, to=2-1]
            \arrow["{\calDS_{\Delta_1}}", from=1-2, to=2-2]
            \arrow[dashed, from=2-1, to=2-2]
        \end{tikzcd}
    \end{equation*}
    Here, the dashed arrow is defined on a finite index subgroup containing the image of $\pi_1(\Fr_{D_1}(Y_1))^{\lift}$, and sends a diffeomorphism to its unique lift that fixes $\del X_1$.
\end{lemma}

Before proving the lemma, we will have to introduce \emph{unordered} embedding spaces and multiframe bundles.

Let $d = \v G \v$, so $X_1 \rightarrow Y_1$ is a degree $d$ cover.
Note that since $\Delta_1$ is simply the preimage of $D_1$ in $X_1$, we can identify $\Emb(\Delta_1, X_1)$ with the space of embeddings $\bigsqcup_{i=1}^d D_1 \hookrightarrow \Int(X_1)$.
In particular, $\Emb(\Delta_1, X_1)$ admits a natural action of the symmetric group $S_d$ by precomposition; we define the quotient
\begin{equation*}
    \UEmb(\Delta_1, X_1) \coloneqq \Emb(\Delta_1, X_1)/S_d,
\end{equation*}
so $\UEmb(\Delta_1, X_1)$ is the space of unordered emebddings.

Similarly, we have a natural action of $S_d$ on $\Conf_d(X_1)$; we let $\UConf_d(X_1)$ denote the quotient by this action, so $\UConf_d(X_1)$ is the configuration space of $d$ unordered points in $\Int(X_1)$.
Then the multiframe bundle $\Fr_{\Delta_1}(X_1)$ descends to a bundle
\begin{equation*}
    \UFr_{\Delta_1}(X_1) \rightarrow \UConf_d(X_1),
\end{equation*}
where the fiber over an unordered set of points $\{x_1, \ldots, x_d\}$ is an unordered set  $\{\omega_1, \ldots, \omega_d\}$ where $\omega_i$ is a basis for $T_{x_i}X_1$.

Note that we have the following commutative diagram with exact rows:
\begin{equation}\label{eqn:uembd-and-ufr}
    \begin{tikzcd}
        1 & {\pi_1(\Emb(\Delta_1,X_1))} & {\pi_1(\UEmb(\Delta_1,X_1))} & {S_d} & 1 \\
        1 & {\pi_1(\Fr_{\Delta_1}(X_1))} & {\pi_1(\UFr_{\Delta_1}(X_1))} & {S_d} & 1
        \arrow[from=1-1, to=1-2]
        \arrow[from=1-2, to=1-3]
        \arrow[from=1-2, to=2-2]
        \arrow[from=1-3, to=1-4]
        \arrow[from=1-3, to=2-3]
        \arrow[from=1-4, to=1-5]
        \arrow[from=1-4, to=2-4]
        \arrow[from=2-1, to=2-2]
        \arrow[from=2-2, to=2-3]
        \arrow[from=2-3, to=2-4]
        \arrow[from=2-4, to=2-5]
    \end{tikzcd}
\end{equation}
By the 5 Lemma, we conclude that $\pi_1(\UEmb(\Delta_1,X_1)) \cong \pi_1(\UFr_{\Delta_1}(X_1))$.

Now, we can prove the lemma.

\begin{proof}[Proof of Lemma \ref{lem:lifting-slides}]
    Observe that by lifting an embedding $D_1 \hookrightarrow Y_1$, we get a natural map
    \begin{equation*}
        \Emb(D_1, Y_1) \rightarrow \UEmb(\Delta_1, X_1).
    \end{equation*}
    This map fits into the following commutative diagram:
    \begin{equation}\label{eqn:fibration-diagram}
        \begin{tikzcd}
            {\LDiff_\del(Y_1 \rel D_1)} & {\LDiff_\del(Y_1)} & {\Emb(D_1, Y_1)} \\
            {\Diff_{\del}(X_1 \rel_{S_d} \Delta_1)} & {\Diff_{\del}(X_1)} & {\UEmb(\Delta_1, X_1)}
            \arrow[from=1-1, to=1-2]
            \arrow[from=1-1, to=2-1]
            \arrow[from=1-2, to=1-3]
            \arrow[from=1-2, to=2-2]
            \arrow[from=1-3, to=2-3]
            \arrow[from=2-1, to=2-2]
            \arrow[from=2-2, to=2-3]
        \end{tikzcd}
    \end{equation}
    Here, the middle vertical map sends a liftable diffeomorphism to the unique lift that fixes $\del X_1$ (recall that $X_1$ and $Y_1$ both have a boundary component coming from $B$).
    The group $\Diff_{\del}(X_1 \rel_{S_d} \Delta_1)$ is the group of diffeomorphisms of $X_1$ fixing $\del X_1$ pointwise and fixing the natural embedding $\Delta_1 \hookrightarrow X_1$ up to precomposition by an element of $S_d$, and the left vertical map is the restriction of the middle vertical map.

    Observe that each row in the diagram (\ref{eqn:fibration-diagram}) is a fibration.
    Indeed, Lemma \ref{lem:disk-embedding-fibration} tells us (by replacing $M$ with $Y_1$) that $\Diff_\del(Y_1) \rightarrow \Emb(D_1, Y_1)$ is a fibration; since $\LDiff_\del(Y_1)$ is a union of components of $\Diff_\del(Y_1 )$, the same conclusion holds when we restrict to $\LDiff_\del(Y_1)$.
    Similarly, Lemma \ref{lem:disk-embedding-fibration} tells us (by replacing $M$ with $X_1$) that $\Diff_\del(X_1) \rightarrow \Emb(\Delta_1, X_1)$ is a fibration.
    Since the quotient map $\Emb(\Delta_1,X_1) \rightarrow \UEmb(\Delta_1, X_1)$ is a covering map (as $S_d$ acts freely on $\Emb(\Delta_1,X_1)$), it follows that the composite map $\Diff_\del(X_1) \rightarrow \UEmb(\Delta_1, X_1)$ is a fibration as well.

    Now, observe that by mapping a frame in $Y_1$ to the set of frames above it in $X_1$, we have a natural map
    \begin{equation*}
        \Fr_{D_1}(Y_1) \rightarrow \UFr_{\Delta_1}(X_1).
    \end{equation*}
    Combining this map with the diagram (\ref{eqn:uembd-and-ufr}) and the long exact sequences of the fibrations in the diagram (\ref{eqn:fibration-diagram}), we get the following commutative diagram:
    \begin{equation*}
        \begin{tikzcd}
            {\pi_1(\Fr_{D_1}(Y_1))} & {\pi_1(\UFr_{\Delta_1}(X_1))} \\
            {\pi_1(\Emb(D_1, Y_1))} & {\pi_1(\UEmb(\Delta_1, X_1))} \\
            {\LMod(Y_1 \rel D_1)} & {\Mod(X_1 \rel_{S_d} \Delta_1)}
            \arrow[from=1-1, to=1-2]
            \arrow["\cong"', from=1-1, to=2-1]
            \arrow["\cong", from=1-2, to=2-2]
            \arrow[from=2-1, to=2-2]
            \arrow[from=2-1, to=3-1]
            \arrow[from=2-2, to=3-2]
            \arrow[from=3-1, to=3-2]
        \end{tikzcd}
    \end{equation*}
    Here $\Mod(X_1 \rel_{S_d} \Delta_1) = \pi_0(\Diff_{\del}(X_1 \rel_{S_d} \Delta_1))$.

    Now, we can use the above diagram to deduce the desired commutative square.
    Following the proof of Proposition \ref{prop:frame-bundle-and-disk-slides}, we see that the composition of the left vertical maps is precisely $\calDS_{D_1}$.
    Next, recall from (\ref{eqn:uembd-and-ufr}) that $\pi_1(\Fr_{\Delta_1}(X_1))$ is naturally a subgroup of $\pi_1(\UFr_{\Delta_1}(X_1))$ (and the analogous statement holds for $\Emb$ and $\UEmb$).
    Composing the right vertical maps in the above diagram yields a map $\pi_1(\UFr_{\Delta_1}(X_1)) \rightarrow \Mod(X_1 \rel_{S_d} \Delta)$, and the restriction to the subgroup $\pi_1(\Fr_{\Delta}(X_1))$ is precisely the map $\calDS_{\Delta_1}$.
    Finally, it follows by the definition of $\pi_1(\Fr_{D_1}(Y_1))^{\lift}$ that the image of $\pi_1(\Fr_{D_1}(Y_1))^{\lift}$ under the map $\pi_1(\Fr_{D_1}(Y_1)) \rightarrow \pi_1(\UFr_{\Delta_1}(X_1))$ lands in the subgroup $\pi_1(\Fr_{\Delta_1}(X_1))$.
    This yields the desired commutative square, where the dashed arrow is defined on the preimage of the finite index subgroup $\Mod(X_1 \rel \Delta_1) \leq \Mod(X_1 \rel_{S_d} \Delta_1)$.
    The fact that $\pi_1(\Fr_{D_1}(Y_1))^{\lift} \rightarrow \pi_1(\Fr_{\Delta_1}(X_1))$ is a section follows by construction.

\end{proof}

\subsection{Lifting and including}\label{subsec:lifting-and-including}

Now, we have the following commutative diagram:
\begin{equation}\label{eqn:lifting-diagram}
    \begin{tikzcd}
        \pi_1(\Fr_{D_1}(Y_1))^{\lift} & {\pi_1(\Fr_{\Delta_1}(X_1))} & {\pi_1(\Fr_{\Delta_1}(M_1))} \\
        {\LMod(Y_1 \rel D_1)} & {\Mod(X_1 \rel \Delta_1)} & {\Mod(M_1 \rel \Delta_1)} \\
        {\LMod(Y)} & {\Mod(X)} & {\Mod(M)}
        \arrow[from=1-1, to=1-2]
        \arrow["\calDS_{D_1}", from=1-1, to=2-1]
        \arrow[from=1-2, to=1-3]
        \arrow["\calDS_{\Delta_1}", from=1-2, to=2-2]
        \arrow["\calDS_{\Delta_1}", from=1-3, to=2-3]
        \arrow[dashed, from=2-1, to=2-2]
        \arrow[from=2-1, to=3-1]
        \arrow[from=2-2, to=2-3]
        \arrow[from=2-2, to=3-2]
        \arrow[from=2-3, to=3-3]
        \arrow[from=3-1, to=3-2]
        \arrow[from=3-2, to=3-3]
    \end{tikzcd}
\end{equation}
The top-left square comes directly from Lemma \ref{lem:lifting-slides}.
The top right square is induced by the inclusion $X_1 \hookrightarrow M_1$.
In the bottom right square, the horizontal maps are induced by the inclusions $X_1 \hookrightarrow M_1$ and $X \hookrightarrow M$, and the vertical maps are induced by the inclusions $X_1^* \hookrightarrow X$ and $M_1^* \hookrightarrow M$.
In the bottom left square, the bottom horizontal map is induced by sending a diffeomorphism to its unique lift that fixes $\del X$, and the vertical maps are again induced by inclusions.

\subsection{A meridian of the fixed set}\label{subsec:fixed-set-meridian}

Next, we will identify a loop in $M$ that we will use to construct the summand slides needed to prove part (i) of Theorem \ref{mainthm:ker-not-fg}.
This loop will be a \emph{meridian} of $B$, i.e.\ the boundary of a fiber of a tubular neighborhood of $B$.

\begin{lemma}\label{lem:fixed-set-meridian}
    The submanifold $B$ has a meridian $\mu$ which lies in $X$ and represents a nontrivial element of $\pi_1(X)$.
\end{lemma}

To prove Lemma \ref{lem:fixed-set-meridian}, we appeal to the following general fact about meridians of codimension $2$ submanifolds.

\begin{lemma}\label{lem:meridian-vs-homology-class}
    Let $W$ be a closed oriented smooth manifold of dimension $n \geq 3$, and let $Z \subseteq W$ be a compact oriented codimension $2$ submanifold without boundary.
    Let $\mu$ be a meridian of $Z$.
    Then for any prime $p$, the meridian $\mu$ represents a nontrivial element of $H_1(W \setminus Z; \Z/p\Z)$ if and only if the homology class $[Z] \in H_{n-2}(W;\Z/p\Z)$ is trivial.
\end{lemma}
\begin{proof}
    Let $\nu(Z)$ denote the normal bundle of $Z$, so $\nu(Z)$ is a $D^2$-bundle over $Z$.
    By excision, we have an isomorphism
    \begin{equation*}
        H_2(W, W \setminus Z; \Z/p\Z) \cong H_2(\nu(Z), \del \nu(Z); \Z/p\Z).
    \end{equation*}
    By the Thom isomorphism and Poincar\'e duality, we have an isomorphism $H_2(\nu(Z), \del \nu(Z); \Z/p\Z) \cong \Z/p\Z$, and a homology class in $H_2(\nu(Z), \del \nu(Z); \Z/p\Z)$ is determined by its $\Z/p\Z$-intersection number with $Z$.
    This implies in particular that $H_2(\nu(Z), \del\nu(Z); \Z/p\Z)$ is generated by a fiber of $\nu(Z)$.

    Now, from the long exact sequence of the pair $(W, W \setminus Z)$, we have an exact sequence
    \begin{equation*}
        H_2(W; \Z/p\Z) \rightarrow H_2(\nu(Z), \del\nu(Z); \Z/p\Z) \rightarrow H_1(W \setminus Z; \Z/p\Z).
    \end{equation*}
    The second map sends a fiber of $\nu(Z)$ to its boundary, i.e.\ to the homology class of $\mu$.
    The first map sends a homology class of $W$ 
    to its $\Z/p\Z$-intersection number with $Z$.
    Thus, the homology class $[Z] \in H_{n-2}(W;\Z/p\Z)$ is trivial if and only if the first map is trivial, which occurs if and only if the second map is injective, which occurs if and only if the meridian $\mu$ is nontrivial in $H_1(W \setminus Z; \Z/p\Z)$.
\end{proof}

Now, we can prove Lemma \ref{lem:fixed-set-meridian}.

\begin{proof}[Proof of Lemma \ref{lem:fixed-set-meridian}]
    Fix $b \in B$, and let $G_b$ denote the $G$-stabilizer of $b$.
    We may assume that $G$ acts by isometries with respect to some Riemannian metric on $M$.
    Then the exponential map $T_bM \rightarrow M$ is a $G_b$-equivariant diffeomorphism near the origin.
    Since $G$ acts by orientation-preserving diffeomorphisms, we know that $G_b$ does not fix any codimension $1$ subspace of $T_bM$.
    Thus near $b$, the preimage of $M^\circ$ in $T_bM$ is the complement of a finite union of subspaces of codimension at least $2$.
    It follows that we can find a meridian $\mu$ of $B$ contained in $M^\circ$, and up to homotopy, we may assume $\mu$ lies in $X$.
    Then by Lemma \ref{lem:meridian-vs-homology-class} (applied to some prime $p$ dividing $\ell$), the class of $\mu$ in $\pi_1(X)$ maps to a nontrivial element under the composition
    \begin{equation*}
        \pi_1(X) \rightarrow \pi_1(M \setminus B) \rightarrow H_1(M \setminus B; \Z) \rightarrow H_1(M\setminus B; \Z/\ell\Z),
    \end{equation*}
    and is therefore nontrivial.
\end{proof}

\subsection{Constructing the subgroup}\label{subsec:constructing-subgroup}

Our next goal is to construct a subgroup $\Lambda \leq \pi_1(\Fr_{D_1}(Y_1))^{\lift}$ which maps onto a virtual free product in $\LMod(Y)$, and which maps to a finite subgroup in $\Mod(M)$.
From here, we will be able to deduce part (i) of Theorem \ref{mainthm:ker-not-fg}.

\begin{figure}[h]
    \centering
    \includegraphics[scale=.6]{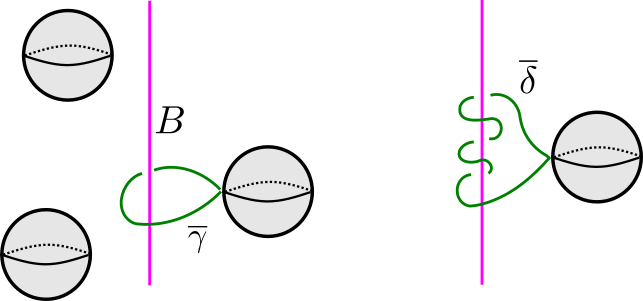}
    \caption{The curves $\overline{\gamma}$ and $\overline{\delta}$.
    The $3$-disk on the right represents $D_1$, and the $3$-disks on the left represent a subset of $\Delta_1$.}
    \label{fig:meridian-and-image}
\end{figure}

We construct the subgroup $\Lambda$ as follows.
Let $\mu$ be the meridian of $B$ from Lemma \ref{lem:fixed-set-meridian}.
Since $B$ lies in $M_1$, we may assume without loss of generality that $\mu$ is contained in $X_1$ and represents a nontrivial conjugacy class in $\pi_1(X_1)$.
Let $\overline{\gamma} \in \pi_1(X_1)$ be an element in this conjugacy class.
Let $\overline{\delta} \in \pi_1(Y_1)$ be the image of $\overline{\gamma}$; by construction, $\overline{\delta}$ lies in the subgroup $\pi_1(Y_1)^{\lift}$.
Let $\delta \in \pi_1(\Fr_{D_1}(Y_1))^{\lift}$ be an element in the preimage of $\overline{\delta}$.
Recall that we have a free splitting
\begin{equation*}
    \pi_1(Y_1) \cong \pi_1(Q_0) * \pi_1(Q_1).
\end{equation*}
Up to conjugacy, we may assume $\overline{\delta}$ lies in the subgroup $\pi_1(Q_0)$; this follows from our assumption in Section \ref{subsec:preliminaries} that the component $\del U_B \subseteq \del X$ projects into $(Q_0')^* \subseteq Y_1$.
Let $\overline{\delta}'$ be a conjugate of $\overline{\delta}$ by a nontrivial element of $\pi_1(Q_1)$, and let $\delta' \in \pi_1(\Fr_{D_1}(Y_1))^{\lift}$ be an element in the preimage of $\overline{\delta}'$.
We define $\Lambda \leq \pi_1(\Fr_{D_1}(Y_1))^{\lift}$ to be the subgroup generated by $\delta$ and $\delta'$.
See Figure \ref{fig:meridian-and-image} for an illustration of the curves $\overline{\gamma}$ and $\overline{\delta}$.

Now, we can prove that $\Lambda$ has the desired properties using the commutative diagram (\ref{eqn:lifting-diagram}).

\begin{lemma}\label{lem:free-product-downstairs}
    Let $\Lambda' \leq \LMod(Y)$ be the image of $\Lambda$, and let $\Lambda'' \leq \Out(\pi_1(Y))$ be the image of $\Lambda'$.
    Then $\Lambda''$ is a nontrivial free product of cyclic groups, and the surjection $\Lambda' \rightarrow \Lambda''$ has a finite kernel.
\end{lemma}
\begin{proof}
    Let $\alpha_{\delta}$ and $\alpha_{\delta'}$ denote the images of $\delta$ and $\delta'$ in $\LMod(Y)$, so $\Lambda'$ is generated by $\alpha_\delta$ and $\alpha_{\delta'}$,
    and $\Lambda''$ is generated by the images of $\alpha_{\delta}$ and $\alpha_{\delta'}$ in $\Out(\pi_1(Y))$.

    First, we can identify the group $\Lambda''$ more explicitly.
    Let $\Inn(\overline{\delta})$ and $\Inn(\overline{\delta}')$ denote the inner automorphisms of $\pi_1(Y_1)$ determined by $\overline{\delta}$ and $\overline{\delta}'$.
    We have a nontrivial free splitting
    \begin{equation*}
        \pi_1(Y) \cong \pi_1(Y_1) * \pi_1(Y_2).
    \end{equation*}
    Each element of $\Inn(\pi_1(Y_1))$ extends to an automorphism of $\pi_1(Y)$ which acts trivially on $\pi_1(Y_2)$, and these automorphisms are never inner (except for the identity), so we get an embedding $\Inn(\pi_1(Y_1)) \hookrightarrow \Out(\pi_1(Y))$.
    Now, $\alpha_{\delta}$ and $\alpha_{\delta'}$ are summand slides, and they act on $Y_1$ by $\calDS_{D_1}(\overline{\delta})$ and $\calDS_{D_1}(\overline{\delta}')$ respectively.
    Thus by Proposition \ref{prop:summand-slide-action-on-pi1}, we see that $\Lambda''$ is the image of $\la \Inn(\overline{\delta}), \Inn(\overline{\delta}')\ra$ under the emebdding $\Inn(\pi_1(Y_1)) \hookrightarrow \Out(\pi_1(Y))$.

    Now, we can complete the proof.
    Since $\pi_1(Y_1)$ is a free product, it has trivial center, and hence the map $\pi_1(Y_1) \rightarrow \Inn(\pi_1(Y_1))$ is an isomorphism.
    We therefore know that $\Lambda''$ is a free product since $\overline{\delta}$ and $\overline{\delta}'$ generate a free product in $\pi_1(Y_1)$.
    It remains just to check that the map $\Lambda' \rightarrow \Lambda''$ has a finite kernel.
    But this follows because the map $\Lambda \rightarrow \Lambda''$ has a finite kernel by Proposition \ref{prop:disk-slides-twists-and-inner-auts} (applied to  $M = Y_1$ and $\Delta = D_1$).
\end{proof}

\begin{lemma}\label{lem:finite-upstairs}
    The group $\Lambda$ maps onto a finite subgroup in $\Mod(M)$.
\end{lemma}
\begin{proof}
    From the diagram (\ref{eqn:lifting-diagram}), it's enough to show that $\Lambda$ maps to a finite subgroup of $\pi_1(\Fr_{\Delta_1}(M_1))$.
    Let $\widetilde{\delta}$ and $\widetilde{\delta}'$ be the images of $\delta$ and $\delta'$ in $\pi_1(\Fr_{\Delta_1}(X_1))$.
    Let $d = \v G \v$.
    By Proposition \ref{prop:disk-slides-twists-and-inner-auts}, we have an exact sequence
    \begin{equation*}
        (\Z/2\Z)^d \rightarrow \pi_1(\Fr_{\Delta_1}(X_1)) \rightarrow \pi_1(X_1)^{d} \rightarrow 1
    \end{equation*}
    (note that in the product $\pi_1(X_1)^{d}$, each factor has a different basepoint, which we omit for simplicity).
    By construction of the map $\pi_1(\Fr_{D_1}(Y_1))^{\lift} \rightarrow \pi_1(\Fr_{\Delta_1}(X_1))$, the image of $\widetilde{\delta}$ in $\pi_1(X_1)^d$ is a $d$-tuple comprised of the $G$-orbit of $\overline{\gamma}$, and in particular it is the $G$-orbit of a loop in the free homotopy class of $\mu$.
    Similarly, since $\overline{\delta}'$ is a conjugate of $\overline{\delta}$, the image of $\widetilde{\delta}'$ in $\pi_1(X_1)^d$ is a $d$-tuple comprised of the $G$-orbit of a loop in the free homotopy class of $\mu$.
    
    Proposition \ref{prop:disk-slides-twists-and-inner-auts} also gives us an exact sequence
    \begin{equation*}
        (\Z/2\Z)^d \rightarrow \pi_1(\Fr_{\Delta_1}(M_1)) \rightarrow \pi_1(M_1)^{d} \rightarrow 1.
    \end{equation*}
    Since $\mu$ is homotopically trivial in $M_1$, it follows from the discussion above that under the map $\pi_1(\Fr_{\Delta_1}(X_1)) \rightarrow \pi_1(\Fr_{\Delta_1}(M_1))$, the elements $\widetilde{\delta}$ and $\widetilde{\delta}'$ are sent to the kernel of the map $\pi_1(\Fr_{\Delta_1}(M_1)) \rightarrow \pi_1(M_1)^d$, and hence lie in the image of $(\Z/2\Z)^d$ in $\pi_1(\Fr_{\Delta_1}(M_1))$.
\end{proof}

\subsection{Completing the proof}\label{subsec:completing-pf-part-i}

Finally, we can complete the proof of part (i) of Theorem \ref{mainthm:ker-not-fg} as follows.

By Lemmas \ref{lem:free-product-downstairs} and \ref{lem:finite-upstairs}, the subgroup $\Lambda \leq \pi_1(\Fr_{D_1}(Y_1))^{\lift}$ maps to a subgroup $\Lambda' \leq \LMod(Y)$ which is a finite extension of a free product, and $\Lambda$ maps to a finite group in $\Mod(M)$.
Then $\Lambda'$ has a finite index subgroup $\Omega$ which maps to the identity in $\Mod(M)$.

We have a lifting map $\LDiff_\del(Y) \rightarrow \Diff_\del(X)$ that maps a diffeomorphism to the unique lift fixing $\del X$ (which is nonempty).
The image of this map is contained in the normalizer of $G$ in $\Diff_\del(X)$, which contains the centralizer $\Diff_\del(X)^G$ as a finite index subgroup.
Let $\LDiff_\del^G(Y)$ denote the preimage of $\Diff_\del(X)^G$.
Then $\pi_0(\LDiff_\del^G(Y))$ has finite index in $\LMod(Y)$.
Thus, up to passing to a further finite index subgroup of $\Omega$, we may assume that each diffeomorphism of $Y$ representing an element of $\Omega$  lifts to an equivariant diffeomorphism of $X$.
This implies that $\Omega$ embeds into $\pi_0(\Diff_{\del}(X)^G)$; indeed, if $f$ represents an element of $\Omega$ and $\widetilde{f}$ is its equivariant lift, then any equivariant isotopy $\widetilde{f} \simeq \id_{X}$ will descend to an isotopy $f \simeq \id_Y$. 

We have a natural inclusion $\Diff_\del(X) \hookrightarrow \Diff(M)$ by extending a diffeomorphism by the identity, and this restricts to an inclusion $\Diff_\del(X)^G \hookrightarrow \Diff(M)^G$.
We claim that the composition
\begin{equation*}
    \Omega \hookrightarrow \pi_0(\Diff_{\del}(X)^G) \rightarrow \pi_0(\Diff(M)^G)
\end{equation*}
has a finite kernel.
To see this, suppose $f \in \Diff_{\del}(Y)$ represents an element of $\Omega$ and $\widetilde{f} \in \Diff_{\del}(X)^G$ is its lift.
Let $\widetilde{f}' \in \Diff(M)^G$ be the natural extension of $\widetilde{f}$ (so $\widetilde{f}'$ represents the image of $[f]$ in $\pi_0(\Diff(M)^G)$).
Then $\widetilde{f}'$ restricts to a diffeomorphism of $M^\circ$, and thus descends to a diffeomorphism $f'$ of $M^\circ/G$ which extends $f$ (recall that $Y$ is a deformation retract of $M^\circ/G$).
Suppose now that $f$ represents an element in the kernel of the map $\Omega \rightarrow \pi_0(\Diff(M)^G)$, meaning there is an equivariant isotopy $\widetilde{f}' \simeq \id_M$.
This isotopy must preserve the set $R \subseteq M$ of points with nontrivial $G$-stabilizers, and thus this isotopy restricts to $M^\circ$ and descends to an isotopy $f' \simeq \id_{M^\circ/G}$.
Since $f'$ is an extension of $f$ and $Y$ is a deformation retract of $M^\circ/G$, this implies that $f$ is homotopic to $\id_Y$.
Then $f$ must represent an element of the kernel of $\Omega \rightarrow \Out(\pi_1(Y))$, which we know is finite by Lemma \ref{lem:free-product-downstairs}.

Thus $\Omega$ virtually embeds into the kernel of the natural map
\begin{equation*}
    \P_G:\pi_0(\Diff(M)^G) \rightarrow \pi_0(\Diff(M)),
\end{equation*}
which completes the proof.

\section{A Model of the Action, and an Infinite Cover}\label{sec:model-of-action}

We now begin towards the proof of part (ii) of Theorem \ref{mainthm:ker-not-fg}.
Fix a manifold $M$ and a group $G$ satisfying the assumptions of Theorem \ref{mainthm:ker-not-fg} with $k=3$.
By assumption, we have a codimension $2$ submanifold $B \subseteq M$ which is fixed by some element of $G$ and which represents a trivial homology class in $H_{n-2}(M;\Z/\ell\Z)$ for some $\ell > 1$.
In this section, we construct local and global models of the action of $G$ on $M$, and use them to construct an infinite sheeted cover $\pi:\widetilde{M} \rightarrow M$ that we will use in the proof of Theorem \ref{mainthm:ker-not-fg}.

In Section \ref{subsec:local-model}, we describe the local action of $G$ near a generic point of $B$.
In Section \ref{subsec:global-model}, we construct a global model of the action by decomposing $M$ according to the connected summands of $M^\circ/G$.
In Sections \ref{subsec:fundamental-group} and \ref{subsec:invariant-subgroup}, we describe $\pi_1(M)$ in terms of this global model.
In Section \ref{subsec:infinite-cover}, we construct the cover $\widetilde{M}$ using this description of $\pi_1(M)$, and in Section \ref{subsec:homology-of-cover}, we show that the homology of $\widetilde{M}$ is infinite dimensional.

\subsection{A local model of the action}\label{subsec:local-model}

We begin by describing the local behavior of the action of $G$ near a generic point of $B$.

\begin{lemma}\label{lem:local-model-of-action}
  There exists a point $b_0 \in B$ whose stabilizer $G_0 \coloneqq \Stab_G(b_0)$ is cyclic and fixes $B$ pointwise.
  Moreover, each point $b$ in the $G$-orbit of $b_0$ has a closed neighborhood $Z_b$ which is invariant under the stabilizer $G_b \coloneqq \Stab_G(b)$ and satisfies the following properties:
  \begin{enumerate}[label=(\roman*)]
    \item There is a diffeomorphism $Z_b \cong D^{n-2} \times D^2$ under which $G_b$ fixes the $D^{n-2}$-factor pointwise and acts on the $D^2$-factor by a $\frac{2\pi}{m}$-rotation, where $m = \v G_b \v$.
    \item For any two distinct points $b$ and $b'$ in the $G$-orbit of $b_0$, the neighborhoods $Z_b$ and $Z_{b'}$ are disjoint.
    \item For any $g \in G$, $g(Z_b) = Z_{g(b)}$.
  \end{enumerate}
\end{lemma}

Note that Lemma \ref{lem:local-model-of-action} requires the hypothesis that $G$ acts by orientation-preserving diffeomorphisms.
Indeed, the lemma fails if we let $M = S^3$, let $G = \Z_2 \times \Z_2$ where the two generators act by reflections in the $xy$- and $xz$-planes (viewing $S^3$ as $\R^3 \cup \{\infty\})$, and let $B$ be the $x$-axis, which is the fixed set of the element $g_0 = (1,1)$.

\begin{proof}[Proof of Lemma \ref{lem:local-model-of-action}]
    Given any $x \in M$, we let $G_x$ denote its stabilizer $\Stab_G(x)$.
    Assume $G$ fixes some Riemannian metric on $M$.
    Then for any $x \in M$, the stabilizer $G_x$ acts on the tangent space $T_xM$ as a finite subgroup of $\mathrm{SO}(n)$.
    Moreover, the exponential map $T_xM \rightarrow M$ is a $G_x$-equivariant local diffeomorphism at $0 \in T_xM$.
    The submanifold $B$ is a closed totally geodesic submanifold of $M$, since it is a component of $\Fix(g_0)$ (see \cite[Thm~5.1]{kobayashi}).
  
    Let $H \leq G$ be the subgroup of elements that fix $B$ pointwise.
    Let $T$ be an $H$-invariant closed tubular neighborhood of $B$ (see \cite[Thm~VI.2.2]{bredon} for the existence of such a neighborhood).
    Since $B$ has codimension 2, each point $x \in B$ has a closed neighborhood $Z$ with a diffeomorphism $Z \cong D^{n-2} \times D^2$ taking $Z \cap B$ to $D^{n-2} \times \{0\}$, and taking each fiber of $T$ to a fiber $\{y\} \times D^2$.
  
    For $x \in B$, the action of $H$ on $T_xM$ fixes a codimension $2$ subspace $W \subseteq T_xM$, and the exponential map yields an isomorphism between the action of $H$ on the $D^2$-fibers of $T$ and the linear action of $H$ on a small $2$-disk $D \subseteq W^\perp$.
    Note that $H$ must act faithfully on $T_xM$; indeed, if $h \in H$ fixes $T_xM$ pointwise, then the fixed set $\Fix(h) \subseteq M$ must be a closed submanifold of codimension 0, and thus $\Fix(h) = M$.
    Since $H$ acts trivially on $W$, it must be that $H$ acts faithfully on $W^\perp$.
    Since $G$ is orientation-preserving, we know that no element of $H$ can act on $W^\perp$ by a reflection, and thus $H$ acts on $W^\perp$ by a finite subgroup of $\SO(2)$.
    Thus $H$ must be cyclic and act on the $D^2$-factor of $Z$ by rotations.
    
    Now, we claim that there exists a point $b_0 \in B$ whose stabilizer $G_0 \coloneqq G_{b_0}$ is equal to $H$.
    Assuming this claim, we can prove the lemma.
    Indeed, let $Z_{b_0}$ be the neighborhood $Z$ described above, and for $b = gb_0$, we define $Z_{b} \coloneqq g(Z_{b_0})$.
    This is well-defined since $Z_{b_0}$ is $G_0$-invariant.
    Then the fact that $G_0$ is cyclic and statement (i) follow from the discussion above, statement (ii) follows by choosing $Z_{b_0}$ sufficiently small, and statement (iii) is automatic.
  
    Thus, it remains to show that there exists a point $b_0 \in B$ such that $G_{b_0} = H$.
    Choose any $x \in B$.
    Observe that $H \leq G_x$ automatically, so suppose that there exists $g \in G_x \setminus H$.
    Let $W \subseteq T_xM$ be the codimension $2$ subspace spanned by $\mathrm{exp}^{-1}(V)$, where $V \subseteq B$ is a small neighborhood of $x$.
    Note that $g$ cannot fix $W$ pointwise; otherwise $g$ would fix an open subset of $B$ pointwise, which would imply that $g$ fixes $B$ pointwise,\footnote{
        To see this, observe that $\{x \in B \mid \text{$g$ fixes an open neighborhood of $x$ in $B$}\}$ is both open and closed in $B$, and so if this set is nonempty, then it is all of $B$.
    } contradicting that $g \not\in H$.
    Thus if we let $Q \subseteq T_xM$ denote the fixed set of $g$, the intersection $W \cap Q$ has positive codimension in $W$.
    Thus, there is an open subset $U \subseteq B$ such that $g \not\in G_y$ for any $y \in U$.
    We can repeat the argument inductively; choosing $y \in U$, if there exists $g' \in G_y \setminus H$, we can pass to a smaller neighborhood where $g'$ does not stabilize any point.
    Since $G$ is finite, we must eventually find a point $b_0 \in B$ such that $G_{b_0} \setminus H$ is empty.
  \end{proof}

\subsection{A global model of the action}\label{subsec:global-model}
Next, we can describe how $G$ acts on the connected summands of $M$.
This model will be slightly different from the one in Section \ref{subsec:preliminaries}.
Recall that we let $M^\circ$ denote the set of points whose $G$-stabilizer is trivial.
Let $N = M^\circ/G$.
By assumption, $N \cong Q_0 \# Q_1 \# Q_2$ where $\pi_1(Q_i)$ is nontrivial for each $i$.

We begin by constructing a model of $N$.
Let $D_1^+$ and $D_2^+$ be embedded $n$-disks in $Q_0$.
Similarly, for $i \in \{1,2\}$, let $D_i^-$ be an embedded $n$-disk in $Q_i$.
Then, we view $N$ as the manifold obtained by gluing each $D_i^+$ to $D_i^-$ (see Section \ref{subsec:summand-slides} for our gluing conventions).
We let $Q_0^*$ denote $Q_0$ minus the centers of each $D_i^+$, and for $i \in \{1,2\}$ we let $Q_i^*$ denote $Q_i$ minus the center of $D_i^-$.
For $i \in \{0,1,2\}$ we identify $Q_i^*$ with its image in $N$.
We let $(D_i^+)^* = D_i^+ \cap Q_0^*$ and $(D_i^-)^* = D_i^- \cap Q_i^*$.

Next, we can construct a model of $M$.
For each $i \in \{0,1,2\}$, we let $(P_i^\circ)^* \subseteq M^\circ$ denote the preimage of $Q_i^*$ in $M^\circ$, and we let $P_i^\circ$ be the (possibly disconnected) manifold obtained by filling in the punctures.
Similarly, we let $(\Delta_i^+)^* \subseteq (P_0^\circ)^*$ denote the preimage of $(D_i^+)^*$ and let $(\Delta_i^-)^* \subseteq (P_i^\circ)^*$ denote the preimage of $(D_i^-)^*$, and we let $\Delta_i^+ \subseteq P_0^\circ$ and $\Delta_i^- \subseteq P_i^\circ$ denote the multidisks obtained by filling in the punctures.
Then we can view $M^\circ$ as the manifold obtained by gluing $P_0^\circ$ to $P_i^\circ$ along the multidisks $\Delta_i^\pm$ for each $i \in \{1,2\}$.
Finally, the inclusion $M^\circ \hookrightarrow M$ induces a way to fill in the ends of $P_i^\circ$ (resp.\ $(P_i^\circ)^*$) to obtain a (possibly disconnected) manifold $P_i$ (resp.\ $P_i^*$), so $M$ is obtained by gluing $P_0$ to $P_i$ for $i \in \{1,2\}$ along the multidisks $\Delta_i^\pm$.
Note that the action of $G$ on $M$ induces an action of $G$ on each $P_i$.

Note that for each $i \in \{1,2\}$, the group $G$ acts freely and transitively on the components of $\Delta_i^+$ and $\Delta_i^-$.
We can therefore label the components of $\Delta_i^{\pm}$ by the elements of $G$, i.e.\ $\Delta_i^\pm = (D_{i,g}^\pm)_{g \in G}$.
We choose these labels so that $gD_{i,h}^\pm = D_{i,gh}^\pm$, and so that the disk $D_{i,g}^+$ is glued to $D_{i,g}^-$.
We call the disks $D_{i,g}^\pm$ the \emph{gluing disks} of $M$.
We define $\Sigma_{i,g} \coloneqq \del D_{i,g}^+$, so $\Sigma_{i,g}$ is an $(n-1)$-sphere in $P_0$.

Fix a component $M_0$ of $P_0$, and let $M_0^* \subseteq M$ denote the corresponding component of $P_0^*$.
We assume without loss of generality that $B$ lies on $M_0$, and that $M_0$ contains the gluing disks $D_{i,\id}^+$ for $i \in \{1,2\}$.
We fix a basepoint $b_0 \in B$ as in Lemma \ref{lem:local-model-of-action}, and let $G_0 \leq G$ denote the stabilizer of $b_0$.

\begin{figure}
    \centering
    \includegraphics[scale=.55]{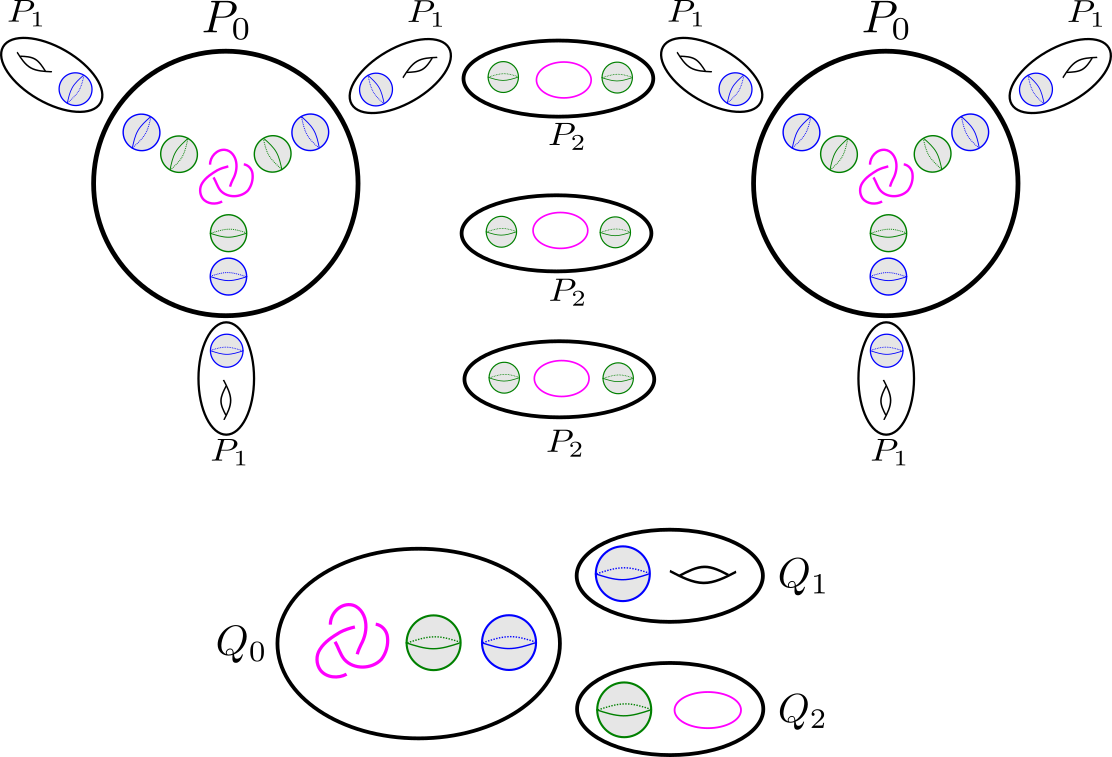}
    \caption{A schematic illustration of $M$ and $N$.  Here the trefoils and unknots (in magenta) represent the points with nontrivial stabilizers and their images in $M/G$, so $N$ and $M^\circ$ are the complement of these curves.}
    \label{fig:global-model}
\end{figure}

See Figure \ref{fig:global-model} for a schematic illustration of $M$.

\subsection{The fundamental group}\label{subsec:fundamental-group}

Our next goal is to describe the fundamental group of $M$ in terms of the pieces $P_i$.
It will be convenient to describe $\pi_1(M,b_0)$ in terms of a graph of groups; we will refer to \cite[\S I.5]{serre} and \cite{higgins} for background material on graphs of groups.
We adapt the convention that a \emph{graph} $\calG$ consists of 
\begin{itemize}
    \item a \emph{vertex set} $V(\calG)$, and
    \item an \emph{edge set} $E(\calG) \subseteq V(\calG) \times V(\calG)$ such that for each edge $e = (v_1,v_2) \in E(\calG)$, the set $E(\calG)$ contains the unique \emph{reverse edge} $\overline{e} = (v_2,v_1)$.
\end{itemize}
The \emph{valence} of a vertex $v$ is the number of edges of the form $(v, v')$ for some $v' \in V(\calG)$.
A \emph{graph of groups} is a graph equipped with a group $\calG_v$ for each vertex $v \in V(\calG)$ (typically one allows edge groups as well, but this will not be relevant for us).

For $i \in \{0,1,2\}$, let $\calC_i$ denote the set of components of $P_i$, and let $\calC = \calC_0 \cup \calC_1 \cup \calC_2$.
For each $C \in \calC$, choose a point $p_C \in C$ away from any gluing disks; we assume that $p_{M_0} = b_0$.
For each $i \in \{1,2\}$ and $g \in G$, we also fix a point $p_{i,g} \in (D_{i,g}^+)^* \subseteq M$ (we can equivalently view $p_{i,g}$ as lying in $(D_{i,g}^-)^*$).
We define $\calG$ to be the graph whose vertex set $V(\calG)$ is the set of all the points $p_C$ and $p_{i,g}$, and whose edge set $E(\calG)$ consists of pairs $(p_C,p_{i,g})$ and $(p_{i,g}, p_C)$ where $p_{i,g} \in C$.
Note that each vertex $p_{i,g}$ has valence 2.
We equip $\calG$ with the structure of a graph of groups where the vertex group at $p_C$ is $\calG_C = \pi_1(C, p_C)$ and whose vertex group at $p_{i,g}$ is trivial. 

Following Higgins \cite{higgins}, we can define the fundamental groupoid $\Pi(\calG)$ of $\calG$.
Namely, $\Pi(\calG)$ is the groupoid whose objects are precisely the elements of $V(\calG)$, and whose morphisms are generated by the edges of $\calG$ and the vertex groups of $\calG$, subject to the relation that for each $e \in \calG$, the morphism associated to $\overline{e}$ is the inverse of the morphism associated to $e$.
We define the fundamental group $\pi_1(\calG, M_0)$ to be the automorphism group of the object $b_0 = p_{M_0}$ of $\Pi(\calG)$.

There is a more concrete description of $\pi_1(\calG, M_0)$ following Serre \cite{serre}.
Given an edge $e = (v, v')$, we write $i(e) = v$ and $t(e) = v'$.
We define a \emph{path} in $\calG$ to be a pair $(\widehat{e}, \widehat{r})$, where either $\widehat{e} = (e_1, \ldots, e_m)$ is a tuple of edges with $t(e_j) = i(e_{j+1})$ and $\widehat{r} = (r_0, \ldots, r_m)$ is a tuple with $r_0 \in \calG_{i(e_1)}$ and $r_j \in \calG_{t(e_j)}$ for $j \geq 1$, or $\widehat{e} = \varnothing$ and $\widehat{r} = (r_0)$ where $r_0 \in \calG_v$ for some $v \in V(\calG)$.
A \emph{loop} in $\calG$ based at a vertex $v$ is a path $(\widehat{e}, \widehat{r})$ such that either $\widehat{e} = \varnothing$ and $r_0 \in \calG_v$, or $i(e_1) = t(e_m) = v$.
Given a path $(\widehat{e}, \widehat{r})$, we get an associated morphism in $\Pi(\calG)$ given by
\begin{equation*}
    \v (\widehat{e},\widehat{r}) \v \coloneqq r_0e_1r_1 \cdots e_mr_m. 
\end{equation*}
Then $\pi_1(\calG, M_0)$ is the group generated by elements of the form $\v (\widehat{e}, \widehat{r}) \v$, where $(\widehat{e}, \widehat{r})$ is a loop based at the vertex $M_0$.
Note that we can naturally view $\pi_1(M_0,b_0)$ as the subgroup of $\pi_1(\calG, M_0)$ generated by paths with $\widehat{e} = \varnothing$.

There is yet another characterization which will be useful for us.
Let $T$ be a spanning tree of the graph underlying $\calG$.
Let $\pi(\calG)$ be the group obtained by identifying all the objects of $\Pi(\calG)$ (equivalently, $\pi(\calG)$ is the free product of the groups $\pi_1(C,p_C)$ with the free group on $E(\calG)$, subject to the relations $\overline{e} = e^{-1}$).
Then we can define $\pi_1(\calG,T)$ to be the group obtained by starting with the group $\pi(\calG)$ and adding the relation $e = 1$ for each edge in $T$.
By \cite[Prop~I.20]{serre}, the natural map $\pi_1(\calG,M_0) \rightarrow \pi_1(\calG,T)$ is an isomorphism for any spanning tree $T$.

Now, we can describe $\pi_1(M,b_0)$ in terms of $\calG$.
Given a component $C$ of $P_i$, we let $C^* \subseteq M$ denote the corresponding component of $P_i^*$.
For each edge $e = (p_C,p_{i,g})$ of $\calG$, choose an oriented arc $a_e$ in $C^*$ joining $p_C$ to $p_{i,g}$; we define $a_{\overline{e}}$ by reversing the arc $a_e$.
Then, given a path $(\widehat{e},\widehat{r})$ in $\calG$, we can construct a path $\gamma_{(\widehat{e},\widehat{r})}$ in $M$ by concatenating representatives of the elements $r_j$ with the arcs $a_{e_j}$.
The path $\gamma_{(\widehat{e},\widehat{r})}$ is well-defined up to homotopy rel the points $p_C$ and $p_{i,g}$.
Then we have the following.

\begin{lemma}\label{lem:pi1-computation}
    There is an isomorphism
    \begin{equation*}
        \pi_1(\calG, M_0) \cong \pi_1(M,b_0)
    \end{equation*}
    which sends the element $\v (\widehat{e},\widehat{r})\v$ to the loop $\gamma_{(\widehat{e},\widehat{r})}$ for any loop $(\widehat{e},\widehat{r})$ in $\calG$ based at $M_0$.
\end{lemma}
\begin{proof}
    Let $\mathbf{p} \subseteq M$ be the set of points $p_C$ for $C \in \calC$ and $p_{i,g}$ for $i \in \{1,2\}$ and $g \in G$.
    Let $\Pi(M, \mathbf{p})$ be the groupoid whose objects are the points in $\mathbf{p}$ and whose morphisms are homotopy classes of paths with endpoints in $\mathbf{p}$.
    We can apply the groupoid version of van Kampen's theorem (see \cite{brown-razak}) to show that there is an isomorphism $\Pi(M,\mathbf{p}) \cong \Pi(\calG)$ which induces the desired isomorphism on fundamental groups.
    Namely, for each $C \in \calC$, we have a natural map $\Pi(C^*, \mathbf{p} \cap C^*) \rightarrow \Pi(\calG)$ which maps the fundamental group $\pi_1(C^*, p_C)$ to the vertex group at $p_C$, and which maps each arc $a_e$ for $e = (C, p_{i,g})$ to the morphism corresponding to $e$.
    Each $2$-fold intersection of the sets $C^*$ for $C \in \calC$ is a disjoint union of simply connected subsets, and each $3$-fold intersection is empty.
    Van Kampen's theorem then implies that the maps $\Pi(C^*, \mathbf{p} \cap C^*) \rightarrow \Pi(\calG)$ assemble into a map $\Pi(M, \mathbf{p}) \rightarrow \Pi(\calG)$.
    This map has a natural inverse, and is therefore an isomorphism.
\end{proof}

\subsection{The invariant subgroup of the fundamental group}\label{subsec:invariant-subgroup}

Recall that we let $G_0 \leq G$ denote the stabilizer of the basepoint $b_0$.
We have a well-defined action of $G_0$ on $\pi_1(M, b_0)$.
Let $\pi_1(M,b_0)^{G_0} \leq \pi_1(M,b_0)$ denote the fixed points of this action.
The following lemma constrains the subgroup $\pi_1(M,b_0)^{G_0}$.
We will use this result in Section \ref{sec:obstruction-map} in the construction of our obstruction map (see Lemma \ref{lem:lifting-loops-of-submanifolds}).

\begin{lemma}\label{lem:pi1-fixed-points}
    The subgroup $\pi_1(M,b_0)^{G_0}$ is contained in the subgroup $\pi_1(M_0,b_0) \leq \pi_1(M,b_0)$.
\end{lemma}
\begin{proof}
    Let $\widehat{M} \rightarrow M$ denote the universal cover of $M$, and fix a base point $\widehat{b}_0 \in \widehat{M}$ above $b_0$.
    Let $\widehat{P}_i^* \subseteq \widehat{M}$ denote the preimage of $P_i^*$, and let $\widehat{P}_i$ be obtained by filling in the punctures.
    We call each component of each $\widehat{P}_i$ a \emph{sheet} of $\widehat{M}$.
    Let $\widehat{M}_0$ denote the sheet containing $\widehat{b}_0$ (so $\widehat{M}_0$ is a component of $\widehat{P}_0$, and is diffeomorphic to the universal cover of $M_0$).

    We will prove the lemma by contrapositive.
    Suppose $\gamma \in \pi_1(M, b_0)$ is an element that is not contained in the subgroup $\pi_1(M_0, b_0) \leq \pi_1(M,b_0)$, and let $g \in G_0$ be a generator (recall that $G_0$ is cyclic).
    Our goal is then to show that $g(\gamma) \neq \gamma$.
    We will show that a lift of ${\gamma}$ to $\widehat{M}$ travels through a unique sequence of sheets up to homotopy, and then show that a lift of $g$ cannot preserve this sequence of sheets.

    Let $\widehat{\gamma}$ denote the lift of $\gamma$ to $\widehat{M}$ based at $\widehat{b}_0$.
    Since $\widehat{M}$ is simply connected, we may assume up to homotopy rel endpoints that for any sheet $W$ of $\widehat{M}$, the intersection $\widehat{\gamma} \cap W$ is connected, i.e.\ $\widehat{\gamma}$ travels through each sheet at most once.
    We claim that $\widehat{\gamma}$ must end on a different sheet than $\widehat{M}_0$.
    Indeed, suppose for the sake of contradiction that $\widehat{\gamma}$ ends at a point on $\widehat{M}_0$.
    Then we can assume that $\widehat{\gamma}$ is entirely contained in the sheet $\widehat{M}_0$.
    Then $\widehat{\gamma}$ projects to a loop in $M$ entirely contained within $M_0$, contradicting that $\gamma$ lies outside of the subgroup $\pi_1(M_0,b_0)$.
    Thus, $\widehat{\gamma}$ travels through a sequence of distinct sheets $\widehat{M}_0 = W_0, W_1, \ldots, W_m$ with $m > 0$, where $W_i$ is adjacent to $W_{i+1}$ and $W_m$ is a component of $\widehat{P}_0$ (since the endpoint of $\widehat{\gamma}$ must live above $b_0$, which lies in $M_0$).

    Next, we claim that this sequence of sheets is unique; that is, if $\widehat{\gamma}$ is homotopic rel endpoints to an arc $\widehat{\gamma}'$ traveling through a sequence of distinct sheets $\widehat{M}_0 = W_0, W_1', \ldots, W_{m'}'$, then $m' = m$ and $W_j = W_j'$ for all $j$.
    Suppose for contradiction that for some $k < m$, we have that $W_j = W_j'$ for $j \leq k$ but $W_{k+1} \neq W_{k+1}'$.
    Then the curves $\widehat{\gamma}$ and $\widehat{\gamma}'$ leave through two different disks on the sheet $W_k = W_k'$, say $D$ and $D'$.
    Let $\Sigma = \del D$ and $\Sigma' = \del D'$.
    Equip $\Sigma$ and $\Sigma'$ with orientations such that $\widehat{\gamma}$ and $\widehat{\gamma}'$ intersect them positively, and let $\overline{\Sigma'}$ denote $\Sigma'$ equipped with the opposite orientation.
    Let $\delta$ be the loop obtained by concatenating $\widehat{\gamma}$ with the reverse of $\widehat{\gamma}'$.
    Then the loop $\delta$ and the submanifold $\Sigma \sqcup \overline{\Sigma'}$ have a positive signed intersection number.
    Since this is a homotopy invariant of $\delta$, it follows that the loop $\delta$ is not homotopic to a point, contradicting that $\widehat{M}$ is simply connected.

    Now, lift $g$ to a diffeomorphism $\widehat{g}:\widehat{M} \rightarrow \widehat{M}$ that fixes the point $\widehat{b}_0$.
    Since $g$ preserves $M_0^* \subseteq M$ (as it fixes $b_0$), it follows that $\widehat{g}$ preserves the sheet $\widehat{M}_0$.
    However, since $g$ does not fix any of the gluing disks in $M_0$, it follows that $\widehat{g}$ does not preserve any of the sheets adjacent to $\widehat{M}_0$.
    Thus, $\widehat{g}(\widehat{\gamma})$ travels through a different sequence of sheets than $\widehat{\gamma}$.
    By the uniqueness statement above, this means that $\widehat{g}(\widehat{\gamma})$ is not homotopic rel endpoints to $\widehat{\gamma}$, and thus we can conclude that $g(\gamma) \neq \gamma$.
\end{proof}

\subsection{An infinite cover}\label{subsec:infinite-cover}

Next, we will define the cover $\widetilde{M}$ that we will use to construct our obstruction map.
Let $\calG$ denote the graph of groups defined above, and let $\overline{\calG}$ denote the graph of groups obtained from $\calG$ by replacing the vertex group $\calG_{M_0} = \pi_1(M_0,b_0)$ with the trivial group.
By Lemma \ref{lem:pi1-computation}, the natural map $\Pi(\calG) \rightarrow \Pi(\overline{\calG})$ induces a surjection
\begin{equation*}
    \pi_1(M,b_0) \rightarrow \pi_1(\overline{\calG},M_0).
\end{equation*}
The kernel of this map is the normal closure $\la \la \pi_1(M_0,b_0) \ra \ra \leq \pi_1(M,b_0)$ (this is easily seen using the isomorphisms $\pi_1(\calG,M_0) \cong \pi_1(\calG,T)$ and $\pi_1(\overline{\calG},M_0) \cong \pi_1(\overline{\calG}, T)$, where $T$ is a spanning tree of the graphs underlying $\calG$ and $\overline{\calG}$).
We define $\pi:\widetilde{M} \rightarrow M$ to be cover associated to the kernel of this map.
We fix a basepoint $\widetilde{b}_0 \in \widetilde{M}$ above $b_0 \in M$, so $\pi_1(\widetilde{M},\widetilde{b}_0) \cong \la \la \pi_1(M_0,b_0) \ra \ra$, and the deck group of $\widetilde{M}$ is the group $\Lambda \coloneqq \pi_1(\overline{\calG}, M_0)$.

Observe that the cover $\widetilde{M}$ is trivial over $M_0^*$, i.e.\ each component of $\pi^{-1}(M_0^*)$ projects diffeomorphically onto $M_0^*$.
On the other hand, if $C^*$ is a component of $P_1^*$ or $P_2^*$, or if $C^*$ is a component of $P_0^* \setminus M_0^*$, then each component of $\pi^{-1}(C^*)$ is the universal cover of $C^*$.
Note that it is possible that $\pi_1(M_0,b_0)$ is trivial, in which case $\overline{\calG} = \calG$ and $\widetilde{M}$ is simply the universal cover of $M$.

The cover $\widetilde{M}$ will be infinite-sheeted; we can construct an infinite subgroup of the deck group $\Lambda = \pi_1(\overline{\calG},M_0)$ as follows.
Fix $i \in \{1,2\}$. 
Let $\overline{\calG}_i'$ be the induced subgraph of $\overline{\calG}$ whose vertices are the points $p_C$ for $C \in \calC_0 \cup \calC_i$ together with the points $p_{i,g}$ lying on each $C \in \calC_0 \cup \calC_i$, and we define $\overline{\calG}_i$ to be the component of $\overline{\calG}_i'$ containing the vertex $p_{M_0} = b_0$ (we equip $\overline{\calG}_i$ with the same vertex groups as $\overline{\calG}$).
We define $\Lambda_i \coloneqq \pi_1(\overline{\calG}_i, M_0)$, and observe that $\Lambda_i$ is naturally a subgroup of $\Lambda$.

\begin{lemma}\label{lem:infinite-subgroup}
    At least one of the following is true:
    \begin{itemize}
        \item $\overline{\calG}_i$ has at least two nontrivial vertex groups $\calG_C$ for $C \in \calC_i$, or
        \item $\overline{\calG}_i$ contains a cycle.
    \end{itemize}
    In particular, the subgroup $\Lambda_i \leq \Lambda$ is infinite.
\end{lemma}
\begin{proof}
    Observe first that each component of $P_i$ contains the same number of gluing disks.
    Thus, we analyze two cases: the case that each component has a single gluing disk, and the case that each component has at least two gluing disks.
  
    Suppose first that each component of $P_i$ has a single gluing disk.
    In this case, we claim that $\overline{\calG}_i$ contains at least two nontrivial vertex groups $\calG_C$ for $C \in \calC_i$.
    First, observe that since each component has a single gluing disk, the cover $P_i^\circ \rightarrow Q_i$ must be trivial.
    This means that $G$ acts freely on the components of $P_i^\circ$ and each component is diffeomorphic to $Q_i$.
    Then $G$ must act freely on $P_i$, so $P_i = P_i^\circ$.
    Thus each component of $P_i$ is not simply connected, since $Q_i$ is not simply connected.
    To complete the proof, we need to check that $\overline{\calG}_i$ contains at least two vertices of the form $p_C$ for $C \in \calC_i$.
    But this follows since each component of $P_0$ contains at least two of the gluing disks $\Delta_{i,g}^+$ for $g \in G$, since the basepoint $b_0 \in M_0$ has a nontrivial $G$-stabilizer. 
    In particular, we see that $\Lambda_i$ is infinite since it contains a free product $\pi_1(C) * \pi_1(C')$ for two distinct components $C,C' \in P_i$.
  
    Otherwise, suppose each component of $P_i$ has at least two gluing disks.
    In this case, we claim that $\overline{\calG}_i$ contains a cycle; the fact that $\Lambda_i$ is infinite then follows from the isomorphism $\pi_1(\overline{\calG}_i, M_0) \cong \pi_1(\overline{\calG}_i, T)$ for any spanning tree $T$.
    As in the previous case, we observe that each component of $P_0$ also contains at least two of the gluing disks $\Delta_{i,g}^+$ for $g \in G$.
    Therefore, since each component of $P_0$ and $P_i$ contains at least two gluing disks, each vertex of $\overline{\calG}_i$ must have valence at least 2.
    Thus $\overline{\calG}_i$ must contain a cycle, as any finite tree contains a vertex of valence $0$ or $1$.
\end{proof}

Observe that the cover $\widetilde{M}$ inherits a natural connected sum decomposition from $M$.
Similar to $M$, we can record this connected sum decomposition using a graph.
Let $\scrG$ denote the graph underlying $\calG$.
Using the arcs $a_e$ constructed in Section \ref{subsec:fundamental-group}, we can view $\scrG$ as embedded in $M$.
Let $\widetilde{\scrG}$ denote the preimage of $\scrG$ in $\widetilde{M}$, and observe that $\widetilde{\scrG}$ naturally inherits the structure of a graph.
We record the following observation, which we will use in the proof of Proposition \ref{prop:infinite-lin-ind-subset}.

\begin{lemma}\label{lem:cpt-graph-a-tree}
    The graph $\widetilde{\scrG}$ is a tree.
\end{lemma}
\begin{proof}
    We have the following commutative square:
    \begin{equation*}
        \begin{tikzcd}
            {\pi_1(\widetilde{\scrG})} & {\pi_1(\widetilde{M})} \\
            {\pi_1(\scrG)} & {\pi_1(M)}
            \arrow[from=1-1, to=1-2]
            \arrow[from=1-1, to=2-1]
            \arrow[from=1-2, to=2-2]
            \arrow[from=2-1, to=2-2]
        \end{tikzcd}
    \end{equation*}

    Observe first that every map in this square is injective.
    Indeed, the injectivity of the vertical maps follows because they are induced by covering maps, the injectivity of the bottom map can be deduced from the isomorphism $\pi_1(\calG, M_0) \cong \pi_1(\calG, T)$ described in Section \ref{subsec:fundamental-group}, and the injectivity of the top map follows from the injectivity of the other three.

    On the other hand, the images of the maps $\pi_1(\scrG) \rightarrow \pi_1(M)$ and $\pi_1(\widetilde{M}) \rightarrow \pi_1(M)$ intersect trivially.
    This follows from the fact that the map $\pi_1(\scrG) \rightarrow \pi_1(\overline{\calG}, M_0)$ is injective (again via the isomorphism $\pi_1(\overline{\calG}, M_0) \cong \pi_1(\overline{\calG}, T)$).
    Thus it must be that $\pi_1(\widetilde{\scrG})$ is trivial.
\end{proof}



\subsection{The homology of the cover}\label{subsec:homology-of-cover}

The key property of the cover $\widetilde{M}$ is that the homology group $H_{n-1}(\widetilde{M};\Z/\ell\Z)$ is not finitely generated.
In fact, we will describe a non-finitely generated subgroup of $H_{n-1}(\widetilde{M};\Z/\ell\Z)$.
In order to do so, we establish some notation.

Let $G_{M_0} \leq G$ denote the stabilizer of $M_0$ (with respect to the action of $G$ on the components of $P_0$).
Then for each $i \in \{1,2\}$, since $M_0$ contains the gluing disk $D_{i,\id}^+$, it follows that the disks $D_{i,g}^+$ for $g \in G_{M_0}$ are precisely the gluing disks contained in $M_0$.
Note that $G_0 \leq G_{M_0}$.

Let $W_{\id}$ be the component of $\pi^{-1}(M_0^*) \subseteq \widetilde{M}$ containing the base point $\widetilde{b}_0$.
Then, for $u \in \Lambda$, let $W_u$ be the image of $W_{\id}$ under $u$.
We will call each $W_u$ a \emph{sheet} of $\widetilde{M}$.
For each $u \in \Lambda$, the cover $\pi:\widetilde{M} \rightarrow M$ restricts to a diffeomorphism $W_u \rightarrow M_0^*$.
Recall that for $g \in G_{M_0}$, we let $\Sigma_{i,g} \subseteq M_0$ denote the boundary sphere of the disk $D_{i,g}^+$.
On $W_u$, let $\Sigma_{i,g,u}$ be the copy of $\Sigma_{i,g}$ for $i \in \{1,2\}$ and $g \in G_{M_0}$.

Now, we can construct an infinite linearly independent subset of $H_{n-1}(\widetilde{M};\Z/\ell\Z)$ indexed by the subgroup $\Lambda_2$ defined above.


\begin{proposition}\label{prop:infinite-lin-ind-subset}
    For any two distinct pairs $(g,u),(g',u') \in G_{M_0} \times \Lambda_2$, the spheres $\Sigma_{1,g,u}$ and $\Sigma_{1,g',u'}$ represent nonequal classes in $H_{n-1}(\widetilde{M};\Z/\ell\Z)$.
    Moreover, the set 
    \begin{equation*}
        \mathcal{S} = \{[\Sigma_{1,g,u}] \in H_{n-1}(\widetilde{M};\Z/\ell\Z) \mid g \in G_{M_0}, u \in \Lambda_2\}
    \end{equation*}
    is an infinite linearly independent subset of $H_{n-1}(\widetilde{M};\Z/\ell\Z)$.
\end{proposition}

To prove Proposition \ref{prop:infinite-lin-ind-subset}, we need the following preliminary lemma.


  
  

\begin{lemma}\label{lem:separating-spheres}
    For any $i \in \{1,2\}$, $g_0 \in G_{M_0}$, and $u_0 \in \Lambda$, the complement $\widetilde{M} \setminus \Sigma_{i,g_0,u_0}$ has two components, and each component is noncompact.
    Moreover, if $u_0 \in \Lambda_j$ where $j \neq i$, then each sheet $W_{u} \subseteq \widetilde{M}$ for $u \in \Lambda_j \setminus \{u_0\}$ lies in a single component of $\widetilde{M} \setminus \Sigma_{i,g_0,u_0}$, namely the component containing the point $u_0(\widetilde{b}_0) \in W_{u_0}$.
\end{lemma}
\begin{proof}
    We begin by showing that $\widetilde{M} \setminus \Sigma_{i,g_0,u_0}$ is disconnected (and hence has two components).
    If $\widetilde{M} \setminus \Sigma_{i,g_0,u_0}$ were connected, then by cutting $\widetilde{M}$ along $\Sigma_{i,g_0,u_0}$, we could construct a simple closed curve that intersects $\Sigma_{i,g_0,u_0}$ once.
    Since mod $2$ intersection numbers are a homotopy invariant, it's therefore enough to show that every loop in $\widetilde{M}$ has an even intersection number with $\Sigma_{i,g_0,u_0}$.

    Let $\gamma$ be a loop in $\widetilde{M}$.
    Up to homotopy, we may assume that $\gamma$ is based at $\widetilde{b}_0$.
    Let $\overline{\gamma}$ be the image of $\gamma$ in $M$, so $\overline{\gamma}$ is a loop based at $b_0$.
    Since $\pi_1(\widetilde{M}, \widetilde{b}_0)$ is the normal closure $\la \la \pi_1(M_0,b_0) \ra \ra$, we know that that $\overline{\gamma}$ homotopic to a concatenation of loops $\overline{\gamma}_1, \ldots, \overline{\gamma}_m$ where each $\overline{\gamma}_j$ is freely homotopic into the subspace $M_0^* \subseteq M$.
    This means that $\gamma$ is homotopic to a concatenation of loops $\gamma_1, \ldots, \gamma_m$ where each $\gamma_j$ is freely homotopic into $W_{u_j}$ for some $u_j \in \Lambda$.
    This implies that the intersection number of $\gamma$ with any sphere of the form $\Sigma_{i,g,u}$ for $(g,u) \in G \times \Lambda$ must be even as desired.

    Now, we can show that each component is noncompact.
    First, we claim that each component $\widetilde{M} \setminus \Sigma_{i,g_0,u_0}$ contains at least one sheet $W_{u}$ for some $u \in \Lambda$.
    To prove this, it's enough to find a loop $\gamma$ in $M$ based at $b_0$ which lifts to an arc in $\widetilde{M}$ that passes through $\Sigma_{i,g_0,u_0}$ once.
    Let $C$ be the component of $P_i$ containing $D_{i,g_0}^-$.
    If $C$ is not simply connected, then we can construct $\gamma$ by conjugating a nontrivial element of $\pi_1(C,p_C)$ with a path from $b_0$ to $p_C$.
    Otherwise, assume each component of $P_i$ is simply connected.
    Then it's enough to show that the sphere $\Sigma_{i,g_0}$ is non-separating in $M$, since then we can take $\gamma$ to be a loop that intersects $\Sigma_{i,g_0}$ once.
    To prove this, suppose for contradiction that $\Sigma_{i,g_0}$ separates $M$.
    Then for any $g \in G$, the sphere $\Sigma_{i,gg_0}$ is also separating.
    This implies that every edge of the graph $\overline{\calG}_i$ is separating, which implies that the graph underlying $\overline{\calG}_i$ is a tree, contradicting Lemma \ref{lem:infinite-subgroup}.

    Now, we can show that each component of $\widetilde{M} \setminus \Sigma_{i,g_0,u_0}$ contains infinitely-many sheets, and is therefore noncompact.
    From above, we know that each component contains a sheet $W_{u}$ for some $u \in \Lambda$.
    Then we can choose infinitely-many distinct elements of $\Lambda_j$ (where $j \neq i$), represent them by loops in $M$ which are disjoint from the spheres $\Sigma_{i,g}$ for $g \in G$, and lift these loops to arcs based at the sheet $W_{u}$.
    Then the endpoint of each lift will be a distinct sheet in the same component of $\widetilde{M} \setminus \Sigma_{i,g_0,u_0}$.

    Finally, we can prove the ``moreover'' statement.
    Namely, assume $u_0,u \in \Lambda_j$ where $j \neq i$.
    Then we can represent $u_0^{-1}u$ by a loop in $M$ which is disjoint from the spheres $\Sigma_{i,g}$ for $g \in G$.
    This loop lifts to an arc from $W_{id}$ to $W_{u_0^{-1}u}$ which is disjoint from each sphere $\Sigma_{i,g,v}$ for $g \in G$ and $v \in \Lambda$.
    Applying the deck transformation $u_0$ to this arc, we get an arc from $W_{u_0}$ to $W_{u}$ which is disjoint from $\Sigma_{i,g_0,u_0}$.

\end{proof}

\begin{proof}[Proof of Proposition \ref{prop:infinite-lin-ind-subset}]
    Throughout this proof, we use homology with $\Z/\ell\Z$-coefficients.
    Lemma \ref{lem:infinite-subgroup} says that the group $\Lambda_2$ is infinite, so it's enough to show that $\Sigma_{1,g,u}$ and $\Sigma_{1,g',u'}$ represent distinct classes and that $\mathcal{S}$ is linearly independent.
    In particular, we must show that for any finite subset 
    \begin{equation*}
        \{(g_1,u_1), \ldots, (g_m, u_m)\} \subseteq G_{M_0} \times \Lambda_2,
    \end{equation*}
    the homology classes $[\Sigma_{1,g_1,u_1}], \ldots [\Sigma_{1,g_m,u_m}]$ in $H_{n-1}(\widetilde{M})$ are linearly independent.
    To do so, we will construct two subsets $X,Y \subseteq \widetilde{M}$ such that $X \cup Y = \widetilde{M}$ and such that each $\Sigma_{1,g_i,u_u}$ is a component of $X \cap Y$, and then analyze the resulting Mayer-Vietoris sequence. \\

    \textbf{Step 1: Constructing the subsets.}
    Let $V \subseteq M_0$ denote the complement of $\Int(D_{i,g}^+)$ for each $i \in \{1,2\}$ and $g \in G_{M_0}$, so in particular $V \subseteq M_0^*$.
    For $u \in \Lambda$, let $V_u = W_u \cap \pi^{-1}(V)$, so $V_u$ is a compact oriented $n$-manifold with
    \begin{equation*}
        \del V_u = \bigcup_{g \in G_{M_0}} \Sigma_{1,g,u} \cup \Sigma_{2,g,u}.
    \end{equation*}
    Now, let $X \subseteq \widetilde{M}$ be the union of the sets $V_{u_j}$ for $1 \leq j \leq m$.
    Then $X$ is a (possibly disconnected) compact oriented $n$-manifold with boundary
    \begin{equation*}
        \del X = \bigcup_{j=1}^m \bigcup_{g \in G_{M_0}} \Sigma_{1,g,u_j} \cup \Sigma_{2,g,u_j}.
    \end{equation*}
    Let $Y \subseteq \widetilde{M}$ be the closure of $\widetilde{M} \setminus X$, so $Y$ is a (possibly disconnected) non-compact oriented $n$-manifold with $\del Y = \del X = X \cap Y$.

    Note that $H_{n-1}(X \cap Y)$ is freely generated by the fundamental classes of its components.
    The components of $X \cap Y$ fall into two types: spheres of the form $\Sigma_{1,g,u_j}$ and spheres of the form $\Sigma_{2,g,u_j}$ for $1 \leq j \leq m$.
    This yields a decomposition
    \begin{equation*}
        H_{n-1}(X \cap Y) = L_1 \oplus L_2
    \end{equation*}
    where $L_i$ is freely generated by the components of the form $\Sigma_{i,g,u}$.

    Note also that we have a commutative diagram
    \begin{equation*}
        \begin{tikzcd}
            & {H_{n-1}(X)} & \\
            {H_{n-1}(X \cap Y)} && {H_{n-1}(\widetilde{M})} \\
            & {H_{n-1}(Y)}
            \arrow["{j_X}", from=1-2, to=2-3]
            \arrow["{i_X}", from=2-1, to=1-2]
            \arrow["{i_Y}"', from=2-1, to=3-2]
            \arrow["{j_Y}"', from=3-2, to=2-3]
        \end{tikzcd}
    \end{equation*}
    where the maps are induced by inclusion. \\

    \textbf{Step 2: Analyzing the components of $Y$.}
    Let $Y_1, \ldots, Y_r$ denote the components of $Y$.
    We claim that if a component $Y_s$ is compact, it must be the closure of a component of the preimage of $P_2^* \subseteq M$.
    Indeed, Lemma \ref{lem:separating-spheres} implies that any sphere of the form $\Sigma_{1,g,u_j} \subseteq X \cap Y$ must bound a non-compact component of $Y$, as each $u_j \in \Lambda_2$.
    Thus if $Y_s$ is compact, then each component of $\del Y_s$ is a sphere of the form $\Sigma_{2,g,u_j}$.
    Thus $Y_s$ must contain the closure of a component of the preimage of $P_2^*$.
    Then the closure of this component must be all of $Y_s$.
    Indeed, if it weren't, then $Y_s$ would contain a component of the preimage of $P_0^* \subseteq M$ (this follows since, by the same argument as in Lemma \ref{lem:infinite-subgroup}, any component of $P_2^*$ is either non-simply connected or intersects at least two components of $P_0^*$; thefore any component of the preimage of $P_2^*$ intersects at least two components of the preimage of $P_0^*$).
    Then we could find infinitely many components of the preimage of $P_0^*$ in $Y_s$ by lifting loops in $\Lambda_1$, contradicting that $Y_s$ is compact. \\

    \textbf{Step 3: Analyzing the kernel of $i_Y$.}
    Next, we claim that $\Ker(i_Y) \subseteq L_2$.
    To see this, we consider the exact sequence
    \begin{equation*}
        H_n(Y, \del Y) \xrightarrow{\delta_Y} H_{n-1}(\del Y) \xrightarrow{i_Y} H_{n-1}(Y).
    \end{equation*}
    Let $Y_1, \ldots, Y_r$ denote the components of $Y$.
    Then $H_n(Y,\del Y) \cong \oplus_{s=1}^r H_n(Y_s, \del Y_s)$.
    If $Y_s$ is compact, then the map $\delta_Y$ sends the fundamental class in $H_n(Y_s, \del Y_s)$ to the sum of the boundary components of $Y_s$.
    Otherwise, if $Y_s$ is non-compact, then the group $H_n(Y_s, \del Y_s)$ is trivial.
    As in Step 2, Lemma \ref{lem:separating-spheres} implies that each $\Sigma_{1,g_j,u_j}$ is the boundary of some noncompact component $Y$.
    Thus we can conclude that $\Ker(i_Y) = \Im(\delta_Y)$ intersects the subspace $L_1$ trivially. \\

    \textbf{Step 4: Proving linear independence.}

    Let $h = \sum_{j=1}^m c_j[\Sigma_{1,g_j,u_j}] \in H_{n-1}(X \cap Y)$ be any linear combination of the homology classes $[\Sigma_{1,g_1,u_1}], \ldots, [\Sigma_{1,g_m,u_m}] \in H_{n-1}(X \cap Y)$.
    Our goal is to show that if $j_X(i_X(h)) = 0$ (i.e.\ if this linear combination is trivial in $H_{n-1}(\widetilde{M})$), then $c_j = 0$ for each $1 \leq j \leq m$.

    If $j_X(i_X(h)) = 0$, then Mayer-Vietoris tells us that $(i_X(h), 0) \in H_{n-1}(X) \oplus H_{n-1}(Y)$ is equal to $(i_X(k), i_Y(k))$ for some $k \in H_{n-1}(X \cap Y)$.
    Since $i_Y(k) = 0$, it follows from Step 2 that $k \in L_2$.
    We also have that $h - k \in \Ker(i_X)$.

    Consider the exact sequence
    \begin{equation*}
        H_n(X, \del X) \xrightarrow{\delta_X} H_{n-1}(\del X) \xrightarrow{i_X} H_{n-1}(X).
    \end{equation*}
    Let $X_1, \ldots, X_t$ denote the components of $X$.
    The group $H_n(X, \del X)$ is freely generated by the fundamental classes $[X_s]$, and the map $\delta_X$ sends $[X_s]$ to $[\del X_s]$, i.e.\ to the sum of the fundamental classes of the boundary components of $X_s$.
    Thus since $h - k \in \Ker(i_X) = \Im(\delta_X)$, we have that 
    \begin{equation*}
        h - k = \sum_{s=1}^t d_s[\del X_s]
    \end{equation*}
    for some $d_s \in \Z/\ell\Z$.
    Since $h \in L_1$ and $k \in L_2$, we have that $h = 0$ if and only if each $d_s = 0$.
    So, we want to show that each $d_s = 0$.

    To see this, consider the following graph $\mathscr{H}$:
    \begin{itemize}
        \item The vertices of $\mathscr{H}$ correspond to the components $X_s$ with $d_s \neq 0$, and the compact components of $Y$ intersecting a component $X_s$ with $d_s \neq 0$.
        \item Two vertices are joined by an edge if the corresponding components intersect.
    \end{itemize}
    Thus each $d_s = 0$ if and only if $\mathscr{H}$ is empty.
    We can naturally view $\mathscr{H}$ (after subdividing the edges) as a subgraph of the graph $\widetilde{\mathscr{G}}$ introduced in Section \ref{subsec:infinite-cover}.
    Since $\widetilde{\mathscr{G}}$ is a tree by Lemma \ref{lem:cpt-graph-a-tree}, it follows that $\mathscr{H}$ is a finite tree.
    Therefore, to show that $\mathscr{H}$ is empty, it's enough to show that any vertex must have valence at least 2 (since any finite tree has a vertex of valence 0 or 1).

    First, we can show that each vertex corresponding to a component of $X$ has valence at least 2.
    Let $X_s$ be a component of $X$ with $d_s \neq 0$.
    Recall that $H_{n-1}(X \cap Y)$ has a natural basis (as a $\Z/\ell\Z$-module) given by the components of $X \cap Y$.
    Since $d_s \neq 0$, each component of $\del X_s$ has a nontrivial coefficient in the basis expression of $h - k$.
    Since $h \in L_1$ and $k \in L_2$, we in fact have that each component of $\del X_s$ of the form $\Sigma_{2,g,u_j}$ has a nontrivial coefficient in the basis expression of $k$.
    Since $k \in \Ker(i_Y)$, and (as observed in Step 2) $\Ker(i_Y)$ is spanned by linear combinations of boundaries of compact components of $Y$, we conclude that each boundary sphere of $X_s$ of the form $\Sigma_{2,g,u_j}$ intersects a compact component of $Y$, and therefore the vertex corresponding to $X_s$ has valence at least 2.

    Next, we can show that any vertex corresponding to a compact component $Y_s$ of $Y$ has valence at least 2.
    As observed in Step 2, $Y_s$ must be the closure of a component of the preimage of $P_2^* \subseteq M$, and $Y_s$ intersects multiple of components of $X$.
    By assumption, some boundary component of $Y_s$ has a nontrivial coefficient in the basis expression of $h - k$, and in particular has a nontrivial coefficient in $k$.
    But as observed in Step 2, $k$ is a linear combination of boundaries of compact components of $Y$.
    Thus each other boundary component of $Y_s$ also has a nontrivial coefficient in $k$, and hence in $h - k$.
    This means that for each component $X_{s'}$ intersecting $Y_s$, we have that $d_{s'} \neq 0$.
    Thus the vertex corresponding to $Y_s$ has valence at least 2.

\end{proof}

\subsection{Adaptations to the topological setting}
Lemma \ref{lem:local-model-of-action} carries over the topological setting under the assumption that the action of $G$ is locally linear.
In particular, the tangent space $T_xM$ and its action by $G$ can be replaced by $\R^n$ and the representation $\rho_x$ described in Section \ref{subsec:topological-setting}.
The remainder of this section carries over to the topological setting verbatim.

\section{The Obstruction Map}\label{sec:obstruction-map}
In this section, we will construct the obstruction map that we use to prove part (ii) of Theorem \ref{mainthm:ker-not-fg}.
Fix a manifold $M$ and a group $G$ satisfying the assumptions of Theorem \ref{mainthm:ker-not-fg}.
Then by assumption we have an oriented submanifold $B \subseteq M$ which is fixed by some element of $G$, and whose homology class in $H_{n-2}(M; \Z/\ell\Z)$ is trivial for some $\ell > 1$.
We fix a basepoint $b_0 \in B$ as in Lemma \ref{lem:local-model-of-action}.
Let $\pi:\widetilde{M} \rightarrow M$ be the infinite-sheeted cover described in Section \ref{subsec:infinite-cover}, and let $\widetilde{b}_0$ be a basepoint above $b_0$.
Recall that the cover $\widetilde{M}$ is trivial over the subset $M_0^* \subseteq M$, and so in particular it is trivial over $B$.
We let $\widetilde{B}$ be the component of $\pi^{-1}(B)$ containing the basepoint $\widetilde{b}_0$ (so $\pi$ restricts to a diffeomorphism $\widetilde{B} \cong B$).

Our obstruction map is not defined on $\Ker(\P_G)$, but a closely related group.
Let $\Diff(M,B)$ denote the group of orientation-preserving diffeomorphisms of $M$ that restrict to an orientation-preserving diffeomorphism of $B$,
and let $\Mod(M,B) = \pi_0(\Diff(M,B))$.
Let $\Mod_G(M,B) \leq \Mod(M,B)$ and $\Mod_G(M) \leq \Mod(M)$ denote the subgroups of mapping classes with an equivariant representative.
Then there is a forgetful map $\F:\Mod(M,B) \rightarrow \Mod(M)$ which restricts to a map
\begin{equation*}
    \F_G:\Mod_G(M,B) \rightarrow \Mod_G(M).
\end{equation*}
Our goal in this section is construct a map
\begin{equation*}
    \varphi:\Ker(\F_G) \rightarrow H_{n-1}(\widetilde{M},\widetilde{B}; \Z/\ell\Z).
\end{equation*}

We give a simple description of the map $\varphi$ in Section \ref{subsec:simple-description-of-obstruction-map}; the rest of this section is devoted to proving this map is well-defined.
In Section \ref{subsec:loops-of-submanifolds}, we describe how to obtain a homology class of $\widetilde{M}$ from an isotopy.
In Sections \ref{subsec:isotopies-equivariant-diffeos} and \ref{subsec:lifting-isotopies}, we describe how a subgroup of isotopies on $M$ lift to $\widetilde{M}$.
We formally construct $\varphi$ in Section \ref{subsec:defining-obstruction-map}.
Throughout this section, we implicitly use homology with $\Z/\ell\Z$-coefficients, but our arguments apply verbatim for other coefficients.

\subsection{A simple description of the map}\label{subsec:simple-description-of-obstruction-map}

The obstruction map $\varphi$ works as follows.
Let $\alpha \in \Ker(\F_G)$, and let $f$ be an equivariant diffeomorphism of $M$ that represents $\alpha$.
Then $f$ preserves $B$ setwise, and there is an isotopy $h:M \times [0,1] \rightarrow M$ with $h_0 = \id_M$ and $h_1 = f$.
This induces a map
\begin{align*}
    r:B \times [0,1] &\rightarrow M \\
    (b,t) &\mapsto h_t(b).
\end{align*}
In particular, $r\v_{B \times \{0\}} = \id_B$ and $r\v_{B \times \{1\}}$ is a diffeomorphism of $B$.
Now, lift $r$ to a map
\begin{equation*}
    \widetilde{r}:B \times [0,1] \rightarrow \widetilde{M}
\end{equation*}
which sends $B \times \{0,1\}$ to $\widetilde{B}$.
Then by pushing forward the fundamental class of $B \times [0,1]$, the map $\widetilde{r}$ yields an element $x \in H_{n-1}(\widetilde{M}, \widetilde{B})$.
We define $\varphi(\alpha) = x$.

It is not obvious that the map $\varphi$ is well-defined with this description.
One has to show that $\varphi$ is independent of the choice of representative $f$ and the choice of isotopy $h$.
One also has to show that the map $r$ lifts as desired; we can choose a lift $\widetilde{r}$ sending $B \times \{0\}$ to $\widetilde{B}$, but its not clear that this lift will map $B \times \{1\}$ to $\widetilde{B}$ (as opposed to some other component of $\pi^{-1}(B)$).

Thus, in this section we will construct $\varphi$ more formally so as to ensure that it is well-defined.
The above description will still be useful to compute $\varphi$ in practice (which we will do in Section \ref{sec:nontrivial-obstruction}).

\subsection{Homology classes from isotopies}\label{subsec:loops-of-submanifolds}

We begin by describing how an isotopy of $\widetilde{B}$ in $\widetilde{M}$ gives rise to a homology class.

There is a natural evaluation map
\begin{align*}
    \ev_{\widetilde{B}}:\Diff(\widetilde{M}) \times \widetilde{B} &\rightarrow \widetilde{M} \\
    (f,b) &\mapsto f(b).
\end{align*}
This map takes the subspace $\Diff(\widetilde{M},\widetilde{B}) \times \widetilde{B}$ to $\widetilde{B}$.
By composing the relative Eilenberg-Zilber map with the evaluation map, we get a map
\begin{equation*}
  H_1\left(\Diff(\widetilde{M}), \Diff(\widetilde{M},\widetilde{B})\right) \otimes H_{n-2}\left(\widetilde{B}\right) \rightarrow H_{n-1}\left(\widetilde{M}, \widetilde{B}\right).
\end{equation*}
By pairing an element of $H_1\left(\Diff(\widetilde{M}), \Diff(\widetilde{M},\widetilde{B})\right)$ with the fundamental class of $\widetilde{B}$, we get a map
\begin{equation*}
    \theta_{\widetilde{B}}:H_1\left(\Diff(\widetilde{M}), \Diff(\widetilde{M},\widetilde{B})\right) \rightarrow H_{n-1}\left(\widetilde{M}, \widetilde{B}\right).
\end{equation*}

Now, we define the group
\begin{equation*}
    \Iso(\widetilde{M},\widetilde{B}) \coloneqq \pi_1\left(\Diff(\widetilde{M}), \Diff(\widetilde{M},\widetilde{B}), \id_{\widetilde{M}}\right).
\end{equation*}
So, an element of $\Iso(\widetilde{M},\widetilde{B})$ is the homotopy class of an isotopy $h:\widetilde{M} \times [0,1] \rightarrow \widetilde{M}$ where $h_0 = \id_{\widetilde{M}}$ and $h_1$ restricts to a diffeomorphism of $\widetilde{B}$.
A priori, $\Iso(\widetilde{M},\widetilde{B})$ is just a (pointed) set, but it inherits a group structure from $\Diff(\widetilde{M})$ (cf.\ \cite{goldsmith}): given two elements $[f],[g] \in \Iso(\widetilde{M},\widetilde{B})$ represented by isotopies $f,g:\widetilde{M} \times [0,1] \rightarrow \widetilde{M}$ with $f_0 = g_0 = \id_{\widetilde{M}}$ and $f_1,g_1 \in \Diff(\widetilde{M},\widetilde{B})$, the product $[g] \cdot [f]$ is represented by the isotopy
\begin{equation*}
    (g \cdot f)_t \coloneqq
    \begin{cases}
        f_{2t} & 0 \leq t \leq 1/2, \\
        g_{2t-1} \circ f_1 & 1/2 \leq t \leq 1.
    \end{cases}
\end{equation*}

Observe that there is a natural function 
\begin{equation*}
    \widehat{\psi}:\Iso(\widetilde{M},\widetilde{B}) \rightarrow H_1\left(\Diff(\widetilde{M}), \Diff(\widetilde{M},\widetilde{B})\right).
\end{equation*}
With respect to the group structure on $\Iso(\widetilde{M},\widetilde{B})$ defined above, the function $\widehat{\psi}$ is \emph{not} a homomorphism.
Rather, it is a \emph{(right-)crossed homomorphism}: it satisfies the property that
\begin{equation*}
    \widehat{\psi}([g] \cdot [f])
    = \widehat{\psi}([g]) \cdot [f_1] + \widehat{\psi}([f]).
\end{equation*}
Here $\widehat{\psi}([g]) \cdot [f_1]$ refers to the right action of $\pi_0\left(\Diff(\widetilde{M},\widetilde{B})\right)$ on $H_1\left(\Diff(\widetilde{M}), \Diff(\widetilde{M}, \widetilde{B})\right)$ given by precomposition.

We define $\psi:\Iso(\widetilde{M},\widetilde{B}) \rightarrow H_{n-1}\left(\widetilde{M}, \widetilde{B}\right)$ to be the composite function
\begin{equation*}
  \Iso(\widetilde{M},\widetilde{B})
  \xrightarrow{\widehat{\psi}} H_1\left(\Diff(\widetilde{M}), \Diff(\widetilde{M},\widetilde{B})\right) \xrightarrow{\theta_{\widetilde{B}}} H_{n-1}\left(\widetilde{M}, \widetilde{B}\right).
\end{equation*}
We can interpret the function $\psi$ as follows.
Let $[h] \in \Iso(\widetilde{M},\widetilde{B})$, so $h:\widetilde{M} \times [0,1] \rightarrow \widetilde{M}$ is an isotopy with $h_0 = \id_{\widetilde{M}}$ and $h_1 \in \Diff(\widetilde{M},\widetilde{B})$.
We then get a map
\begin{align*}
    \widetilde{r}:\widetilde{B} \times [0,1] &\rightarrow \widetilde{M} \\
    (b,t) &\mapsto h_t(b).
\end{align*}
We have an induced map on relative homology
\begin{equation*}
    \widetilde{r}_*:H_{n-1}(\widetilde{B} \times [0,1], \widetilde{B} \times \{0,1\})
    \rightarrow H_{n-1}(\widetilde{M}, \widetilde{B}).
\end{equation*}
Then $\psi([h])$ is the image of the fundamental class under $\widetilde{r}_*$.

Since $\widehat{\psi}$ is a crossed homomorphism, the function $\psi$ is not obviously a homomorphism.
However, by appealing to the fact that elements of $\Diff(\widetilde{M},\widetilde{B})$ are required to preserve the orientation of $\widetilde{B}$, we can verify that it is indeed a homomorphism.

\begin{lemma}\label{lem:obstruction-map-is-hom}
    The function $\psi:\Iso(\widetilde{M},\widetilde{B}) \rightarrow H_{n-1}\left(\widetilde{M},\widetilde{B}\right)$ is a homomorphism.
\end{lemma}
\begin{proof}
    Let $[f],[g] \in \Iso(\widetilde{M},\widetilde{B})$ be two elements represented by isotopies $f,g:\widetilde{M} \times [0,1] \rightarrow \widetilde{M}$ such that $f_0 = g_0 = \id$ and $f_1,g_1 \in \Diff(\widetilde{M},\widetilde{B})$.
    Then
    \begin{equation*}
        \widehat{\psi}([g] \cdot [f]) = [g \circ (f_1 \times \id)] + [f] \in H_1\left(\Diff(\widetilde{M}),\Diff(\widetilde{M},\widetilde{B})\right).
    \end{equation*}
    Since $f_1$ preserves the orientation of $\widetilde{B}$, it follows that $f_1 \times \id$ acts trivially on the homology group $H_{n-1}(\widetilde{B} \times [0,1], \widetilde{B} \times \{0,1\})$, and so 
    \begin{equation*}
        \theta_{\widetilde{B}}([g \circ (f_1 \times \id)]) = \theta_{\widetilde{B}}([g]).
    \end{equation*}
    Thus we conclude that
    \begin{align*}
        \psi([g] \cdot [f])
        &= \theta_{\widetilde{B}}\left([g \circ (f_1 \times \id)] + [f]\right) \\
        &= \theta_{\widetilde{B}}([g \circ (f_1 \times \id)]) + \theta_{\widetilde{B}}([f]) \\
        &= \theta_{\widetilde{B}}([g]) + \theta_{\widetilde{B}}([f]) \\
        &= \psi([g]) + \psi([f]).
    \end{align*}
\end{proof}

\begin{remark}\label{rem:embedding-spaces}
    Let $\Emb(\widetilde{B},\widetilde{M})$ denote the space of embeddings $\widetilde{B} \hookrightarrow \widetilde{M}$, and let $\Sub(\widetilde{B},\widetilde{M})$ denote quotient of $\Emb(\widetilde{B},\widetilde{M})$ by $\Diff(\widetilde{B})$ (which acts by precomposition).
    Then a theorem of Cerf \cite{cerf} says that the map $\Diff(\widetilde{M}) \rightarrow \Emb(\widetilde{B},\widetilde{M})$ is a fibration, and a theorem of Binz-Fischer \cite{binz-fischer} says that the map $\Emb(\widetilde{B},\widetilde{M}) \rightarrow \Sub(\widetilde{B},\widetilde{M})$ is a fibration.
    Thus the map $\Diff(\widetilde{M}) \rightarrow \Sub(\widetilde{B},\widetilde{M})$ is a fibration with fiber $\Diff(\widetilde{M}, \widetilde{B}^\pm)$, where $\Diff(\widetilde{M},\widetilde{B}^\pm)$ is the group of diffeomorphisms preserving $B$ setwise (not necessarily preserving the orientation of $B$).
    It follows that $\pi_1(\Sub(\widetilde{B},\widetilde{M})) \cong \pi_1(\Diff(\widetilde{M}), \Diff(\widetilde{M},\widetilde{B}^\pm))$, so the function $\psi$ can be viewed as an invariant of loops in $\Sub(\widetilde{B},\widetilde{M})$.
    This was in fact our original motivation in defining the function $\psi$, but for the arguments in this paper it is generally simpler to work with the pair $(\Diff(\widetilde{M}),\Diff(\widetilde{M},\widetilde{B}))$ instead.
\end{remark}

\subsection{Isotopies and equivariant diffeomorphisms}\label{subsec:isotopies-equivariant-diffeos}

Next, we will identify the relationship between $\Ker(\F_G)$ and isotopies of $B$ in $M$.

Similar to above, we define the group
\begin{equation*}
    \Iso(M,B) \coloneqq \pi_1(\Diff(M),\Diff(M,B),\id_M),
\end{equation*}
where $\Diff(M,B)$ is the group of orientation-preserving diffeomorphisms of $M$ that restrict to an orientation-preserving diffeomorphism of $B$.
From the long exact sequence of homotopy groups for the pair $(\Diff(M), \Diff(M,B))$, we have an exact sequence
\begin{equation*}
  \pi_1\left(\Diff(M)\right) \rightarrow \Iso(M,B) \rightarrow \Mod(M, B) \xrightarrow{\F} \Mod(M).
\end{equation*}
Equivalently, we have an exact sequence
\begin{equation*}
  \pi_1\left(\Diff(M)\right) \rightarrow \Iso(M,B) \rightarrow \Ker(\F) \rightarrow 1.
\end{equation*}
The group $\Ker(\F_G)$ is a subgroup of $\Ker(\F)$; let $\Iso(M,B)_G \leq \Iso(M,B)$ denote the preimage of $\Ker(\F_G)$.
Thus we have an exact sequence
\begin{equation*}
  \pi_1\left(\Diff(M)\right) \rightarrow \Iso(M,B)_G \rightarrow \Ker(\F_G) \rightarrow 1.
\end{equation*}
The map $\Iso(M,B)_G \rightarrow \Ker(\F_G)$ is given by $[h] \mapsto [h_1]$, where $h:M \times [0,1] \rightarrow M$ is an isotopy with $h_0 = \id_M$ and $h_1 \in \Diff(M,B)$.
Any element of $\Iso(M,B)_G$ can be represented by an isotopy $h$ such that $h_1$ is equivariant.

\subsection{Lifting isotopies}\label{subsec:lifting-isotopies}

Our next goal is to show that every element of the group $\Iso(M,B)_G$ lifts to an element of $\Iso(\widetilde{M},\widetilde{B})$.

\begin{lemma}\label{lem:lifting-loops-of-submanifolds}
    There is a well-defined group homomorphism $\lambda_G:\Iso(M,B)_G \rightarrow \Iso(\widetilde{M},\widetilde{B})$ such that $\lambda_G([h])$ is represented by an isotopy $\widetilde{h}:\widetilde{M} \times [0,1] \rightarrow \widetilde{M}$ such that $\pi \circ \widetilde{h}_t = h_t$.
\end{lemma}


\begin{proof}
    Let $h:M \times [0,1] \rightarrow M$ be an isotopy with $h_0 = \id_M$ and $h_1 \in \Diff(M,B)$.
    Then $h$ lifts to an isotopy $\widetilde{h}:\widetilde{M} \times [0,1] \rightarrow \widetilde{M}$ with $\widetilde{h}_0 = \id_{\widetilde{M}}$.
    Since $\widetilde{h}$ is a lift of $h$, we know that $\widetilde{h}_1$ preserves the set $\pi^{-1}(B)$.
    Moreover, a homotopy of $h$ lifts to a homotopy of $\widetilde{h}$.
    We therefore have a map
    \begin{equation*}
        \Iso(M,B) \rightarrow \pi_1\left(\Diff(\widetilde{M}), \Diff(\widetilde{M},\pi^{-1}(B))\right).
    \end{equation*}
    Now, let $\Diff(\widetilde{M},\pi^{-1}(B), \widetilde{B})$ be the group of diffeomorphisms that restrict to an orientation-preserving diffeomorphism of $\pi^{-1}(B)$ and preserve the component $\widetilde{B}$ setwise.
    Then we have a natural inclusion
    \begin{equation*}
        \pi_1\left(\Diff(\widetilde{M}), \Diff(\widetilde{M},\pi^{-1}(B), \widetilde{B})\right) \hookrightarrow
        \pi_1\left(\Diff(\widetilde{M}), \Diff(\widetilde{M},\pi^{-1}(B))\right).
    \end{equation*}
    We claim that the image of $\Iso(M,B)_G$ is contained in the image of this inclusion.
    Assuming this, we have a map
    \begin{equation*}
        \Iso(M,B)_G \rightarrow \pi_1\left(\Diff(\widetilde{M}), \Diff(\widetilde{M},\pi^{-1}(B), \widetilde{B})\right),
    \end{equation*}
    and we can obtain the map $\lambda_G$ by post-composing with the natural map
    \begin{equation*}
        \pi_1\left(\Diff(\widetilde{M}), \Diff(\widetilde{M},\pi^{-1}(B), \widetilde{B})\right)
        \rightarrow \Iso(\widetilde{M},\widetilde{B}).
    \end{equation*}

    So, it remains to prove our claim.
    Let $h:M \times [0,1] \rightarrow M$ represent an element of $\Iso(M,B)_G$.
    This means that $h_1$ is an equivariant diffeomorphism $f:M \rightarrow M$ that retricts to an orientation-preserving diffeomorphism of $B$.
    Let $a:[0,1] \rightarrow M$ be the path $a(t) = h_t(b_0)$, i.e.\ $a$ is the path traveled by the basepoint $b_0$ under $h$.
    Then the lift $\widetilde{h}$ restricts to a lift $\widetilde{a}:[0,1] \rightarrow \widetilde{M}$ with $\widetilde{a}(0) = \widetilde{b}_0$.
    To check that $\widetilde{h}$ preserves the component $\widetilde{B}$, it is enough to check that $\widetilde{a}(1) \in \widetilde{B}$.
    We will first prove this in the case that $a$ is a loop, i.e.\ $f(b_0) = b_0$.
    Then we will address the general case.

    So, assume that $f(b_0) = b_0$.
    Our goal is to show that the loop $a$ lifts to a loop in $\widetilde{M}$.
    It's enough to show that $a$ represents an element of the subgroup $\pi_1(M_0,b_0) \leq \pi_1(M,b_0)$.
    Let $G_0 \leq G$ denote the stabilizer of $b_0$ (see Lemma \ref{lem:local-model-of-action} for a description of $G_0$).
    Then by Lemma \ref{lem:pi1-fixed-points}, it's enough to show that $G_0$ fixes the element $[a] \in \pi_1(M,b_0)$.
    Since $\pi_1(M,b_0)$ splits as a free product by Lemma \ref{lem:pi1-computation}, it has a trivial center, and therefore it's enough to show that $[a]^{-1}g_0([a])$ is central for any $g_0 \in G_0$.
  
    We can show that $[a]^{-1}g_0([a])$ is central as follows (we adapt an argument of Birman--Hilden \cite[Lem~1.3]{birman-hilden}).
    Let $c:S^1 \rightarrow M$ be any loop based at $b_0$.
    Applying the isotopy $h$ to $c$, we get that $c \simeq a \cdot f(c) \cdot \overline{a}$, where $\cdot$ denotes concatenation and $\overline{(-)}$ denotes the reverse path.
    Applying $g_0$, we conclude
    \begin{equation*}
        g_0(c) \simeq g_0(a) \cdot g_0(f(c)) \cdot \overline{g_0(a)}
        \simeq g_0(a) \cdot f(g_0(c)) \cdot \overline{g_0(a)}.
    \end{equation*}
    On the other hand, applying the isotopy $h$ directly to the element $g_0(c)$, we get 
    \begin{equation*}
        g_0(c) \simeq a \cdot f(g_0(c)) \cdot \overline{a}.
    \end{equation*}
    Thus, we conclude that 
    \begin{equation*}
        (\overline{a} \cdot g_0(a)) \cdot f(g_0(c)) \cdot \overline{(\overline{a} \cdot g_0(a))}
        \simeq f(g_0(c)).
    \end{equation*}
    Since any element of $\pi_1(M,b_0)$ is represented by $f(g_0(c))$ for some $c$ (as $f \circ g_0$ induces an automorphism of $\pi_1(M,b_0$)), we conclude that $[a]^{-1}g_0([a])$ is central as desired.
    
    Finally, we must address the case that $f(b_0) \neq b_0$.
    Note that in the argument above, we only used the fact that $f$ was $G_0$-equivariant.
    So, it's enough to find a $G_0$-equivariant ambient isotopy of $M$ taking $f(b_0)$ to $b_0$; we can then concatenate $h$ with this isotopy to reduce to the case $f(b_0) = b_0$.

    To build this isotopy, fix a $G_0$-invariant tubular neighborhood of $B$ (see \cite[Thm~VI.2.2]{bredon}).
    Fix an $(n-2)$-disk $D \subseteq B$ containing $b_0$ and $f(b_0)$, and choose an isotopy $k:D \times [0,1] \rightarrow D$ such that $k_0 = \id$, $k_1(f(b_0)) = b_0$, and $k_t(x) = x$ for all $x$ in a neighborhood of $\del D$ and $t \in [0,1]$.
    As in Lemma \ref{lem:local-model-of-action}, over $D$ the tubular neighborhood is diffeomorphic to $D \times D^2$, and $G_0$ acts by rotations on each $D^2$-fiber.
    Choose a smooth function $\mu:[0,1] \rightarrow [0,1]$ which equals $1$ in a neighborhood of $0$ and equals $0$ in a neighborhood of $1$.
    Then we define an isotopy $\widehat{k}:D \times D^2 \times [0,1] \rightarrow D \times D^2$ by $\widehat{k}_t(x, re^{i\theta}) = (k_{\mu(r)t}(x), re^{i\theta})$.
    Then $\widehat{k}_t$ is $G_0$-equivariant since it is independent of $\theta$.
    Moreover, since $\widehat{k}_t$ is the identity on a neighborhood of $\del(D \times D^2$), it extends to a $G_0$-equivariant ambient isotopy of $M$.
\end{proof}

\subsection{Defining the obstruction map}\label{subsec:defining-obstruction-map}

Finally, we can define our obstruction map by combining the maps $\lambda_G$ and $\psi$.
For emphasis, we make the $\Z/\ell\Z$-coefficients explicit moving forward.

We define the map $\widehat{\varphi}:\Iso(M,B)_G \rightarrow H_{n-1}(\widetilde{M},\widetilde{B};\Z/\ell\Z)$ to be the composition
\begin{equation*}
  \Iso(M,B)_G \xrightarrow{\lambda_G} \Iso\left(\widetilde{M}, \widetilde{B}\right)
  \xrightarrow{\psi} H_{n-1}(\widetilde{M},\widetilde{B};\Z/\ell\Z)
\end{equation*}
The main result of this section is then the following.

\begin{proposition}\label{prop:obstruction-map-descends}
    The map $\widehat{\varphi}:\Iso(M,B)_G \rightarrow H_{n-1}(\widetilde{M},\widetilde{B};\Z/\ell\Z)$ factors through a map
    \begin{equation*}
      \varphi:\Ker(\F_G) \rightarrow H_{n-1}(\widetilde{M},\widetilde{B};\Z/\ell\Z).
    \end{equation*}
\end{proposition}
\begin{proof}
    Recall that we have an exact sequence
    \begin{equation*}
      \pi_1\left(\Diff(M)\right) \rightarrow \Iso(M,B)_G \rightarrow \Ker(\F_G) \rightarrow 1.
    \end{equation*}
    Let $[h] \in \Iso(M,B)_G$ be an element in the image of $\pi_1\left(\Diff({M})\right)$.
    Our goal is to show that $\widehat{\varphi}([h])$ is trivial.
    Applying the map $\lambda_G$, we lift $[h]$ to an element $[\widetilde{h}] \in \Iso(\widetilde{M},\widetilde{B})$.
    Then $\widehat{\varphi}([h]) = \psi([\widetilde{h}])$, so it's enough to show that $\psi([\widetilde{h}]) = 0$.

    First, since $[h] \in \Iso(M,B)_G$ lies in the image of $\pi_1(\Diff(M))$, we may assume that $h_1 = \id_M$.
    Then, the lift $\widetilde{h}_1$ will be a deck transformation of $\widetilde{M}$.
    Since $\widetilde{h}_1$ preserves the component $\widetilde{B} \subseteq \pi^{-1}(B)$, and the deck group of $\widetilde{M}$ acts freely on the components of $\pi^{-1}(B)$ (as $B \subseteq M_0^*$ and the cover is trivial over $M_0^*$), it must be that $\widetilde{h}_1 = \id_{\widetilde{M}}$.
    Thus, $\widetilde{h}$ in fact defines an element of $H_1\left(\Diff(\widetilde{M});\Z/\ell\Z\right)$.

    Recall from Section \ref{subsec:loops-of-submanifolds} that we have an evaluation map
    \begin{equation*}
        \ev_{\widetilde{B}}:\Diff(\widetilde{M}) \times \widetilde{B} \rightarrow \widetilde{M},
    \end{equation*}
    and since the subspace $\Diff(\widetilde{M},\widetilde{B}) \times \widetilde{B}$ maps into $\widetilde{B}$, we get a pairing
    \begin{equation*}
        H_1\left(\Diff(\widetilde{M}), \Diff(\widetilde{M},\widetilde{B});\Z/\ell\Z\right) \otimes H_{n-2}(\widetilde{B};\Z/\ell\Z) \rightarrow
        H_{n-1}(\widetilde{M},\widetilde{B};\Z/\ell\Z).
    \end{equation*}
    On the other hand, we have have an evaluation map
    \begin{equation*}
        \ev_{\widetilde{M}}:\Diff(\widetilde{M}) \times \widetilde{M} \rightarrow \widetilde{M},
    \end{equation*}
    and hence a pairing
    \begin{equation*}
        H_1\left(\Diff(\widetilde{M});\Z/\ell\Z\right) \otimes H_{n-2}(\widetilde{M};\Z/\ell\Z) \rightarrow
        H_1(\widetilde{M};\Z/\ell\Z).
    \end{equation*}
    So, we have the following commutative diagram:
    \begin{equation*}
        \begin{tikzcd}
            {H_1(S^1; \Z/\ell\Z)} & {H_1\left(\Diff(\widetilde{M});\Z/\ell\Z\right)} & {H_{n-1}(\widetilde{M};\Z/\ell\Z)} \\
            {H_1([0,1],\{0,1\};\Z/\ell\Z)} & {H_1\left(\Diff(\widetilde{M}), \Diff(\widetilde{M},\widetilde{B});\Z/\ell\Z\right)} & {H_{n-1}(\widetilde{M},\widetilde{B};\Z/\ell\Z)}
            \arrow[from=1-1, to=1-2]
            \arrow["\cong"', from=2-1, to=1-1]
            \arrow[from=1-2, to=1-3]
            \arrow[from=1-3, to=2-3]
            \arrow[from=2-1, to=2-2]
            \arrow[from=2-2, to=2-3]
          \end{tikzcd}
    \end{equation*}
    The upper-left and lower-left horizontal maps are both induced by the isotopy $\widetilde{h}$.
    The upper right map comes from pairing a class with $[\widetilde{B}] \in H_{n-2}(\widetilde{M};\Z/\ell\Z)$ via the map $\ev_{\widetilde{M}}$, while the lower right map comes from pairing a class with the fundamental class $[\widetilde{B}] \in H_{n-2}(\widetilde{B};\Z/\ell\Z)$ via $\ev_{\widetilde{B}}$.

    The key observation is that the upper-right horizontal map is in fact the zero map.
    This is because the class $[B] \in H_{n-2}(M;\Z/\ell\Z)$ is trivial by assumption, which implies that the class $[\widetilde{B}] \in H_{n-2}(\widetilde{M};\Z/\ell\Z)$ is trivial (to see this, use the fact that $\pi:\widetilde{M} \rightarrow M$ restricts to a diffeomorphism $W_{\id} \cong M_0^*$, and apply Mayer-Vietoris to the subspaces $W_{\id}$ and $\widetilde{M} \setminus W_{\id}$).
    On the other hand, the composition of the two lower horizontal maps is $\psi([\widetilde{h}])$ by definition.
    Thus, we conclude that $\psi([\widetilde{h}])$ is zero as desired.
\end{proof}

\subsection{Adaptations to the topological setting}
The definition of the obstruction map carries over to the topological setting essentially verbatim, just by replacing the word ``diffeomorphism'' with ``homeomorphism'' wherever it appears.
The only place that does not immediately adapt is the end of the proof of Lemma \ref{lem:lifting-loops-of-submanifolds}. 
Here we make use of the assumption that the action of $G$ is locally linear. 
While we cannot necessarily choose a single disk containing $b_0$ and $f(b_0)$ on which $G$ has the desired local behavior, we can cover a path from $b_0$ to $f(b_0)$ with finitely many of the neighborhoods described in Section \ref{subsec:topological-setting} and apply the same argument finitely many times in order to equivariantly isotope $f$ to fix $b_0$.

\section{An Element with Nontrivial Obstruction}\label{sec:nontrivial-obstruction}

In Section \ref{sec:obstruction-map}, we defined an obstruction map
\begin{equation*}
    \varphi:\Ker(\F_G) \rightarrow H_{n-1}(\widetilde{M},\widetilde{B};\Z/\ell\Z).
\end{equation*}
Our next goal is to construct an element $\alpha_0 \in \Ker(\F_G)$ for which $\varphi(\alpha_0)$ is nontrivial.

First, we recall some notation from Section \ref{subsec:homology-of-cover}.
Let $G_{M_0} \leq G$ denote the stabilizer of $M_0$ (with respect to the action of $G$ on the components of $P_0$).
Then the gluing disks $D_{i,g}^+$ for $g \in G_{M_0}$ are precisely the gluing disks that lie on $M_0$.
We let $\Sigma_{i,g} \subseteq M$ denote the boundary of the gluing disk $D_{i,g}^+$, and for $g \in G_{M_0}$ and $u \in \Lambda$, we let $\Sigma_{i,g,u}$ denote the copy of $\Sigma_{i,g}$ on the sheet $W_u$ above $M_0$.

Let $G_B \leq G$ denote the setwise stabilizer of $B$, so $G_B \leq G_{M_0}$.
For $(g,u) \in G_{M_0} \times \Lambda$, we let $[\Sigma_{1,g,u}]$ denote the homology class of $\Sigma_{1,g,u}$ in $H_{n-1}(\widetilde{M};\Z/\ell\Z)$.
Since the map $H_{n-1}(\widetilde{M};\Z/\ell\Z) \rightarrow H_{n-1}(\widetilde{M},\widetilde{B};\Z/\ell\Z)$ is injective (as $\widetilde{B}$ has dimension $n-2$), we can equivalently view $[\Sigma_{1,g,u}]$ as an element of $H_{n-1}(\widetilde{M},\widetilde{B};\Z/\ell\Z)$.
Then we have the following.

\begin{proposition}\label{prop:nontrivial-obstruction}
    There exists an element $\alpha_0 \in \Ker(\F_G)$ such that
    \begin{equation*}
      \varphi(\alpha_0) = \sum_{g \in G_B} [\Sigma_{1,g,\id}].
    \end{equation*}
\end{proposition}

We will compute $\varphi(\alpha_0)$ using the description of $\varphi$ in Section \ref{subsec:simple-description-of-obstruction-map}.

\subsection{An informal description of the construction}\label{subsec:informal-description}
In order to motivate the construction of this section, we start with an informal picture of the mapping class $\alpha_0$.
Assume for simplicity that $G_B = G$.
Assume also that the gluing disks on $M_0$ are positioned in a tubular neighborhood of $B$.
We will construct an equivariant diffeomorphism $f$ as a summand slide which rotates all the gluing disks $D_{1,g}^+$ around $B$, and which fixes the remaining gluing disks.
The path traveled by the center of each $D_{1,g}^+$ is a meridian of $B$.

\begin{figure}[h]
    \centering
    \includegraphics[scale=.55]{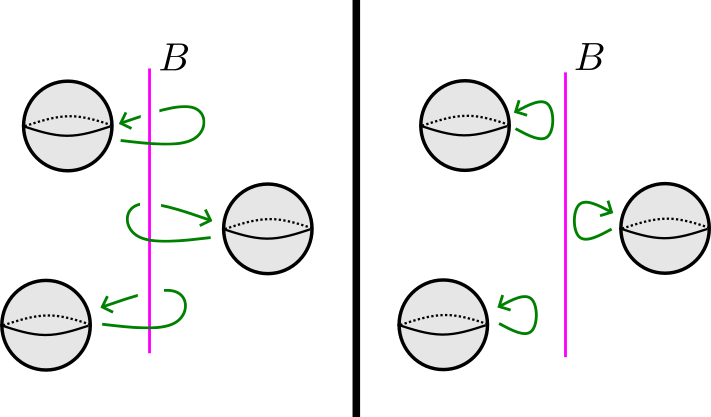}
    \caption{We isotope $f$ to $\id_M$ by sliding the disks $D_{1,g}^+$ along shorter and shorter loops.}
    \label{fig:informal-isotopy}
\end{figure}

We can then (non-equivariantly) isotope to $f$ to $\id_M$ by modifying the path traveled by each $D_{1,g}^+$.
In particular, we pull the path across $B$ and then shrink it down to the trivial loop at the center of each $D_{1,g}^+$, as in Figure \ref{fig:informal-isotopy} (in our actual construction, this will isotope $f$ to a product of sphere twists, so we actually define $\alpha_0$ by multiplying $f$ with a product of sphere twists).

\begin{figure}[h]
    \centering
    \includegraphics[scale=.75]{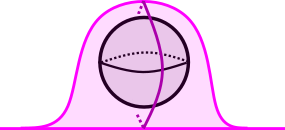}
    \caption{During the isotopy $f_t$, the sliding path of the disks $D_{1,g}^+$ runs into $B$, creating a ``bubble'' around each disk (cf.\ Figures \ref{fig:disk-scc} and \ref{fig:disk-scc-isotopy}).}
    \label{fig:bubble}
\end{figure}

This gives us an isotopy $f_t$ with $f_0 = f$ and $f_1 = \id_M$, where each $f_t$ is a summand slide.
For some range of times $t_1 \leq t \leq t_2$, the summand slide $f_t$ will push the gluing disk through $B$ (this occurs as we pull the sliding path across $B$).
The end result is that in the image of the map $B \times [0,1] \rightarrow M$ given by $(b,t) \mapsto f_t(b)$, we get a ``bubble'' around each gluing disk $D_{1,g}^+$ (see Figure \ref{fig:bubble}).
Thus in homology, the map $B \times [0,1] \rightarrow M$ yields the sum of all the classes $[\Sigma_{1,g}]$ for $g \in G$.
Since the cover $\pi:\widetilde{M} \rightarrow M$ restricts to a diffeomorphism $W_{\id} \rightarrow M_0$, we see the same phenomenon when we lift this isotopy to $\widetilde{M}$.

\subsection{A family of diffeomorphisms of the 2-disk}\label{subsec:diffeos-of-disk}

Let $D^2$ denote the closed unit $2$-disk in $\R^2$, and let $C$ denote the group generated by a $2\pi/m$-rotation.
We begin by defining a family of $C$-equivariant diffeomorphisms $F_s:D^2 \rightarrow D^2$ for $s \in [0,1]$.
In our construction of $\alpha_0$, we will apply these diffeomorphisms to the fibers of a tubular neighborhood of $B$.

\begin{figure}[h]
    \centering
    \includegraphics[scale=.5]{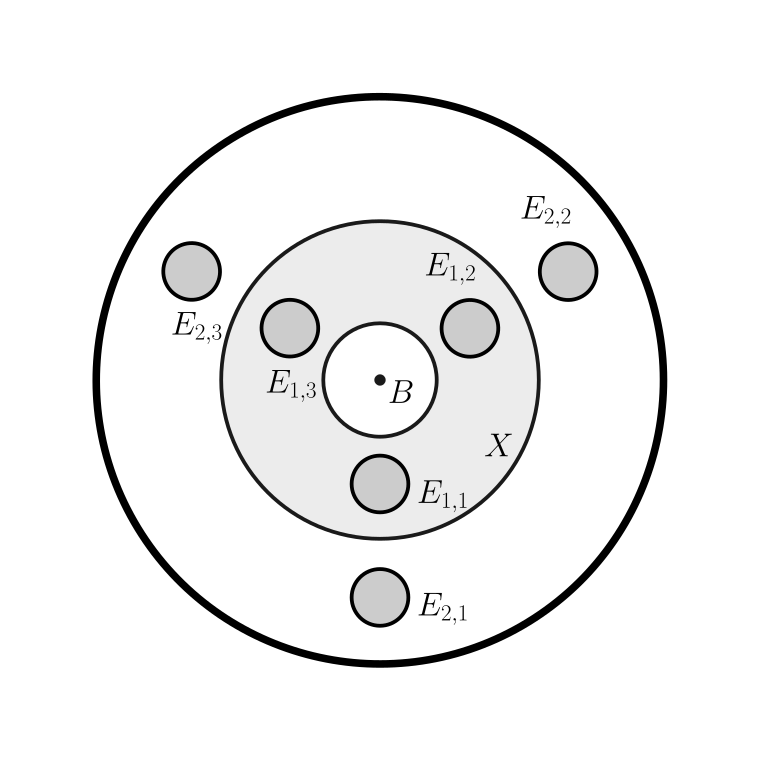}
    \caption{The disks $E_{i,j}$ and the annulus $X$.  The submanifold $B$ intersects $D^2$ at the origin.}
    \label{fig:labeled-disk}
\end{figure}

First, for $i \in \{1,2\}$ and $j \in \{1, \ldots, m\}$, let $E_{i,j}$ be a smaller $2$-disk in $D^2$ such that for each $i$, the disks $E_{i,1}, \ldots, E_{i,m}$ are permuted by $C$.
We choose the disks small enough so that $E_{i,j}$ is disjoint from $E_{i',j'}$ for all distinct pairs $(i,j)$ and $(i',j')$, and so that no $E_{i,j}$ contains the origin.
We arrange the disks so that for each $j$, the centers of the disks $E_{1,j}$ and $E_{2,j}$ lie on the same radial line in $D^2$, with $E_{1,j}$ closer to the origin and $E_{2,j}$ farther.
See Figure \ref{fig:labeled-disk} for an illustration.

Let $(r,\theta)$ denote polar coordinates on $D^2 - \{0\}$.
Choose $r_1 < r_2$ such that the annulus $X \coloneqq [r_1, r_2] \times [0,2\pi)$ contains each disk $E_{1,j}$ for $1 \leq j \leq m$ within its interior, and such that the disks $E_{2,j}$ for $1 \leq j \leq m$ are disjoint from $X$ (see Figure \ref{fig:labeled-disk}).
Choose $\ep > 0$ such that the regions $[r_1, r_1+\ep] \times [0,2\pi)$ and $[r_2-\ep, r_2] \times [0,2\pi)$ are disjoint from all the disks $E_{i,j}$.
Let $\xi:[0,\ep] \rightarrow [0,1]$ be a smooth function which equals $0$ in a neighborhood of $0$ and equals $1$ in a neighborhood of $\ep$.
Finally, let $S_r \subseteq D^2 - \{0\}$ be the circle of radius $r \in (0, 1]$.
Then for $s \in [0,1]$, we define a $C$-equivariant diffeomorphism $F_s:D^2 \rightarrow D^2$ as follows:
\begin{itemize}
  \item For $r \in [r_1, r_1+\ep]$, the map $F_s$ acts by a $\xi(r - r_1)2\pi s$-rotation on $S_r$.
  \item For $r \in [r_1 + \ep, r_2 - \ep]$, the map $F_s$ acts by a $2\pi s$ rotation on $S_r$.
  \item For $r \in [r_2 - \ep, r_2]$, the map $F_s$ acts by a $\xi(r_2 - r)2\pi s$-rotation.
  \item Outside of the annulus $X$, the map $F_s$ acts by $\id_{D^2}$.
\end{itemize}
In particular, $F_0 = \id_{D^2}$ and $F_1$ is a left Dehn twist about the inner boundary component of $X$ composed with a right Dehn twist about the outer boundary component of $X$.


\subsection{An equivariant diffeomorphism: Step one}\label{subsec:equiv-diffeo-1}

Our next goal is to use the family of diffeomorphisms $F_s$ to construct an equivariant diffeomorphism $f:M \rightarrow M$.
The diffeomorphism $f$ will be a summand slide supported on $P_0$, but we will need to carefully describe it in coordinates to compute $\varphi$.
As a first step, we will construct a disk slide $f':M_0 \rightarrow M_0$.

First, apply Lemma \ref{lem:local-model-of-action} to the action of $G_{M_0}$ on $M_0$ to obtain a point $b_0 \in B$ and associated closed neighborhood $Z_{b_0} \subseteq M_0$.
We let $G_0$ denote the $G_{M_0}$-stabilizer of $b_0$, so $G_0$ is cyclic.
In the following discussion, we will implicitly identify $Z_{b_0}$ with $D^{n-2} \times D^2$.

Without loss of generality, we assume that the gluing disks $D_{i,\id}^+$ for $i \in \{1,2\}$ lie in the interior of $Z_{b_0}$, and thus so do the gluing disks $D_{i,g}^+$ for $g \in G_0$.
Without loss of generality, we may also assume that for all $g \in G_0$ and $i \in \{1,2\}$, there is a unique $j \in \{1, \ldots, m\}$ such that the gluing disk $D_{i,g}^+$ is contained in the interior of $D^{n-2} \times E_{i,j}$.

Fix a closed $(n-2)$-disk $K \subseteq \Int(D^{n-2})$ such that under the projection map $Z_{b_0} \rightarrow D^{n-2}$, the image of each gluing disk $D_{i,g}^+$ lies in $\Int(K)$.
Choose a larger closed disk $K' \subseteq \Int(D^{n-2})$ such that $K \subseteq \Int(K')$.
Let $\nu:D^{n-2} \rightarrow [0,1]$ be a smooth function which is identically $1$ on $K$ and which is identically $0$ outside of $K'$.
Then, we define $f':M_0 \rightarrow M_0$ to be the following $G_0$-equivariant diffeomorphism:
\begin{itemize}
  \item on $Z_{b_0}$, $f'$ acts on the fiber $\{x\} \times D^2$ by the diffeomorphism $F_{\nu(x)}$,
  \item outside of $Z_{b_0}$, $f'$ acts by $\id_{M_0}$.
\end{itemize}

Note that since $F_1$ fixes the disks $E_{i,j} \subseteq D^2$ pointwise, it follows that $f'$ fixes all of the gluing disks in $M_0$.
In fact, $f'$ represents a disk slide in $\Mod(M_0, \Delta^0)$, where $\Delta^0$ the union of all the gluing disks on $M_0$.
Indeed, we can build an isotopy from $f'$ to $\id_{M_0}$ by smoothly homotoping the function $\nu$ to the constant $0$ map.
Under this isotopy, the $2$-disks $E_{1,1}, \ldots, E_{1,m}$ move a full rotation about the origin in each $D^2$-fiber over $K$.

\subsection{An equivariant diffeomorphism: Step two}\label{subsec:equiv-diffeo-2}

Next, we will use $f'$ to construct the desired equivariant diffeomorphism $f:M \rightarrow M$.
To do so, we will first define a disk slide $f''$ on $P_0$, which will then induce a summand slide on $M$.

Let $\{b_0, \ldots, b_p\} \subseteq P_0$ denote the $G$-orbit of $b_0$, and choose $g_j \in G$ such that $g_jb_0 = b_j$.
Let $Z_{b_j} = g_jZ_{b_0}$ (this set is independent of the choice of $g_j$).
Then we define $f'':P_0 \rightarrow P_0$ to be the following diffeomorphism:
\begin{itemize}
  \item on $Z_{b_j}$, $f''$ acts by $g_j \circ f' \circ g_j^{-1}$,
  \item outside of $\bigcup_{j=0}^p Z_{b_j}$, $f''$ acts by $\id_{P_0}$.
\end{itemize}
This definition is independent of the choice of $g_j$.
Moreover, one can directly check that $f''$ is $G$-equivariant. 

The diffeomorphism $f''$ represents a disk slide in $\Mod(P_0, \Delta_1^+ \cup \Delta_2^+)$ (recall that $\Delta_i^+$ is comprised of the disks $D_{i,g}^+$ for $g \in G$).
Indeed, the isotopy $f' \simeq \id_{M_0}$ induces an isotopy $f'' \simeq \id_{P_0}$.

Since $f''$ fixes each gluing disk pointwise, it extends a diffeomorphism $f:M \rightarrow M$ which acts by the identity outside of $P_0^*$.
By definition, the diffeomorphism $f$ is a summand slide.
Moreover, since $f''$ is equivariant, it follows that $f$ is equivariant.
Note that the isotopy $f'' \simeq \id_{P_0}$ does not fix the gluing disks at all times, and therefore does \emph{not} induce an isotopy from $f$ to $\id_M$.

\subsection{Identifying a product of sphere twists}\label{subsec:identifying-sphere-twists}

As discussed above, the diffeomorphism $f$ is a summand slide, and is induced by a disk slide $f''$ on $P_0$.
Let $d_{i,g}$ denote the center of the gluing disk $D_{i,g}^+$.
We can see from our explicit description of the map $f'$ that the path traveled by $d_{i,g}$ under the isotopy $f'' \simeq \id_{P_0}$ is null-homotopic in $P_0$.
Therefore, by Proposition \ref{prop:disk-slides-twists-and-inner-auts}, the mapping class $[f] \in \Mod(M)$ lies in the subgroup generated by sphere twists about the spheres $\Sigma_{i,g} = \del D_{i,g}^+$ for $(i,g) \in \{1,2\} \times G$.

In fact, we can explicitly determine the class $[f]$ as a product of sphere twists.
For each point $b_j$, we have a diffeomorphism $Z_{b_j} \cong D^{n-2} \times D^2$.
Under this identification, each tangent space in $Z_{b_j}$ is naturally isomorphic to $\R^{n-2} \times \R^2$.
Fix a frame at $d_{i,g}$ corresponding to the standard basis of $\R^{n-2} \times \R^2$.
For $i = 2$, the frame is fixed throughout the isotopy $f'' \simeq \id_{P_0}$.
For $i = 1$, the basis vectors in $\R^{n-2}$ remain fixed, while the basis vectors in $\R^2$ trace out a $2\pi$-rotation.
Thus, we conclude that in $\Mod(M)$, we have an equality
\begin{equation*}
    [f] = \prod_{g \in G} \tau_{1,g},
\end{equation*}
where $\tau_{1,g} \in \Mod(M)$ is the isotopy class of a sphere twist about $\Sigma_{1,g} = \del D_{1,g}^+$.




\subsection{An explicit isotopy}\label{subsec:an-explicit-isotopy}

Next, we construct an isotopy $f_t:M \rightarrow M$ for $t \in [0,1]$ with $f_0 = f$.

\begin{figure}[h]
    \centering
    \includegraphics[scale=.4]{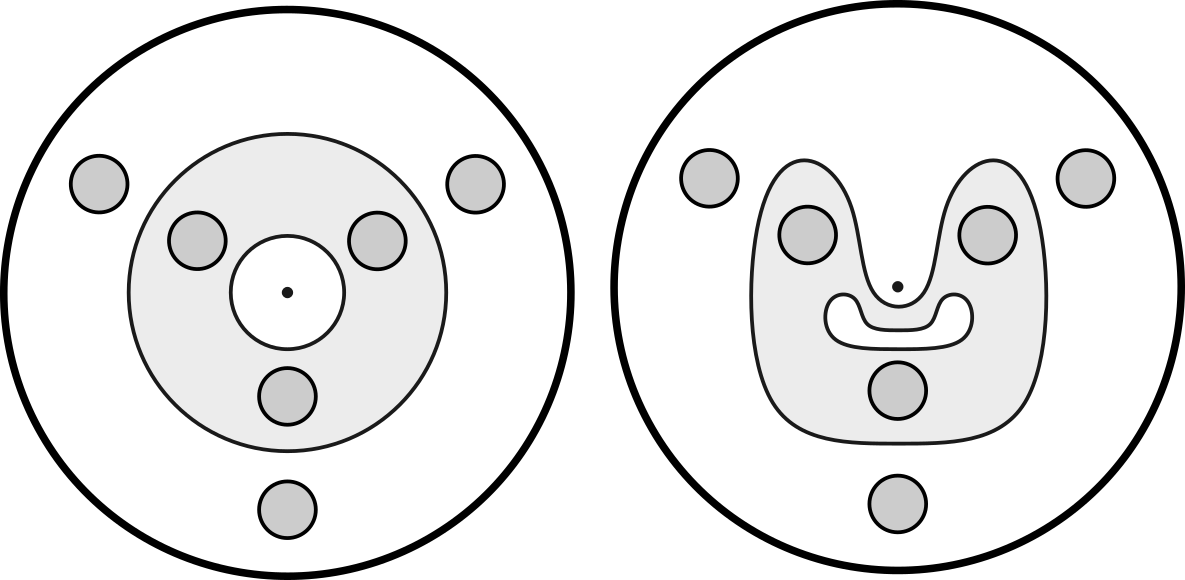}
    \caption{The annulus $X$ at time $t=0$ and $t=1$.}
    \label{fig:disk-isotopy}
\end{figure}

First, we construct an isotopy of the diffeomorphisms $F_s:D^2 \rightarrow D^2$.
Choose an ambient isotopy $k:D^2 \times [0,1] \rightarrow D^2$ with $k_0 = \id_{D^2}$ which fixes the disks $E_{i,j}$ pointwise at all times and which pushes the annulus $X$ across the origin, so that both boundary components of $k_1(X)$ have winding number $0$ about the origin (see Figure \ref{fig:disk-isotopy}).
Then, for $t \in [0,1]$, define $F_{s,t} \coloneqq k_t \circ F_s \circ k_t^{-1}$.
Thus $F_{s,0} = F_s$ for all $s$, $F_{0,t}$ is the identity for all $t$, and $F_{1,t}$ is a product of opposite Dehn twists about the two boundary components of $k_t(X)$.

Next, we can build the isotopy of $f$.
For $t \in [0,1]$, let $f_t':M_0 \rightarrow M_0$ be the diffeomorphism defined the same as $f'$, except on each fiber $\{x\} \times D^2$ of $Z_{b_0}$, $f_t'$ acts by the diffeomorphism $F_{\nu(x),t}$ instead of $F_{\nu(x)}$.
Then we build a disk slide $f_t''$ on $P_0$ exact as how we built $f''$ from $f'$, and $f_t''$ induces a summand slide $f_t:M \rightarrow M$.

Note that for each fixed $t$ we have an isotopy $f_t' \simeq \id_{M_0}$ constructed the same way as the isotopy $f' \simeq \id_{M_0}$, and this induces an isotopy $f_t'' \simeq \id_{P_0}$.
The key observation is that the path traveled by the center $d_{i,g}$ of $D_{i,g}^+$ under the isotopy $f_1'' \simeq \id_{P_0}$ is null-homotopic in $P_0 \backslash B$ (note that when we remove $B$ from $P_0$, this removes the origin from each $D^2$-fiber of each $Z_b$).
Thus in the mapping class group $\Mod(M \setminus B)$, we have the equality
\begin{equation}\label{eqn:sphere-twist-equality}
    [f_1] = \prod_{g \in G} \tau_{1,g}.
\end{equation}

\subsection{Constructing the element}\label{subsec:constructing-the-element}

Now, we will construct the desired element $\alpha_0 \in \Ker(\F_G)$. 
In fact, we will construct an isotopy $h:[0,1] \rightarrow \Diff(M)$ such that $h_0 = \id_M$ and $h_1$ is equivariant and preserves $B$ setwise; we can then let $\alpha_0 = [h_1]$.
The isotopy $f_{1-t}$ almost works for our purposes, but the problem is that the diffeomorphism $f = f_0$ is isotopic to a product of sphere twists, and thus may not be isotopic to $\id_M$.

Let $T_{i,g}:M \rightarrow M$ be a sphere twist about $\del D_{i,g}^-$ supported in a collar neighborhood $\del D_{i,g}^- \times [0,1] \subseteq P_i \backslash \Int(D_{i,g}^-)$.
Since $\del D_{i,g}^+$ and $\del D_{i,g}^-$ are isotopic spheres in $M \setminus B$, we know that $T_{i,g}$ represents the mapping class $\tau_{i,g}$ (in both $\Mod(M)$ and $\Mod(M \setminus B)$).
Define $h':[0,1] \rightarrow \Diff(M)$ to be the isotopy
\begin{equation*}
    h_t' \coloneqq f_{1-t} \circ \prod_{g \in G} T_{1,g}.
\end{equation*}
Recall from Section \ref{subsec:disk-slides-inner-auts} that the mapping classes $\tau_{i,g}$ have order 2, and since the representatives $T_{i,g}$ have disjoint support, these mapping classes also commute.
Thus from the equality (\ref{eqn:sphere-twist-equality}) in $\Mod(M \setminus B)$, it follows that there is an isotopy $\id_M \simeq h_0'$ that fixes $B$ at all times.
Thus, we construct $h$ by concatenating the isotopy $\id_M \simeq h_0'$ with the isotopy $h'$.

\subsection{Analyzing the image}

Let $r:B \times [0,1] \rightarrow M$ be the map $r(b,t) = h_t(b)$.
Our next goal is to describe the image of $r$.
Recall that $b_0, \ldots, b_p$ are the points in the $G$-orbit of $b_0$.
Assume without loss of generality that for some $k \leq p$, the points $b_0, \ldots, b_k$ lie on $B$, and the points $b_{k+1}, \ldots, b_p$ do not lie on $B$.
Let $Z_{b_j}^* = Z_{b_j} \cap P_0^* \subseteq M$.
Note that $Z_{b_j}^*$ is diffeomorphic to $D^{n-2} \times D^2$ with a finite number of points removed (namely, the center of each gluing disk in $Z_{b_j}$).

\begin{figure}[h]
    \centering
    \includegraphics[scale=.4]{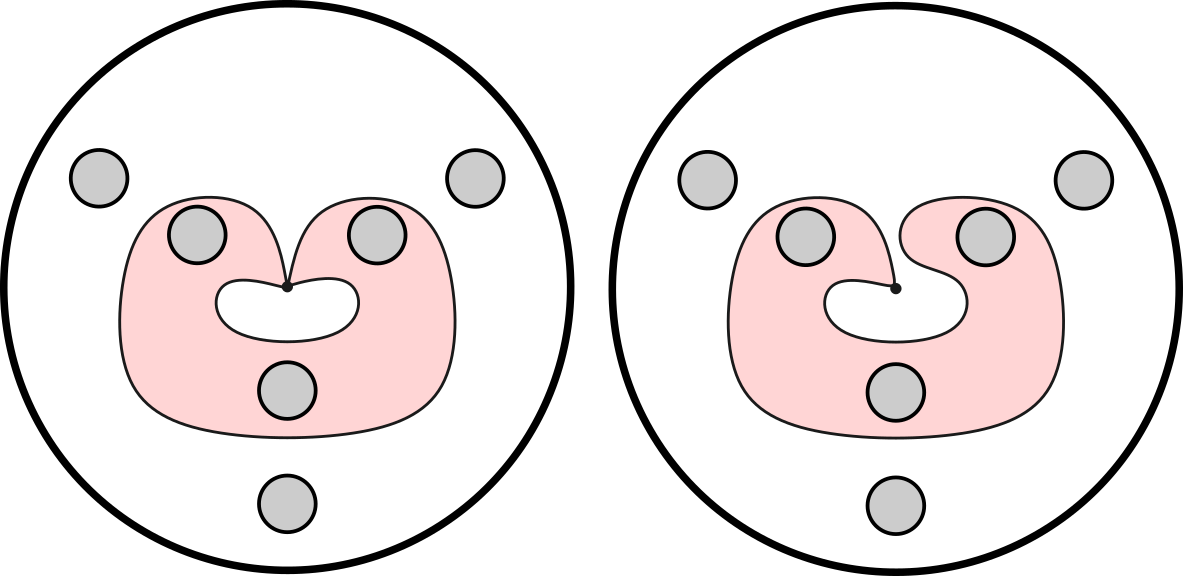}
    \caption{On the left is the loop $r:\{b\} \times [0,1] \rightarrow M$ for $b \in K$.  On the right is a slight perturbation of this loop for clarity.}
    \label{fig:disk-scc}
\end{figure}

\begin{figure}[h]
    \centering
    \includegraphics[scale=.4]{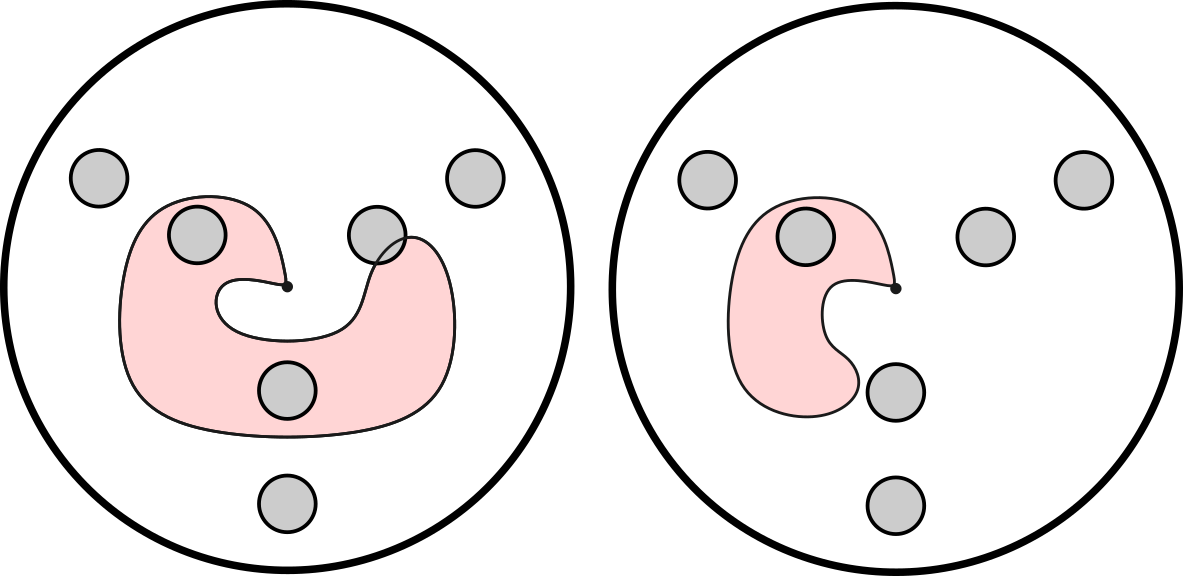}
    \caption{The loop $r:\{b\} \times [0,1] \rightarrow M$ for $b \in K' \setminus K$.  On the left is the case that $b$ is closer to $\del K$, on the right is the case that $b$ is closer to $\del K'$.
    For some values of $b \in K' \setminus K$, the loop will intersect the disks $E_{1,j}$.}
    \label{fig:disk-scc-isotopy}
\end{figure}

Now, let $b \in B$; then we can describe the loop $r:\{b\} \times [0,1] \rightarrow M$ as follows.
Suppose first that $b$ lies in $Z_{b_0}^*$.
Let $K,K' \subseteq D^{n-2}$ be the subsets described in the construction of $f'$; we will identify $K$ and $K'$ as subsets of $B$.
If $b$ lies in the $(n-2)$-disk $K$, then following the description of the isotopy $F_{\nu(b),t} = F_{1,t}$, the restriction $r:\{b\} \times [0,1] \rightarrow M$ is a closed curve in the fiber $(\{b\} \times D^2) \cap P_0^* \subseteq Z_{b_0}^*$. 
One side of this loop contains the disks $E_{1,j}$ for $j \in \{1,\ldots,m\}$, and the other side contains the disks $E_{2,j}$ for $j \in \{1,\ldots, m\}$ (see Figure \ref{fig:disk-scc}).
As we vary from $b$ from $\del K$ towards $\del K'$, this loop shrinks down towards the origin, and for $b$ outside of $K'$, the loop $r:\{b\} \times [0,1] \rightarrow M$ is the constant loop at $b$ (see Figure \ref{fig:disk-scc-isotopy}).

Next, observe that if $b$ lies in $Z_{b_j}^*$ for $1 \leq j \leq k$, then by applying the the element $g_j$ to $Z_{b_0}$, we see the same picture of the loop $r:\{b\} \times [0,1] \rightarrow M$ as in the case that $b \in Z_{b_0}^*$.
Finally, if $b$ lies outside of the sets $Z_{b_j}$ for each $0 \leq j \leq k$, then $r(b,t) = b$ for all $t \in [0,1]$.

\subsection{Computing the obstruction}

Finally, we can compute $\varphi(\alpha_0)$ following the description of $\varphi$ given in Section \ref{subsec:simple-description-of-obstruction-map}.
Let $h:[0,1] \rightarrow \Diff(M)$ be the isotopy constructed above, so $h_0 = \id_M$ and $h_1$ represents $\alpha_0$.
Let $r:B \times [0,1] \rightarrow M$ be the map $r(b,t) = h_t(b)$ as above.
Note that the image of $r$ is contained in $M_0^* \subseteq M$.

Throughout the remainder of this section, we use homology with $\Z/\ell\Z$-coefficients.
To prove Proposition \ref{prop:nontrivial-obstruction}, it's enough to show that the homology class $r_*([B \times [0,1]]) \in H_{n-1}(M_0^*, B)$ is equal to the sum of the classes $[\Sigma_{1,g}]$ for $g \in G_B$.
Indeed, the image of the lift $\widetilde{r}:B \times [0,1] \rightarrow \widetilde{M}$ is contained in the sheet $W_{\id}$.
The diffeomorphism $W_{\id} \cong M_0^*$ induces an isomorphism
\begin{equation*}
    H_{n-1}(M_0^*, B) \cong H_{n-1}(W_{\id}, \widetilde{B})
\end{equation*}
which takes $[\Sigma_{1,g}]$ to $[\Sigma_{1,g,\id}]$.
Then by Mayer–Vietoris, the natural map
\begin{equation*}
    H_{n-1}(W_{\id}, \widetilde{B}) \rightarrow H_{n-1}(\widetilde{M}, \widetilde{B})
\end{equation*}
is an embedding (cf.\ the proof of Proposition \ref{prop:infinite-lin-ind-subset}).

Moreover, since $h_1$ fixes $B$ pointwise, the map $r$ descends to a map $\overline{r}:B \times S^1 \rightarrow M_0^*$.
We have the following commutative diagram:
\begin{equation*}
    \begin{tikzcd}
        {H_{n-1}(B \times [0,1],B\times\{0,1\})} & {H_{n-1}(B \times S^1, B \times \{*\})} & {H_{n-1}(B \times S^1)} \\
        {H_{n-1}(M_0^*,B)} && {H_{n-1}(M_0^*)}
        \arrow["\cong", from=1-1, to=1-2]
        \arrow["{r_*}"', from=1-1, to=2-1]
        \arrow["\cong"', from=1-3, to=1-2]
        \arrow["{\overline{r}_*}", from=1-3, to=2-3]
        \arrow[hook', from=2-3, to=2-1]
    \end{tikzcd}
\end{equation*}
Therefore, it's enough to show that the map $\overline{r}_*$ takes the fundamental class of $B \times S^1$ to the sum of the classes $[\Sigma_{1,g}] \in H_{n-1}(M_0^*)$ for $g \in G_B$.

For $j \in \{0, \ldots, k\}$, let $K_j$ denote the copy of the $(n-2)$-disk $K$ in $Z_{b_j}$, so each $K_j$ is an $(n-2)$-disk in $B$ (recall that $b_0, \ldots, b_k$ are the points in the $G$-orbit of $b_0$ that lie on $B$).
Let $B^\circ$ be the closure of $B \setminus \bigcup_{j=0}^k K_j$.
For each $b \in \del K_j$, we saw above that the loop $r:\{b\} \times [0,1] \rightarrow \{b\} \times D^2$ is a simple closed curve $\gamma_b$ that bounds a $2$-disk $C_b$ containing $E_{1,1}, \ldots, E_{1,m}$.
Let $\sigma_j \in H_{n-1}(M_0^*)$ be the homology class obtained by taking $r(K_j \times [0,1])$ and capping the loop $\gamma_b$ with the $2$-disk $C_b$ for each $b \in \del K_j$.
Let $\sigma^\circ \in H_{n-1}(M_0^*)$ be the homology class obtained capping the boundary of ${r}(B^\circ \times [0,1])$ in the same way.
Then we see that the homology class $\overline{r}_*([B \times S^1]) \in H_{n-1}(M_0^*;\Z/\ell\Z)$ is equal to the sum $\sigma^\circ + \sigma_0 + \cdots + \sigma_k$.

The class ${\sigma}^\circ$ is trivial because the map $\overline{r}:B^\circ \times S^1 \rightarrow M_0^*$ extends to a map $B^\circ \times D^2 \rightarrow M_0^*$.
Indeed, if $b \in B^\circ$ lies outside of each $Z_{b_j}$, then the map $\overline{r}:\{b\} \times S^1 \rightarrow M_0^*$ is the constant loop at $b$.
Otherwise, if $b \in Z_{b_j} \cap B^\circ$, then the loop $\overline{r}:\{b\} \times S^1 \rightarrow M_0^*$ is either the constant loop at $b$, or a closed curve in a $D^2$-fiber of $Z_{b_j}$ which bounds $2$-disk in the same fiber.
\footnote{
    Note that this argument does not apply if $b \in K_j$, since then the loop $\overline{r}:\{b\} \times S^1 \rightarrow M_0^*$ might lie on a $D^2$-fiber of $Z_{b_j}$ that contains the center of a gluing disk, and thus does not bound a $2$-disk in the intersection of that fiber with $M_0^*$.
}

We also see that
\begin{equation*}
  \sigma_j = \sum_{g \in G_0} [\Sigma_{1,g_jg}].
\end{equation*}
for each $0 \leq j \leq k$.
Indeed, after we cap off ${r}(K_j \times S^1)$, we get a manifold with corners that cobounds an $n$-dimensional submanifold of $M_0^*$ with the spheres $\Sigma_{1,g_jg}$ for $g \in G_0$.
Note that $G_B = \{g_jg \mid g \in G_0, 0 \leq j \leq k\}$, as $G_B$ acts transitively on the points $b_0, \ldots, b_k$ with stabilizer $G_0$.\footnote{
    To prove this, observe first that $G_0 \leq G_B$ since $G_0$ fixes $B$ pointwise by Lemma \ref{lem:local-model-of-action}.
    Then one must check that $G_B \cdot b_0 = \{b_0, \ldots, b_k\}$.
    The inclusion ``$\subseteq$'' is immediate. 
    For the inclusion ``$\supseteq$,'' it suffices to check that for each $0 \leq j \leq k$, we have that $g_j \in G_B$, i.e.\ $g_j(B) = B$.
    To see this, note first that since $G_0$ fixes $B$ pointwise, we have that $g_j^{-1}G_0g_j \subseteq \mathrm{Stab}_G(b_0) = G_0$, and thus $g_j^{-1}G_0g_j = G_0$.
    This means that for any $g \in G_0$, $g_j^{-1}gg_j$ fixes $B$ pointwise.
    Then for any $b \in B$ and $g \in G_0$, we have that $g(g_j(b)) = g_j(g_j^{-1}gg_j(b)) = g_j(b)$, and so $g$ fixes $g_j(B_0)$ pointwise.
    This means that $g_j(B)$ is a component of the fixed set of $G_0$.
    Since $b_j \in g_j(B) \cap B$, it must be that $g_j(B) = B$.
}
Thus the proposition follows.

\subsection{Adaptations to the topological setting}
The construction of the element $\alpha_0$ carries over essentially verbatim to the topological setting, so long as we appeal to the topological version of Proposition \ref{prop:disk-slides-twists-and-inner-auts} described in Section \ref{subsec:disk-slides-adaptation}.

\section{Proof of Part (ii) of Theorem \ref{mainthm:ker-not-fg}}\label{sec:pf-part-ii}

In this section, we will complete the proof of part (ii) of Theorem \ref{mainthm:ker-not-fg}.
Fix a manifold $M$ and a finite group $G$ satisfying the hypotheses of Theorem \ref{mainthm:ker-not-fg}.
Recall that we let $\Mod_G(M,B) \leq \Mod(M,B)$ and $\Mod_G(M) \leq \Mod(M)$ denote the subgroups of mapping classes with an equivariant representative, and we have a forgetful map $\F_G:\Mod_G(M,B) \rightarrow \Mod_G(M)$.
In Section \ref{sec:obstruction-map}, we constructed a map
\begin{equation*}
    \varphi:\Ker(\F_G) \rightarrow H_{n-1}(\widetilde{M}, \widetilde{B};\Z/\ell\Z),
\end{equation*}
where $\pi:\widetilde{M} \rightarrow M$ is the cover constructed in Section \ref{sec:model-of-action} and $\widetilde{B}$ is a component of $\pi^{-1}(B)$.
In Section \ref{sec:nontrivial-obstruction}, we constructed an element $\alpha_0 \in \Ker(\F_G)$ for which $\varphi(\alpha_0)$ is nontrivial.

We prove Theorem \ref{mainthm:ker-not-fg} as follows.
In Section \ref{subsec:reduction-to-forgetful}, we show that it's enough to prove the group $\Ker(\F_G)$ is not finitely generated.
In Section \ref{subsec:obstruction-equivariant}, we show that the obstruction map $\varphi$ is equivariant with respect to the subgroup $L \leq \Mod_G(M,B)$ of mapping classes which lift to the cover $\widetilde{M}$.
In Section \ref{subsec:complete-pf-part-ii}, we show that the $L$-orbit of $\varphi(\alpha_0)$ contains an infinite linearly independent set of homology classes, which completes the proof.

\subsection{Reduction to the forgetful kernel}\label{subsec:reduction-to-forgetful}

First, we will show that it's enough to prove the kernel $\Ker(\F_G)$ is not finitely generated.

\begin{lemma}\label{lem:reduction-to-forgetful}
    If $\Ker(\F_G)$ is not finitely generated, then $\Ker(\P_G)$ is not finitely generated.
\end{lemma}
\begin{proof}
    Let $\Gamma_G(M,B) = \pi_0(\Diff(M,B)^G)$, where $\Diff(M,B)^G$ is the group of orientation-preserving diffeomorphisms of $M$ that restrict to an orientation-preserving diffeomorphism of $B$ and commute with $G$.
    Then we have the following commutative diagram:
    \begin{equation*}
      \begin{tikzcd}
        {\Gamma_G(M, B)} & {\Mod_G(M,B)} \\
        {\Gamma_G(M)} & {\Mod_G(M)}
        \arrow["\Psi", two heads, from=1-1, to=1-2]
        \arrow["\Phi"', hook, from=1-1, to=2-1]
        \arrow["{\F_G}", from=1-2, to=2-2]
        \arrow["{\P_G}", two heads, from=2-1, to=2-2]
      \end{tikzcd}
    \end{equation*}
    Here $\Phi$ is induced by the inclusion $\Diff(M,B)^G \hookrightarrow \Diff(M)^G$, and $\Psi$ is induced by the inclusion $\Diff(M,B)^G \hookrightarrow \Diff(M,B)$.
    The maps $\Psi$ and $\P_G$ are surjective by definition of $\Mod_G(M,B)$ and $\Mod_G(M)$.
    To see that the map $\Phi$ is injective, note first that if $f$ is an equivariant diffeomorphism that preserves $B$, then any equivariant isotopy of $f$ must also preserve $B$ setwise, since $B$ is a component of the fixed set of $g_0 \in G$.
    Thus, if $f$ represents an element of $\Ker(\Phi)$, then there is an equivariant isotopy $f \simeq \id_M$, and since this isotopy must preserve $B$, this implies that $f$ represents the identity in $\Gamma_G(M,B)$.
  
    Suppose now $\Ker(\F_G)$ is not finitely generated.
    Since $\Psi$ is surjective, we have a surjection $\Ker(\F_G \circ \Psi) \rightarrow \Ker(\F_G)$, and so $\Ker(\F_G \circ \Psi) = \Ker(\P_G \circ \Phi)$ is not finitely generated.
    Since $\Phi$ is injective, we have that $\Ker(\P_G \circ \Phi) \cong \Ker(\P_G) \cap \Im(\Phi)$.
    Thus it's enough to show that $\Im(\Phi)$ has finite index in $\Gamma_G(M)$.
    To see this, observe that $\Gamma_G(M)$ acts on the set of pairs $(C, \ep)$, where $C$ is a component of $\Fix(g_0)$ and $\ep$ is an orientation of $C$, and $\Im(\Phi)$ is a stabilizer for this action; then note that the set of such pairs is finite, since $\Fix(g_0)$ is a closed submanifold and $M$ is compact.
\end{proof}

\subsection{Equivariance of the obstruction map}\label{subsec:obstruction-equivariant}

Next, we will show that the obstruction map $\varphi$ satisfies a certain equivariance property.
Let $\Diff(M,B,b_0)$ denote the group of orientation-preserving diffeomorphisms of $M$ that restrict to an orientation-preserving diffeomorphism of $B$ and fix the basepoint $b_0 \in B$, and let $\Mod(M,B,b_0) = \pi_0(\Diff(M,B,b_0))$.
Let $\Mod_G(M,B,b_0) \leq \Mod(M,B,b_0)$ denote the subgroup of mapping classes with an equivariant representative.
Then, we let $L_* \leq \Mod_G(M,B,b_0)$ denote the subgroup of mapping classes that lift along the cover $\pi:\widetilde{M} \rightarrow M$, and we let $L \leq \Mod_G(M,B)$ denote the image of $L_*$.

Since every $\beta \in L_*$ fixes the base point $b_0 \in M$, any representative $f$ has a distinguished lift $\widetilde{f}$ which fixes the base point $\widetilde{b}_0 \in \widetilde{M}$.
Furthermore, since $f$ preserves $B$ setwise, it follows that $\widetilde{f}$ preserves $\widetilde{B}$ setwise.
Thus we get a well-defined map $L_* \rightarrow \Mod(\widetilde{M}, \widetilde{B})$.
Post-composing with the action on homology, we have an action of $L_*$ on $H_{n-1}(\widetilde{M}, \widetilde{B}; \Z/\ell\Z)$.

First, we can show that the basepoint is irrelevant to this action.

\begin{lemma}
    The action $L_* \rightarrow \Aut(H_{n-1}(\widetilde{M},\widetilde{B};\Z/\ell\Z))$ factors through $L$.
  \end{lemma}
  \begin{proof}
    Suppose $\beta \in \Ker(L_* \rightarrow L)$, and choose a representative $f$ of $\beta$.
    Then $f$ is isotopic to $\id_M$ through an isotopy preserving $B$ setwise.
    This lifts to an isotopy from $\widetilde{f}$ to a deck transformation.
    Since this lifted isotopy preserves $\widetilde{B}$, and the deck group of $\widetilde{M}$ acts freely on the components of $\pi^{-1}(B)$, we in fact have that $\widetilde{f}$ is isotopic to $\id_{\widetilde{M}}$ through an isotopy preserving $\widetilde{B}$ setwise.
    Thus $\widetilde{f}$ acts trivially on $H_{n-1}(\widetilde{M}, \widetilde{B};\Z/\ell\Z)$.
\end{proof}

Thus, we have a well-defined action of $L$ on $H_{n-1}(\widetilde{M},\widetilde{B};\Z/\ell\Z)$.
On the other hand, $L$ acts on the normal subgroup $\Ker(\F_G) \triangleleft \Mod_G(M,B)$ by conjugation.
We can show that $\varphi$ is equivariant with respect to these two actions. 

\begin{lemma}\label{lem:obstruction-map-equivariant}
    The map $\varphi:\Ker(\F_G) \rightarrow H_{n-1}(\widetilde{M}, \widetilde{B};\Z/\ell\Z)$ is $L$-equivariant.
\end{lemma}
\begin{proof}
    Take any $\alpha \in \Ker(\F_G)$ and $\beta \in L$.
    Our goal is to show that $\varphi(\beta \alpha \beta^{-1}) = \beta \cdot \varphi(\alpha)$.

    First, choose an isotopy $h:M \times [0,1] \rightarrow M$ such that $h_0 = \id$ and $h_1$ represents $\alpha$.
    Let $r:B \times [0,1] \rightarrow M$ be the map $r(b,t) = h_t(b)$, and let $\widetilde{r}:B \times [0,1] \rightarrow \widetilde{M}$ be the lift mapping $B \times \{0,1\}$ to $\widetilde{B}$.
    So, $\varphi(\alpha)$ is the image of the fundamental class under the natural map
    \begin{equation*}
        \widetilde{r}_*:H_{n-1}(B \times [0,1], B \times \{0,1\}; \Z/\ell\Z) \rightarrow H_{n-1}(\widetilde{M}, \widetilde{B}; \Z/\ell\Z).
    \end{equation*}

    Next, let $f$ be a representative of $\beta$; without loss of generality we may assume that $f$ fixes the basepoint $b_0 \in B$.
    Define the isotopy $h^f:M \times [0,1] \rightarrow M$ to be $f \circ h \circ (f^{-1} \times \id)$, so $h_t^f = fh_tf^{-1}$.
    Then the diffeomorphism $h_1^f$ represents $\beta\alpha\beta^{-1}$.
    Let $r^f:B \times [0,1] \rightarrow M$ be the map $r^f(b,t) = h_t^f(b)$, and let $\widetilde{r}^f:B \times [0,1] \rightarrow \widetilde{M}$ be the lift maping $B \times \{0,1\}$ to $\widetilde{B}$.
    Then $\varphi(\beta\alpha\beta^{-1})$ is the image of the fundamental class under the map
    \begin{equation*}
        \widetilde{r}^f_*:H_{n-1}(B \times [0,1], B \times \{0,1\}; \Z/\ell\Z) \rightarrow H_{n-1}(\widetilde{M}, \widetilde{B}; \Z/\ell\Z).
    \end{equation*}

    Now, let $\widetilde{f}$ be the lift of $f$ along the cover $\pi:\widetilde{M} \rightarrow M$ that fixes $\widetilde{b}_0$.
    We claim that $\widetilde{r}^f = \widetilde{f} \circ \widetilde{r} \circ (f^{-1} \times \id)$.
    Indeed, we see that 
    \begin{align*}
        \pi \circ \left(\widetilde{f} \circ \widetilde{r} \circ (f^{-1} \times \id)\right)
        &= f \circ \pi \circ \widetilde{r} \circ (f^{-1} \times \id) \\
        &= f \circ r \circ (f^{-1} \times \id) \\
        &= r^f,
    \end{align*}
    so $\widetilde{r}^f = \widetilde{f} \circ \widetilde{r} \circ (f^{-1} \times \id)$ by the uniqueness of lifts.

    So, we have the following commutative diagram:
    \begin{equation*}
        \begin{tikzcd}
            {H_{n-1}(B \times [0,1], B \times \{0,1\};\Z/\ell\Z)} & {H_{n-1}(B \times [0,1], B \times \{0,1\};\Z/\ell\Z)} \\
            & {H_{n-1}(\widetilde{M},\widetilde{B};\Z/\ell\Z)} \\
            & {H_{n-1}(\widetilde{M},\widetilde{B};\Z/\ell\Z)}
            \arrow["{(f^{-1} \times \id)_*}", from=1-1, to=1-2]
            \arrow["{\widetilde{r}^f_*}"', from=1-1, to=3-2]
            \arrow["{\widetilde{r}_*}", from=1-2, to=2-2]
            \arrow["{\widetilde{f}_*}", from=2-2, to=3-2]
        \end{tikzcd}
    \end{equation*}
    Since $f^{-1}$ is an orientation-preserving diffeomorphism of $B$, the map $f^{-1} \times \id$ acts by the identity on $H_{n-1}(B \times [0,1], B \times \{0,1\}; \Z/\ell\Z)$.
    Thus, the lemma follows.
\end{proof}

\subsection{Completing the proof}\label{subsec:complete-pf-part-ii}

By Proposition \ref{prop:nontrivial-obstruction}, we know that
\begin{equation*}
    \varphi(\alpha_0) = \sum_{g \in G_B} [\Sigma_{1,g,\id}],
\end{equation*}
where $G_B \leq G$ is the subgroup of $G$ preserving $B$ setwise.
By Lemma \ref{lem:obstruction-map-equivariant}, the map $\varphi$ is $L$-equivariant.
Thus, to prove Theorem \ref{mainthm:ker-not-fg}, it's enough to show that the $L$ orbit of $\sum_{g \in G_B} [\Sigma_{1,g,\id}]$ contains an infinite linearly independent subset of $H_{n-1}(\widetilde{M},\widetilde{B};\Z/\ell\Z)$.

Recall from Proposition \ref{prop:infinite-lin-ind-subset} that the set 
\begin{equation*}
        \mathcal{S} = \{[\Sigma_{1,g,u}] \in H_{n-1}(\widetilde{M};\Z/\ell\Z) \mid g \in G_{M_0}, u \in \Lambda_2\}
\end{equation*}
is infinite and linearly independent, where $G_{M_0}$ is the $G$-stabilizer of $M_0$ (so $G_B \leq G_{M_0}$) and $\Lambda_2 \leq \pi_1(M)$ is the subgroup constructed in Section \ref{subsec:infinite-cover}.
The set $\mathcal{S}$ embeds into $H_{n-1}(\widetilde{M},\widetilde{B};\Z/\ell\Z)$.
Thus, we can complete the proof of Theorem \ref{mainthm:ker-not-fg} with the following lemma.

\begin{lemma}\label{lem:orbit-is-not-fg}
    For each $u \in \Lambda_2$, there exists $\beta_u \in L$ such that $\beta_u([\Sigma_{1,\id,\id}]) = [\Sigma_{1,\id,u}]$, and for all $g \in G_{M_0}$, $\beta_u([\Sigma_{1,g,\id}]) = [\Sigma_{1,g,u'}]$ for some $u' \in \Lambda_2$.
\end{lemma}
\begin{proof}
    Throughout this proof, we will use the notation of Sections \ref{sec:disk-slides} and \ref{sec:model-of-action}.
    Recall that we define $\Lambda_2$ as $\pi_1(\overline{\calG}_2, M_0)$, and by Lemma \ref{lem:pi1-computation}, we can represent elements of $\Lambda_2$ by loops in $M$ based at $b_0$ which are disjoint from the gluing disks $D_{1,g}^\pm$ for all $g \in G$.
    Let $M'$ denote the manifold obtained from gluing $P_0$ and $P_2$ along the multidisks $\Delta_{2,g}^\pm$, and let $M_0'$ denote the component of $M'$ containing $M_0$.
    Then we can represent elements of $\Lambda_2$ with loops in $M_0'$.

    Now, fix $u \in \Lambda_2$.
    For $g \in G$, let $d_{1,g}$ denote the center of the gluing disk $D_{1,g}^+$ in $P_0$.
    By choosing an arc in $M_0$ from $b_0$ to $d_{1,\id}$, we get an identification $\pi_1(M_0',b_0) \cong \pi_1(M_0',d_{1,\id})$; let $\gamma_{\id}:S^1 \rightarrow M_0'$ be a loop based at $d_{1,\id}$ representing the element $u$.
    Then for $g \in G$, let $\gamma_g = g\gamma_{\id}$, so $\gamma_g$ is a loop in $M'$ based at $d_{1,g}$.

    Without loss of generality, we may assume that the loops $\gamma_g$ are all disjoint.
    Choose an element $\delta \in \pi_1(\Fr_{\Delta_1^+}(M'))$ which maps onto the tuple $(\gamma_g)_{g \in G}$ under the map $\pi_1(\Fr_{\Delta_1^+}(M')) \rightarrow \prod_{g \in G} \pi_1(M', d_{1,g})$.
    Then $\delta$ determines a disk slide $\calDS_{\Delta_1^+}(\delta)$ on $M'$; since the loops $\gamma_g$ are disjoint and permuted by $G$, we may choose an equivariant representative $f'$ of $\calDS_{\Delta_1^+}(\delta)$.
    Since $f'$ fixes the multidisks $\Delta_{1}^{\pm}$ pointwise, it extends to an equivariant diffeomorphism $f:M \rightarrow M$.
    Since $B$ has codimension $2$, we may assume that the loops $\gamma_g$ are disjoint from $B$, and hence we may assume that $f$ fixes $B$ pointwise.

    Let $\beta_u \in \Mod_G(M,B)$ be the element represented by $f$.
    Note that $\beta$ lies in the subgroup $L \leq \Mod_G(M,B)$; indeed, $\beta$ is liftable since $f$ acts trivially on the subgroup $\pi_1(M_0,b_0) \leq \pi_1(M,b_0)$ and hence preserves the normal closure $\la \la \pi_1(M_0,b_0) \ra \ra$.
    Let $\widetilde{f}:\widetilde{M} \rightarrow \widetilde{M}$ be the lift fixing $\widetilde{b}_0$.
    Then it follows by construction that $\widetilde{f}$ preserves the $\Lambda_2$-orbit of each $[\Sigma_{1,g,\id}]$, and that $\beta_u$ maps $[\Sigma_{1,\id,\id}]$ to $[\Sigma_{1,\id,u}]$.
\end{proof}

\subsection{Adaptations to the topological setting}
The arguments of this section carry over essentially verbatim to prove Theorem \ref{mainthm:top-version}.
Note that in the proof of Lemma \ref{lem:orbit-is-not-fg}, we use the topological version of Proposition \ref{prop:disk-slides-twists-and-inner-auts} described in Section \ref{subsec:disk-slides-adaptation}.

\section{Proof of Theorem \ref{mainthm:symmetric-aut-kernel}}\label{sec:pf-sym-aut-kernel}

In this section, we prove Theorem \ref{mainthm:symmetric-aut-kernel}. 
For this proof, it will be simpler to work with {homeomorphisms}, rather than diffeomorphisms, so we will appeal to Theorem \ref{mainthm:top-version}.

We prove Theorem \ref{mainthm:symmetric-aut-kernel} as follows.
In Section \ref{subsec:unlink-covers}, we recall the relationship between symmetric automorphisms and branched covers of the unlink in $S^3$.
In Section \ref{subsec:lifting-map}, we define a ``lifting map''
whose kernel is virtually isomorphic to (the topological version of) $\Ker(\P_G)$.
In Sections \ref{subsec:orbifold-pi1} and \ref{subsec:restriction-map}, we study the orbifold fundamental group associated to the branched cover.
We complete the proof in Section \ref{subsec:pf-sym-aut-kernel} by using the orbifold fundamental group to show that $\Ker(\calQ_{k,d})$ is virtually isomorphic to the kernel of the lifting map.


\subsection{Branched covers of the unlink and symmetric automorphisms}\label{subsec:unlink-covers}

Fix integers $k,d > 1$, and let $p_{k,d}:M_{k,d} \rightarrow S^3$ be the finite regular cover branched over the $k$-component unlink $C_k \subseteq S^3$ and with deck group $G_d = \Z/d\Z$.
For notational simplicity, we let $p = p_{k,d}$, $M = M_{k,d}$, and $G = G_d$.
We note that the manifold $M$ is homeomorphic to a connected sum of $S^1 \times S^2$'s, though we will not use this fact.
\footnote{
    This can be seen by viewing $S^3 \setminus C_k$ as a $k$-fold connected sum $(S^3 \setminus C_1)^{\# k}$ and constructing a model of $M$ as in Section \ref{sec:model-of-action}.
}

Let $\Homeo(S^3,C_k^{\pm})$ denote the group of orientation-preserving homeomorphisms of $S^3$ preserving $C_k$ setwise (not necessarily preserving the orientation of $C_k$), and let $\Mod(S^3,C_k^{\pm}) = \pi_0(\Homeo(S^3,C_k^{\pm}))$.
We have an isomorphism $\pi_1(S^3 \setminus C_k) \cong F_k$, where $F_k$ is the free group of rank $k$.
We fix a basis $\{x_1, \ldots, x_k\}$ of $F_k$ where each $x_i$ corresponds to a loop that winds once around the $i$th component of $C_k$.
As mentioned in Section \ref{subsec:appliction-to-symmetric-auts}, a theorem of Goldsmith \cite{goldsmith}, building off the work of Dahm \cite{dahm}, says that the action of $\Mod(S^3,C_k^{\pm})$ on $\pi_1(S^3 \setminus C_k)$ induces an isomorphism
\begin{equation*}
    \Mod(S^3, C_k^{\pm}) \cong \SymOut(F_k),
\end{equation*}
where $\SymOut(F_k)$ is the group of outer automorphisms sending each $x_i$ to a conjugate of some $x_j^{\pm 1}$ (Goldsmith actually proves an analogous result about the $\R^3$, but the above isomorphism can easily be deduced, see e.g. \cite[Prop~5.1]{bh-3manifolds}).
We let $\Mod(S^3, C_k)$ denote the finite index subgroup of homeomorphisms preserving the orientation of $C_k$, and let $\SymOut^+(F_k)$ denote its image in $\SymOut(F_k)$.
Then $\SymOut^+(F_k)$ is the group of outer automorphisms sending each $x_i$ to a conjugate of some $x_j$.

\subsection{The lifting map}\label{subsec:lifting-map}
Let $M^\circ = M \setminus p^{-1}(C)$.
Then the unbranched cover $M^\circ \rightarrow S^3 \setminus C_k$ is classified by the natural map $F_k \rightarrow \Z/d\Z$ sending each $x_i$ to $1$.
Since precomposition by an element of $\SymOut^+(F_k)$ preserves this homomorphism, it follows that every element of $\Mod(S^3,C_k)$ lifts along the cover $p$ to an equivariant homeomorphism of $M$ (see \cite[\S 2.4]{bh-3manifolds}).
Let $\Mod(M) = \pi_0(\Homeo(M))$, and let $\Mod_G(M)$ be the subgroup of mapping classes represented by an equivariant homeomorphism.
Since any lift is well-defined up to a deck transformation, we get a \emph{lifting map}
\begin{equation*}
    \L_p:\Mod(S^3,C_k) \rightarrow \Mod_G(M)/G.
\end{equation*}
which sends a mapping class to a lift.

Let $\Gamma_G(M) = \pi_0(\Homeo(M)^G)$.
Observe that the kernel $\Ker(\L_p)$ is virtually isomorphic to the kernel of the map $\P_G:\Gamma_G(M) \rightarrow \Mod_G(M)$.
Indeed, since any equivariant (topological) isotopy of $M$ descends to an isotopy of $S^3$ preserving $C_k$ setwise, we get a map $\Gamma_G(M) \rightarrow \Mod(S^3,C_k^{\pm})$ with a finite kernel and a finite-index image.
This induces a virtual isomorphism between $\Ker(\L_p)$ and $\Ker(\P_G)$.
By Theorem \ref{mainthm:top-version}, we conclude that $\Ker(\L_p)$ is not finitely generated for $k \geq 3$.

\subsection{The orbifold fundamental group}\label{subsec:orbifold-pi1}
The branched cover $p$ induces an orbifold structure $\mathcal{O}_p$ on $S^3$, where each point of $C_k$ is an orbifold point of order $d$.
Then we have an isomorphism
\begin{equation*}
    \pi_1^{\orb}(\calO_p) \cong \underbrace{\Z/d\Z * \cdots * \Z/d\Z}_{\text{$k$ times}} \eqqcolon H_{k,d}.
\end{equation*}
Indeed, the orbifold fundamental group is the quotient of $\pi_1(S^3 \setminus C_k)$ obtained by adding the relations $x_i^d = 1$ for each $i$ (see e.g. \cite{haefliger-du}).
We let $z_i$ denote the image of $x_i$ in $\pi_1^{\orb}(\calO_p)$.
The branched cover $p$ gives us the following commutative diagram with exact rows:
\begin{equation}\label{eqn:orbifold-pi1}
    \begin{tikzcd}
        1 & {\pi_1(M^\circ)} & {\pi_1(S^3 \setminus C_k)} & {G } & 1 \\
        1 & {\pi_1(M)} & {\pi_1^{\orb}(\calO_p)} & G & 1
        \arrow[from=1-1, to=1-2]
        \arrow[from=1-2, to=1-3]
        \arrow[two heads, from=1-2, to=2-2]
        \arrow[from=1-3, to=1-4]
        \arrow[two heads, from=1-3, to=2-3]
        \arrow[from=1-4, to=1-5]
        \arrow[from=1-4, to=2-4, equals]
        \arrow[from=2-1, to=2-2]
        \arrow[from=2-2, to=2-3]
        \arrow[from=2-3, to=2-4]
        \arrow[from=2-4, to=2-5]
    \end{tikzcd}
\end{equation}
where the map $\pi_1^{\orb}(\calO_p) \rightarrow G = \Z/d\Z$ maps each $z_i$ to $1$.

Since $\Mod(S^3,C_k)$ preserves the normal closure of the elements $x_1^d, \ldots, x_k^d$, the action of $\Mod(S^3,C_k)$ on $\pi_1(S^3 \setminus C_k)$ descends to an action of $\Mod(S^3, C_k)$ on $\pi_1^{\orb}(\calO_p)$ by outer automorphisms.
Let $\SymOut(H_{k,d})$ denote the group of outer automorphisms sending each $z_i$ to a conjugate of some $z_j^{\pm}$, and let $\SymOut^+(H_{k,d})$ denote the finite index subgroup sending each $z_i$ to a conjugate of some $z_j$.
The projection $F_k \rightarrow H_{k,d}$ induces a natural map
\begin{equation*}
    \calQ_{k,d}:\SymOut(F_k) \rightarrow \SymOut(H_{k,d}).
\end{equation*}
It follows that the action of $\Mod^+(S^3,C_k)$ on $\pi_1^{\orb}(\calO_p)$ factors through the restriction
\begin{equation*}
    \calQ_{k,d}^+:\SymOut^+(F_k) \rightarrow \SymOut^+(H_{k,d}).
\end{equation*}

\subsection{The restriction map}\label{subsec:restriction-map}

Let $\Aut(\pi_1^{\orb}(\calO_p),\pi_1(M))$ denote the group of automorphisms of $\pi_1^{\orb}(\calO_p)$ that preserve the subgroup $\pi_1(M)$.
Since $\pi_1(M)$ is a normal subgroup of $\pi_1^{\orb}(\calO_p)$, we know that $\Inn(\pi_1^{\orb}(\calO_p)) \leq \Aut(\pi_1^{\orb}(\calO_p),\pi_1(M))$, and so we define
\begin{equation*}
    \Out(\pi_1^{\orb}(\calO_p),\pi_1(M)) \coloneqq \Aut(\pi_1^{\orb}(\calO_p),\pi_1(M))/\Inn(\pi_1^{\orb}(\calO_p)).
\end{equation*}
Let $G_*$ denote the image of $G$ in $\Out(\pi_1(M))$, and let $\Out_{G_*}(\pi_1(M))$ denote its normalizer.
Then we claim there is a map
\begin{equation*}
    \overline{r}:\Out(\pi_1^{\orb}(\calO_p),\pi_1(M)) \rightarrow \Out_{G_*}(\pi_1(M))/G_*.
\end{equation*}
To construct the map $\overline{r}$, we start with the restriction map
\begin{equation*}
    r:\Aut(\pi_1^{\orb}(\calO_p),\pi_1(M)) \rightarrow \Aut(\pi_1(M)),
\end{equation*}
and post-compose to obtain a map
\begin{equation*}
    r':\Aut(\pi_1^{\orb}(\calO_p),\pi_1(M)) \rightarrow \Out(\pi_1(M)).
\end{equation*}
From the short exact sequence $\pi_1(M) \hookrightarrow \pi_1^{\orb}(\calO_p) \twoheadrightarrow G$, we see that the map $r'$ takes $\Inn(\pi_1^{\orb}(\calO_p))$ to the group $G_*$.
Since $\Inn(\pi_1^{\orb}(\calO_p))$ is normal in $\Aut(\pi_1^{\orb}(\calO_p),\pi_1(M))$, it follows that $G_*$ is normal in the image of $r'$.
Thus the map $r'$ descends to the map $\overline{r}$.

Now, applying the isomorphism $\pi_1^{\orb}(\calO_p) \cong H_{k,d}$, we can view $\SymOut^+(H_{k,d})$ as a subgroup of $\Out(\pi_1^{\orb}(\calO_p),\pi_1(M))$.
Then $\overline{r}$ restricts a map
\begin{equation*}
    \overline{R}:\SymOut^+(H_{k,d}) \rightarrow \Out_{G_*}(\pi_1(M))/G_*.
\end{equation*}

\begin{lemma}\label{lem:restriction-map-kernel}
    The map $\overline{R}:\SymOut^+(H_{k,d}) \rightarrow \Out_{G_*}(\pi_1(M))/G_*$ has a finite kernel.
\end{lemma}
\begin{proof}
    Let $\PSymOut(H_{k,d})$ be the finite index subgroup of $\SymOut(H_{k,d})$ consisting of outer automorphisms that send each $z_i$ to a conjugate of itself.
    Then it's enough to show that $\overline{R}$ is injective on $\PSymOut(H_{k,d})$.

    Suppose $f$ represents an element of $\PSymOut(H_{k,d})$ and $\overline{R}([f])$ is trivial.
    Our goal is to show that $f$ is an inner automorphism of $H_{k,d}$.
    Since $\overline{R}([f])$ is trivial, this means there exists $\iota \in \Inn(H_{k,d})$ such that $\iota \circ f$ restricts to the identity map on $\pi_1(M)$.
    It's enough to show that  $\iota \circ f$ is inner.
    Thus, we assume without loss of generality that $f$ itself acts trivially on $\pi_1(M)$.
    
    We will view elements of $H_{k,d}$ as freely reduced words on the elements $z_1, \ldots, z_k$.
    For each $1 \leq i \leq k$, we have that $f(z_i) = w_iz_iw_i^{-1}$ for some reduced word $w_i$.
    Without loss of generality, we may assume that $w_i$ does not end with a power of $z_i$.

    First, we claim that if $w_i$ is trivial for some $i$, then $w_i$ must be trivial for all $1 \leq i \leq k$, in which case we are done.
    To see this, suppose $w_i$ is trivial and let $1 \leq j \leq k$ with $j \neq i$.
    Then the element $z_iz_j^{-1}$ lies in the subgroup $\pi_1(M)$, and so since $f$ is the identity on $\pi_1(M)$, we have that
    \begin{equation*}
        z_iz_j^{-1} = f(z_iz_j^{-1}) = z_iw_jz_j^{-1}w_j^{-1}.
    \end{equation*}
    This is possible only if $w_j$ is trivial.

    Otherwise, assume that $w_i$ is nontrivial for all $i$.
    Let $i$ and $j$ be arbitrary.
    Then we have that
    \begin{equation*}
        z_iz_j^{-1} = f(z_iz_j^{-1})
        = w_iz_iw_i^{-1}w_jz_jw_j^{-1}.
    \end{equation*}
    Given a reduced word $w$ on the elements $z_1, \ldots, z_k$, we say that a reduced word $w'$ is a \emph{left subword} (resp.\ \emph{right subword}) of $w$ if there exists a reduced word $v$ such that $w = w'v$ and the last letter of $w'$ is not a power of the first letter of $v$ (resp.\ $w = vw'$ and the last letter of $v$ is the not a power of the first letter of $w'$).
    We say the left (resp.\ right) subword is \emph{proper} if $v$ is nontrivial.
    Then the above equality implies that
    \begin{equation*}
        w_i^{-1}w_j = w
    \end{equation*}
    where $w$ is either a proper left subword of $w_i^{-1}$ or a proper right subword of $w_j$.
    Now, the equality
    \begin{equation*}
        z_iz_j^{-1} = w_iz_iwz_j^{-1}w_j^{-1}
    \end{equation*}
    implies that $w$ is nontrivial, and 
    either 
    \begin{equation*}
        w_iz_iwz_j^{-1}w_j^{-1}
        = w_i'z_j^{-1}w_j^{-1}
    \end{equation*}
    where $w_i'$ is a left subword of $w_i$, or 
    \begin{equation*}
        w_iz_iwz_j^{-1}w_j^{-1}
        = w_iz_i(w_j^{-1})'
    \end{equation*}
    where $(w_j^{-1})'$ is a right subword of $w_j^{-1}$.
    Then we have that
    \begin{equation*}
        z_iz_j^{-1} = w_i'z_j^{-1}w_j^{-1} \text{ or } z_iz_j^{-1} = w_iz_i(w_j^{-1})',
    \end{equation*}
    and either case implies that $w_i$ starts with $z_i$ or $w_j$ starts with $z_j$.
    But then from the equality 
    \begin{equation*}
        w_i^{-1}w_j = w,
    \end{equation*}
    we see that if $w_i$ starts with $z_i$, then $w_j$ must start with $z_i$, and if $w_j$ starts with $z_j$, then $w_i$ must start with $z_j$.
    Thus, we conclude that $w_i$ and $w_j$ must start with the same letter.
    Since $i$ and $j$ were arbitrary, we conclude that the words $w_1, \ldots, w_k$ all must start with the same letter $z$.
    Then if we let $f' = \Inn(z) \circ f$, we have that $f'(z_i) = v_i z_i v_i^{-1}$ where $v_i$ is a reduced word of shorter length than $w_i$.
    By repeating the argument inductively, we obtain a sequence of inner automorphisms $\iota_1, \ldots, \iota_m$ such that $\iota_1 \circ \cdots \iota_m \circ f$ is trivial, and thus $f$ is inner as desired.
\end{proof}

\subsection{Proof of Theorem \ref{mainthm:symmetric-aut-kernel}}\label{subsec:pf-sym-aut-kernel}

Now, we can prove Theorem \ref{mainthm:symmetric-aut-kernel}.
First, note that if $\Ker(\calQ_{k,d})$ is not finitely generated, then $\Ker(\widehat{\calQ}_{k,d})$ must also be not finitely generated.
Indeed, we have the following commutative diagram:
\begin{equation*}
    \begin{tikzcd}
        {\SymAut(F_k)} & {\SymAut(H_{k,d})} \\
        {\SymOut(F_k)} & {\SymOut(H_{k,d})}
        \arrow["{\widehat{\calQ}_{k,d}}", from=1-1, to=1-2]
        \arrow["\Phi"', two heads, from=1-1, to=2-1]
        \arrow["\Psi", two heads, from=1-2, to=2-2]
        \arrow["{\calQ_{k,d}}"', from=2-1, to=2-2]
    \end{tikzcd}
\end{equation*}
Here $\Phi$ and $\Psi$ are both the quotient map by the group of inner automorphisms.
If $\Ker(\calQ_{k,d})$ is not finitely generated, then since $\Ker(\calQ_{k,d} \circ \Phi)$ surjects onto $\Ker(\calQ_{k,d})$, we conclude that $\Ker(\calQ_{k,d} \circ \Phi) = \Ker(\Psi \circ \widehat{\calQ}_{k,d})$ is not finitely generated.
Now, since the map $F_k \rightarrow H_{k,d}$ is surjective, it follows that $\widehat{\calQ}_{k,d}$ restricts to a surjection $\Inn(F_k) \rightarrow \Inn(H_{k,d}) = \Ker(\Psi)$.
Thus we have a short exact sequence
\begin{equation*}
    1 \rightarrow \Ker(\widehat{\calQ}_{k,d}) \rightarrow \Ker(\Psi \circ \widehat{\calQ}_{k,d}) \rightarrow \Ker(\Psi) \rightarrow 1.
\end{equation*}
Since $\Ker(\Psi) = \Inn(H_{k,d})$ is finitely generated, it must be that $\Ker(\widehat{\calQ}_{k,d})$ is not finitely generated.

Thus, it remains to prove that the group $\Ker(\calQ_{k,d})$ is not finitely generated for $k \geq 3$.
As discussed in Section \ref{subsec:lifting-map}, Theorem \ref{mainthm:ker-not-fg} implies that $\Ker(\L_p)$ is not finitely generated for $k \geq 3$.
Thus Theorem \ref{mainthm:symmetric-aut-kernel} is a consequence of the following proposition.

\begin{proposition}
    The kernels of the maps
    \begin{equation*}
        \L_p:\Mod(S^3,C_k) \rightarrow \Mod_G(M)/G
    \end{equation*}
    and
    \begin{equation*}
        \calQ_{k,d}:\SymOut(F_k) \rightarrow \SymOut(H_{k,d})
    \end{equation*}
    are virtually isomorphic.
\end{proposition}
\begin{proof}
    Let $\overline{R}:\SymOut^+(H_{k,d}) \rightarrow \Out_{G_*}(\pi_1(M))/G_*$ be the restriction map from Section \ref{subsec:restriction-map}.
    Observe that if $f$ is a homeomorphism of $(S^3, C_k)$ representing an element of $\Mod(S^3,C_k)$, then the from the diagram (\ref{eqn:orbifold-pi1}), the action of a lift $\widetilde{f}$ on $\pi_1(M)$ is obtained by restricting the action of $f$ on $\pi_1^{\orb}(\calO_p)$.
    Therefore we have the following commutative diagram:
    \begin{equation*}
        \begin{tikzcd}
            {\Mod(S^3,C_k)} & {\SymOut^+(F_k)} \\
            & {\SymOut^+(H_{k,d})} \\
            {\Mod_G(M)/G} & {\Out_{G_*}(\pi_1(M))/G_*}
            \arrow["\cong", from=1-1, to=1-2]
            \arrow["{\L_p}", from=1-1, to=3-1]
            \arrow["{\calQ_{k,d}^+}", from=1-2, to=2-2]
            \arrow["{\overline{R}}", from=2-2, to=3-2]
            \arrow["{\Phi}",from=3-1, to=3-2]
        \end{tikzcd}
    \end{equation*}
    Here $\Phi$ is induced by the action of $\Mod(M)$ on $\pi_1(M)$ by outer automorphisms.
    To prove the proposition, it's enough to show that the maps $\Phi$ and $\overline{R}$ have finite kernels.
    We know that $\overline{R}$ has a finite kernel by Lemma \ref{lem:restriction-map-kernel}. 
    We also know that $\Phi$ has a finite kernel since the kernel of the map $\Mod(M) \rightarrow \Out(\pi_1(M))$ is the subgroup generated by sphere twists (see e.g.\ \cite[Prop~2.1]{hatcher-wahl}), and this subgroup is finite \cite[Prop~1.2]{mccullough}.
\end{proof}

\printbibliography

\end{document}